\documentclass[10pt,letterpaper]{amsart} 

\usepackage{color}
\usepackage{multicol}
\usepackage[letterpaper, margin=1in]{geometry}
\usepackage[overload]{empheq}

\usepackage{xcolor}

\usepackage{listings}
\usepackage{amsmath, amsthm, amssymb, wasysym, verbatim, bbm, graphics}
\usepackage{enumerate}
\usepackage{graphicx}
\usepackage{url}
\usepackage{mathtools}
\usepackage[colorlinks,linkcolor=blue]{hyperref}
\usepackage{romannum}
\usepackage{xfrac}
\usepackage{float}
\usepackage{amsmath,bm}
\usepackage[title]{appendix}
\usepackage{todonotes}
\usepackage{subcaption}
\usepackage{caption}
\usepackage{enumitem}
\usepackage{algorithm}
\usepackage{algpseudocode}
\usepackage{amsfonts}

\newtheorem{lemma}{Lemma}[section]

\newtheorem*{maintheorem*}{Main Theorem}
\theoremstyle{definition}{}

\newcommand{\Plambda}{P^\lambda}
\newcommand{\Plambdaone}{P^{\lambda_1}}
\newcommand{\Plambdatwo}{P^{\lambda_2}}
\newcommand{\ulambdaone}{u^{\lambda_1}}
\newcommand{\ulambdatwo}{u^{\lambda_2}}
\newcommand{\ylambdaone}{y^{\lambda_1}}
\newcommand{\ylambdatwo}{y^{\lambda_2}}
\newcommand{\Alambdaonelambdatwo}{ A^{\lambda_1,\lambda_2}_\omega}
\newcommand{\Alambdaonelambdatwok}{ A^{\lambda_1,\lambda_2}_{\omega,k}}
\newcommand{\ulambdaonelambdatwoomega}{ u^{\lambda_1,\lambda_2}_\omega}
\newcommand{\ylambdaonelambdatwoomega}{ y^{\lambda_1,\lambda_2}_\omega}
\newcommand{\ulambdaonelambdatwoomegak}{ u^{\lambda_1,\lambda_2}_{\omega,k}}
\newcommand{\ylambdaonelambdatwoomegak}{ y^{\lambda_1,\lambda_2}_{\omega,k}}
\newcommand{\Hlambda}{H^\lambda}
\newcommand{\unn}{u_{NN}}
\newcommand{\ynn}{y_{NN}}
\newcommand{\fnn}{f_{NN}}

\newcommand{\usnn}{u^s_{NN}}
\newcommand{\ysnn}{y^s_{NN}}
\newcommand{\fsnn}{f^s_{NN}}

\newcommand{\udnn}{u^d_{NN}}
\newcommand{\ydnn}{y^d_{NN}}
\newcommand{\fdnn}{f^d_{NN}}

\newtheorem{thm}[lemma]{Theorem}

\newtheorem{prop}[lemma]{Proposition}
\newtheorem{rem}[lemma]{Remark}
\newtheorem{lem}[lemma]{Lemma}

\allowdisplaybreaks

\numberwithin{equation}{section}

\title[Penalty Adversarial Network]{Penalty Adversarial Network (PAN): A neural network-based method to solve PDE-constrained optimal control problems}
\date{\today}

\author[S. Ma]{Shilin Ma}
\address[Shilin Ma]{\newline Acadian Asset Management LLC\footnote{The views expressed herein are those of the author and do not necessarily reflect the views of Acadian Asset Management LLC. The views should not be considered investment advice and do not constitute or form part of any offer to issue or sell, or any solicitation of any offer to subscribe or to purchase, shares, units or other interests in any particular investments.}, Boston, MA, USA.}
\email[]{shilinm@alumni.cmu.edu}

\author[Y. Yue]{Yukun Yue}
\address[Yukun Yue]{ \newline Department of Mathematical Sciences \newline University of Wisconsin Madison, Madison, Wisconsin, USA.}
\email[]{yyue@math.wisc.edu}

\begin{document}
 \pagenumbering{arabic}
\maketitle

\begin{abstract}
In this work, we introduce a novel strategy for tackling constrained optimization problems through a modified penalty method. Conventional penalty methods convert constrained problems into unconstrained ones by incorporating constraints into the loss function via a penalty term. However, selecting an optimal penalty parameter remains challenging; an improper choice, whether excessively high or low, can significantly impede the discovery of the true solution. This challenge is particularly evident when training neural networks for constrained optimization, where tuning parameters can become an extensive and laborious task. To overcome these issues, we propose an adversarial approach that redefines the conventional penalty method by simultaneously considering two competing penalty problems—a technique we term the penalty adversarial problem. Within linear settings, our method not only ensures the fulfillment of constraints but also guarantees solvability, leading to more precise solutions compared to traditional approaches. We further reveal that our method effectively performs an automatic adjustment of penalty parameters by leveraging the relationship between the objective and loss functions, thereby obviating the need for manual parameter tuning. Additionally, we extend this adversarial framework to develop a neural network-based solution for optimal control problems governed by linear or nonlinear partial differential equations. We demonstrate the efficacy of this innovative approach through a series of numerical examples.

\end{abstract}
\section{Introduction}

Optimal control problems are fundamental in various scientific and engineering disciplines. These problems involve finding a control state \(y\) that determines the desired state \(u\) through governing physical constraints, aiming to minimize or maximize a given performance criterion, typically expressed as an objective functional \cite{kirk2004optimal,lewis2012optimal}. In many practical scenarios, the system dynamics are governed by partial differential equations (PDEs), leading to PDE-constrained optimal control problems \cite{lions1971optimal,russell1978controllability}. These problems have gained significant attention due to their critical applications in fields such as aerospace engineering \cite{longuski2014optimal}, environmental marine sciences \cite{strazzullo2018model}, medical treatment planning for radiation therapy \cite{gibou2005partial}, heat transfer \cite{wang2008null}, fluid dynamics \cite{ghattas1997optimal,gunzburger1991analysis,mathew2007optimal}, liquid crystals \cite{surowiec2023optimal}, and wave propagation \cite{benamou1997domain}. Achieving optimal performance while adhering to physical laws and constraints is crucial in these applications.

Mathematically, we can formulate the problem as 
\begin{equation}\label{eq:constrained_problem}
    \min_{u \in \mathrm{U}, y \in \mathrm{Y}} J(u,y), \quad \text{subject to } F(u,y) = 0.
\end{equation}
Here, \( J(u,y) \) represents the performance criterion to be minimized, often referred to as the objective functional \cite{rosseel2012optimal}. The term \( F(u,y) \) contains the constraints that \( u \) and \( y \) must satisfy, including the differential operators in the form of PDEs and the boundary or initial conditions for the PDEs. We denote \( \mathrm{U} \) and \( \mathrm{Y} \) as the appropriate spaces in which \( u \) and \( y \) reside, respectively.

The highly nonlinear nature and multi-scale structure of PDE-constrained optimal control problems \cite{chen2017global,einkemmer2024suppressing,ferrari2016distributed} necessitate using complex numerical methods. Over the years, various approaches have been developed to create robust and accurate numerical algorithms and tools to solve these problems. Among the prevalent methods, adjoint-based techniques are particularly notable for their effectiveness in gradient computation, which is crucial for iterative optimization algorithms \cite{buskens2000sqp,jameson1988aerodynamic}. These methods are often combined with traditional numerical techniques, such as the finite difference method or finite element methods, to handle spatial and temporal discretizations, allowing for the management of complex geometries and boundary conditions \cite{antil2014optimal,brenner2020multigrid,brenner2016c,rees2010optimal}.

Recently, with the rapid development of neural networks, machine learning-based numerical methods have been extensively developed to solve PDEs \cite{Beck2023overview,chen2021solving,ding2022overparameterization,goswami2023physics,goswami2022physics,han2018solving,karniadakis2021physics,kim2022fast,lu2021learning,lu2021deepxde,raissi2019physics,sirignano2018dgm,wang2021learning,zang2020weak}. As an important extension of this work, considerable research efforts have focused on applying these methods to solve PDE-constrained optimal control problems. Several successful examples have emerged in this area \cite{Barry2022physics,demo2023extended,garcia2023control, hwang2022solving,lu2021physics,lye2021iterative,mowlavi2023optimal,yin2024aonn}. Among them, the most prevalent approaches can be classified into three main categories:

\begin{enumerate}[label=\textbullet]
    \item Training surrogate models to obtain solvers for PDEs, then using the trained solver to map inputs to solutions of the PDE, enforcing the PDE constraint while minimizing an objective cost functional \cite{hwang2022solving,lye2021iterative};
    \item Using the Lagrangian approach to reformulate the constrained optimization problem and solve systems associated with the Karush–Kuhn–Tucker (KKT) conditions \cite{karush2013minima,kuhn2013nonlinear} via classical neural network-based methods \cite{Barry2022physics,demo2023extended};
    \item Adding the cost functional to the standard loss induced by the PDEs and minimizing the total loss simultaneously \cite{lu2021physics, mowlavi2023optimal}.
\end{enumerate}

A common challenge these approaches encounter is enforcing the constraints during optimization. The most direct approach is to treat the constraints as penalty terms and add them to the original loss function minimized by the neural network, transforming a constrained problem into an unconstrained one. This transformation makes it easier to apply various iterative optimization tools to solve the problem \cite{boyd2004convex}, as implemented in the third approach above. Though not directly employing this method when enforcing PDE constraint, the other two approaches implicitly address the same challenge while solving the PDE by ensuring that initial or boundary conditions are satisfied \cite{han2018solving, lu2021learning}. For simplicity, we will explain our main idea concerning the last approach as an example. At the same time, the same analysis can be applied to the other two approaches when penalty terms are introduced.

Specifically, instead of solving problem \eqref{eq:constrained_problem}, we consider a penalty problem which can be formulated as:
\begin{equation}
    \label{eq:unconstrained_penalty_problem}
    \min_{u\in\mathrm{U},y\in\mathrm{Y}} \Plambda(u,y),
\end{equation}
where $\Plambda$ is defined as
\begin{equation}
    \label{eq:P_lambda_def}
    \Plambda(u,y) = J(u,y) + \frac{\lambda}{2}\,\|F(u,y)\|^2.
\end{equation}
Here, $\lambda > 0$ is a tunable penalty parameter, and $\|\cdot\|$ denotes the standard $L^2$ norm. It is important to note that problem \eqref{eq:unconstrained_penalty_problem} is not equivalent to the constrained problem \eqref{eq:constrained_problem}. However, as $\lambda$ approaches infinity, the solution of this unconstrained problem converges to the solution of the constrained one \cite{nocedal1999numerical} (We will provide more details on this in the linear case in Section \ref{sec:linear_analysis}). Conversely, when the penalty parameter becomes too large, the problem can become ill-posed and difficult to solve \cite{bertsekas2014constrained}. If a neural network is used to solve this optimization problem, achieving convergence during training can be challenging. On the other hand, if the penalty parameter is too small, the PDE constraints will not be adequately satisfied, resulting in a solution that deviates significantly from the desired solution of the original constrained problem. 

To address this challenge, \cite{mowlavi2023optimal} proposes a two-step line-search approach to determine the optimal penalty parameter that minimizes the cost function while ensuring the PDE constraints are satisfied within an acceptable tolerance. In contrast, \cite{lu2021physics} suggests dynamically adjusting the penalty parameters based on a problem-dependent update rule \cite{bertsekas2014constrained}. These approaches commonly face issues such as the extensive effort required for parameter tuning and the absence of general guidelines for selecting an appropriate penalty parameter for various problems. This complicates the implementation of existing methods for solving optimal control-related problems.

In this paper, we propose a novel framework based on penalty methods and present it using a neural network structure that does not require an artificial selection process for the varying value of the penalty parameter. Instead, we construct two neural networks to minimize $\Plambda(u,y)$ with different fixed values of $\lambda$ and train them simultaneously. One network competes with the other during training to ensure ease of training and convergence while maintaining the constraints within an acceptable tolerance. This concept is inspired by the well-known generative adversarial network (GAN) \cite{creswell2018generative,karras2019style,goodfellow2020generative}, which has extensive applications, for example, in image translation \cite{isola2017image,zhang2017stackgan}, video generation \cite{Vondrick2016GeneratingVW}, and speech synthesis \cite{kong2020hifi}. The adversarial network structure incorporates a generator and a discriminator, with the discriminator providing feedback to the generator to help it produce more realistic information to deceive the discriminator. This idea has also benefited research on numerical methods based on adversarial network structures for solving PDEs. We refer readers to \cite{kadeethum2021framework,thanasutives2021adversarial,wu2020enforcing,yang2020physics,yang2019adversarial,zang2020weak} and the references therein.

In particular, we create a solver network (corresponding to the generator in GAN) and a discriminator network, and we choose two real numbers $\lambda_1, \lambda_2 > 0$ with $\lambda_1$ being much larger than $\lambda_2$, and $\lambda_2$ being relatively small. The discriminator network is set to minimize the objective functional $\Plambdatwo(u,y)$ as defined in \eqref{eq:P_lambda_def}. In practice, the smallness of $\lambda_2$ ensures the convergence of the discriminator network as long as the optimal control problem \eqref{eq:constrained_problem} is well-posed. Conversely, the solver network aims to minimize the sum of $\Plambdaone(u,y)$ and an additional cost based on feedback from the discriminator network. The exact functional form of this extra cost will be detailed in Section \ref{sec:PAN}. 

With feedback from the discriminator network, the solver network can automatically adjust its focus during training without manual tuning. If the solver network emphasizes reducing the objective functional at the expense of not maintaining the PDE constraint, the weight of the PDE constraint will increase accordingly. Conversely, if the solver network focuses too much on satisfying the PDE constraint but fails to reduce the objective cost functional, it will adjust its weights on the objective functionals to correct this imbalance.

To this end, the proposed framework utilizes the strengths of both traditional penalty methods and the adversarial training paradigm, offering a novel solution to the challenges inherent in PDE-constrained optimal control problems. Numerical examples, presented in Section \ref{sec:numerics}, demonstrate that our approach ensures robust convergence and effective enforcement of PDE constraints without requiring extensive parameter tuning. This dual-network strategy not only simplifies the training process but also enhances the overall stability and performance of the optimization.

The rest of this paper is structured as follows: In Section \ref{sec:problem_setup}, we describe various PDE-constrained optimal control problems in general forms that will be the focus of this paper. Section \ref{sec:linear_analysis} serves as a motivation for our methodology, where a linear problem is discussed to illustrate the effectiveness of the proposed penalty adversarial framework in solving problems with penalty formulations. We demonstrate that in the linear case, the solution to the adversarial problem better adheres to the constraints than simply solving the problem with a small penalty parameter under certain conditions. Following this analysis, Section \ref{sec:PAN} provides a detailed construction of a neural network-based method utilizing this concept, culminating in a practical algorithm. Finally, in Section \ref{sec:numerics}, we conduct numerical experiments on both linear and nonlinear problems in 1D and 2D to validate the effectiveness of the proposed approach.

\section{Problem Setup}\label{sec:problem_setup}

This section formally outlines the various types of optimal control problems, including distributed, boundary, and initial value control problems. We consider an open bounded physical domain $\Omega \subset \mathbb{R}^d$, where $d$ is a positive integer denoting the spatial dimension, and a time span $[0, T]$. The problem is governed by the following system:

\begin{subequations}\label{eq:PDE_constraints}
    \begin{align}
        \mathcal{L}\left[u(x,t),y_f(x,t)\right] &= 0, \quad \forall x \in \Omega, \; t \in [0,T],  \\
        \mathcal{B}\left[u(x,t),y_b(x,t)\right] &= 0, \quad \forall x \in \partial \Omega, \; t \in [0,T], \\
        \mathcal{I}\left[u(x,0),y_i(x)\right] &= 0, \quad \forall x \in \Omega.
    \end{align}
\end{subequations}
Here, $x$ and $t$ denote the spatial and temporal variables, respectively. $\mathcal{L}$ is an operator involving differentials that represents the PDE to be satisfied by $u$, $\mathcal{B}$ denotes the boundary condition, and $\mathcal{I}$ represents the initial condition. Common choices for boundary conditions include Dirichlet, Neumann, or Robin conditions, and our framework imposes no specific restrictions on these choices. The spaces $\mathrm{U}$ and $\mathrm{Y}$, in which $u$ and $y$ reside, are selected to ensure the well-posedness of the PDE problem. For instance, if $\mathcal{L}(u,y) = \Delta u - y$, corresponding to a standard second-order elliptic equation with no initial condition and a homogeneous Dirichlet boundary condition, appropriate spaces to consider are $\mathrm{U} = H^1_0(\Omega)$ and $\mathrm{Y} = L^2(\Omega)$ \cite{evans2022partial}.

If we set $y = (y_f, y_b, y_i)$ and define the constraint $F(u,y) = \left(\mathcal{L}(u,y), \mathcal{B}(u,y), \mathcal{I}(u,y)\right)$, we recover the constraint given in \eqref{eq:constrained_problem}. Specifically, the entire system \eqref{eq:PDE_constraints} or any individual equation within it can be viewed as a specific example of the general form of constraints $F(u,y) = 0$.

By setting different components of $y$ in \eqref{eq:PDE_constraints} to be tunable, we obtain various types of control problems. For example, if we consider $y_b$ and $y_i$ to be fixed and take $y_f$ as a tunable control, we obtain a distributed control problem, initially introduced in \cite{lions1971optimal}. The objective function to be minimized in this case is:
\begin{equation}
    \label{eq:distributed_optimal_control_objective}
    J_d(u,y) = \frac{1}{2}\|u - \hat{u}\|^2 + \frac{\rho}{2}\|y_f\|^2,
\end{equation}
where $\hat{u}$ denotes the desired state that we aim for our solution of the PDE system to approximate by tuning the value of $y_f$. The second term in this functional is a Tikhonov regularization term \cite{golub1999tikhonov}. Generally, the problem can be ill-posed without such a regularization term, and the Tikhonov regularization parameter $\rho$ value is typically chosen in advance \cite{haber2001preconditioned,rees2010optimal}. 


If we consider $y_b$ to be tunable and $y_f$ and $y_i$ to be fixed, then we obtain a boundary control problem, which minimizes the following objective function:
\begin{equation}
    \label{eq:boundary_optimal_control_objective}
        J_b(u,y) = \frac{1}{2}\|u - \hat{u}\|^2 + \frac{\rho}{2}\|y_b\|^2,
\end{equation}
with the same notation for $\hat{u}$ and $\rho$. Similarly, one can define an initial value optimal control problem. 

Additionally, to clarify our terminology: from now on, we refer to functionals like $J(u,y)$ as objective functionals, as they represent the goal that we aim to minimize with $(u,y)$ found by our algorithms. On the other hand, we refer to functionals like $\Plambda(u,y)$ as cost functionals or loss functionals.


It is important to note that the optimal control problems listed here only encompass some possible applications of our proposed method. General PDE-constrained optimization problems can be adapted to fit within our framework. Our methodology can be viewed as a variant of the penalty method. As long as a problem can be resolved or approximated using the penalty method, it is feasible to implement our approach. This study will focus on distributed and boundary optimal control problems because they are typically classic and representative examples.

\section{Discretized Problem}\label{sec:linear_analysis}

This section is devoted to discussing our motivation for setting up our method. In PDE-constrained optimization problems, there is ongoing debate on whether to use the discretize-then-optimize or optimize-then-discretize strategy \cite{becker2007optimal,collis2002analysis}. Since our method is based on a neural network, the autodifferentiation technique \cite{baydin2018automatic} naturally leads to a discretized system to solve. We adopt the discretize-then-optimize approach, beginning with an analysis of a discretized problem.

We consider the following discretized constrained optimization problem: find $u \in \mathbb{R}^n$ and $y \in \mathbb{R}^m$ to minimize
\begin{equation}\label{eq:discretized_constrained_problem}
    J(u,y) = \frac{1}{2}\|Au - b\|^2 + \frac{\rho}{2}\|y\|^2, 
\end{equation}
subject to
\begin{equation}
    \label{eq:discretized_constraint}
    Ku = y,
\end{equation}
where $0 < m \leq n$, $A \in \mathbb{R}^{k \times n}$, $K \in \mathbb{R}^{m \times n}$ with $k>0$ are matrices, $b \in \mathbb{R}^k$ is a given vector, and $\rho > 0$ is a given Tikhonov regularization parameter. This is a discretized version of problem \eqref{eq:constrained_problem}, with the objective functional chosen to match the type of optimal control problem set up in Section \ref{sec:problem_setup}. As a natural choice for discretizing the optimal control problem, we can take $A = I$, where $I$ is the $n \times n$ identity matrix. However, for the generality of our analysis, we only require that the row vectors of $A$ and $K$ span $\mathbb{R}^n$. Under this assumption, we know that for any $\alpha > 0$, the matrix $G_\alpha$ is invertible, where $G_\alpha$ is defined as
\begin{equation*}
    G_\alpha = A^T A + \alpha K^T K.
\end{equation*}

We will start from here to discuss how to solve problem \eqref{eq:discretized_constrained_problem}-\eqref{eq:discretized_constraint}.

\subsection{Explicit Solution}\label{subsec:explicit_solution}
We note that problem \eqref{eq:discretized_constrained_problem}-\eqref{eq:discretized_constraint} has an explicit solution that can be computed using a Lagrangian formulation. We consider the Lagrangian form:
\begin{equation}\label{eq:Lagrangian_form}
    L(u,y) = J(u,y) + \zeta (Ku - y) = \frac{1}{2}\|Au - b\|^2 + \frac{\rho}{2}\|y\|^2 + \zeta^T (Ku - y),
\end{equation}
with $\zeta\in\mathbb{R}^m$ as an auxiliary Lagrangian parameter vector. By differentiating \eqref{eq:Lagrangian_form} with respect to $u, y, \zeta$ respectively, the first-order optimality conditions \cite{rockafellar1993lagrange} are given by the following equations:
\begin{equation*}
    \left\{
        \begin{array}{l}
            A^T A u - A^T b + K^T \zeta = 0 ,\\
            \rho y - \zeta = 0, \\
            Ku - y = 0.
        \end{array}
    \right.
\end{equation*}
Solving this system results in:
\begin{equation}
    \label{eq:analytical_solution}
    \left\{
        \begin{array}{l}
           \hat{u} = \left(A^T A + \rho K^T K\right)^{-1} A^T b, \\
            \hat{y} = K\hat{u} = K \left(A^T A + \rho K^T K\right)^{-1} A^T b.
        \end{array}
    \right.
\end{equation}
Thus, we have found $(\hat{u}, \hat{y})$ to be the analytical solution to problem \eqref{eq:discretized_constrained_problem}-\eqref{eq:discretized_constraint}, and this notation will continue to be used throughout this paper. However, computing the inverse matrix can be challenging when dealing with large-scale systems, and its numerical stability may become an issue \cite{trefethen2022numerical}. Therefore, in practice, the penalty method is often preferred for implementation, as it is more suitable for applying iterative methods that do not require direct computation of the inverse \cite{bertsekas2014constrained,hestenes1969multiplier}. As introduced in \eqref{eq:unconstrained_penalty_problem}, the problem is formulated as follows: Find $u \in \mathbb{R}^n$ and $y \in \mathbb{R}^m$ to minimize
\begin{equation}
    \label{eq:penalty_formulation}
    \Plambda(u,y) = J(u,y) + \frac{\lambda}{2} \|Ku - y\|^2 = \frac{1}{2}\|Au - b\|^2 + \frac{\rho}{2}\|y\|^2 + \frac{\lambda}{2} \|Ku - y\|^2,
\end{equation}
where $\lambda > 0$ is a penalty parameter. For simplicity of notation, we will denote the remainder function corresponding to the constraint \eqref{eq:discretized_constraint} as $R(u,y)$, defined as
\begin{equation}
    \label{eq:discretized_constraint_definition}
    R(u,y) = \|Ku - y\|^2.
\end{equation} 
The first-order optimality conditions to minimize \eqref{eq:penalty_formulation} are given by:
\begin{equation*}
    \left\{
        \begin{array}{l}
            A^T A u - A^T b + \lambda K^T K u - \lambda K^T y = 0 ,\\
            \rho y - \lambda Ku + \lambda y = 0.
        \end{array}
    \right.
\end{equation*}
Solving this system results in:
\begin{equation}
    \label{eq:analytical_solution_lambda}
    \left\{
        \begin{array}{l}
            u^\lambda = \left(A^T A + \frac{\rho \lambda}{\rho + \lambda} K^T K\right)^{-1} A^T b ,\\
            y^\lambda = \frac{\lambda}{\rho + \lambda} Ku^{\lambda} = \frac{\lambda}{\rho + \lambda} K \left(A^T A + \frac{\rho \lambda}{\rho + \lambda} K^T K\right)^{-1} A^T b.
        \end{array}
    \right.
\end{equation}
Comparing \eqref{eq:analytical_solution_lambda} with \eqref{eq:analytical_solution}, we observe that as $\lambda$ tends to infinity, $u^\lambda$ will converge to $\hat{u}$ because $\lim\limits_{\lambda \to \infty} \frac{\rho \lambda}{\rho + \lambda} = \rho$, and $y^\lambda$ will converge to $\hat{y}$. Therefore, although problem \eqref{eq:penalty_formulation} is distinct from problem \eqref{eq:discretized_constrained_problem}-\eqref{eq:discretized_constraint}, we can consider \eqref{eq:penalty_formulation} with a sufficiently large $\lambda$ as an acceptable approximation to problem \eqref{eq:discretized_constrained_problem}-\eqref{eq:discretized_constraint}.

However, in practice, as $\lambda$ increases, the difficulty of solving the unconstrained optimization problem associated with \eqref{eq:penalty_formulation} also increases. Mathematically, this can be observed by examining the Hessian matrix for $\Plambda(u,y)$ with respect to $u$ and $y$, denoted as $\Hlambda$. Direct computation shows that
\begin{equation}
    \label{eq:Hessian_lambda}
    \Hlambda = \begin{pmatrix}
        A^TA+\lambda K^TK & -\lambda K^T\\
       -\lambda K & (\rho +\lambda)I_m,
    \end{pmatrix}
\end{equation}
where $I_m$ denotes an $m\times m$ identity matrix. As $\lambda$ tends to infinity, $\Hlambda$ tends to become a singular matrix, making it difficult to find the minimizer of $\Plambda(u,y)$. 

\begin{figure}[ht]
    \centering
    \begin{minipage}{0.48\textwidth}
        \centering
        \includegraphics[width=\textwidth]{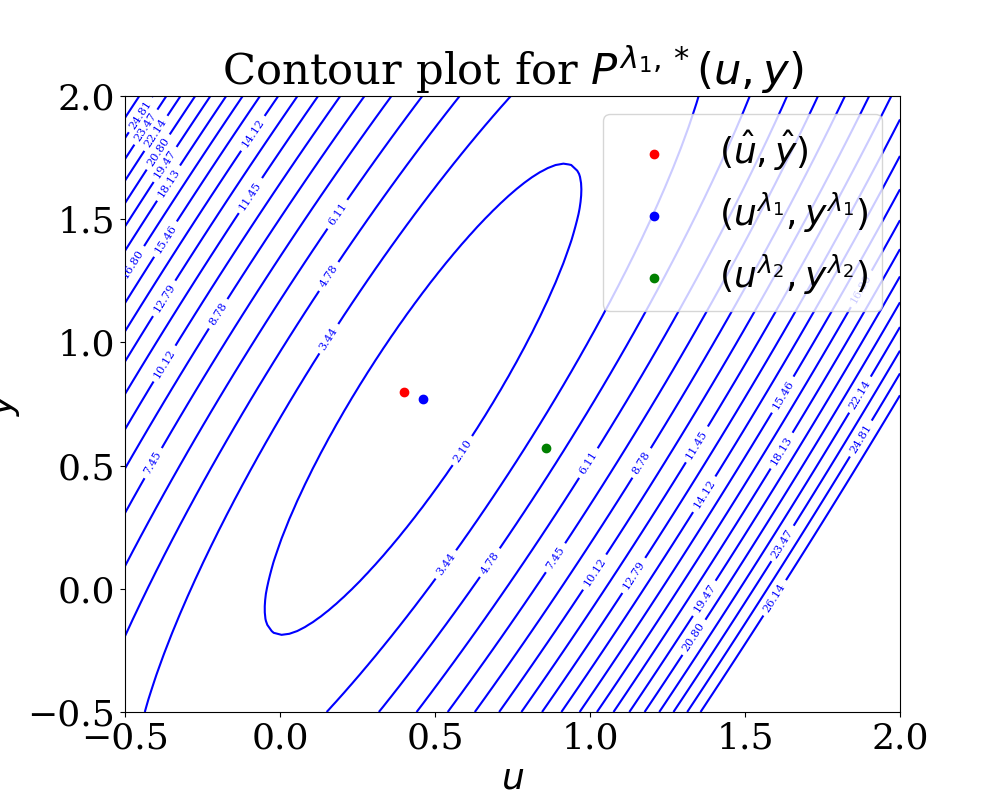}
        \caption*{(a) Level Sets of $P^{\lambda,*}$ when $\lambda = 5$}
    \end{minipage}%
    \hspace{0.03\textwidth}
    \begin{minipage}{0.48\textwidth}
        \centering
        \includegraphics[width=\textwidth]{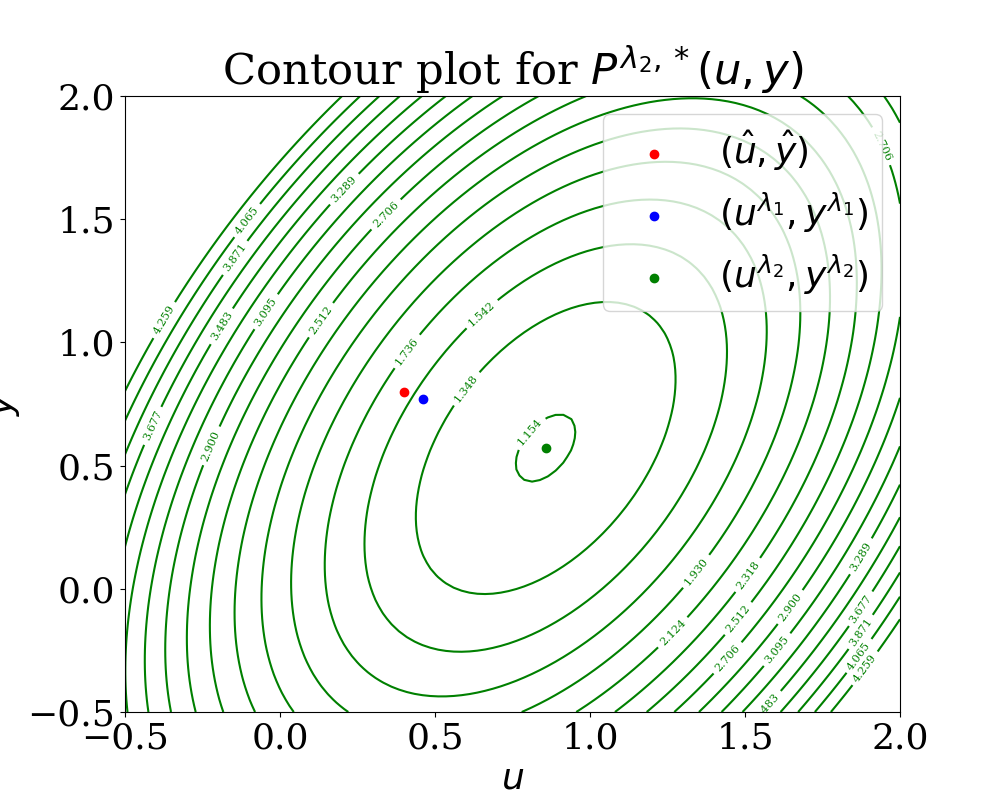}
        \caption*{(b) Level Sets of $P^{\lambda,*}$ when $\lambda = 0.5$}
    \end{minipage}
    \caption{Comparison of Contour Plot for $P^{\lambda,*}(u,y)$ Defined in \eqref{eq:penalty_formulation_contour_example} with Different Values of $\lambda$}
    \label{fig:contour_different_lambda} 
\end{figure}

To illustrate this, we present a graphical example. In Figure \ref{fig:contour_different_lambda}, we consider a one-dimensional problem with $m=n=1$ and aim to find $u, y\in\mathbb{R}$ that minimize
\begin{equation}\label{eq:discretized_constrained_problem_contour_example}
    J^*(u,y) = \frac{1}{2}\lvert u - 2\vert^2 + \frac{1}{2}\lvert y\rvert^2, \quad \text{subject to } 2u = y.
\end{equation}
The corresponding penalty formulation is to find $u, y\in\mathbb{R}$ such that they minimize
\begin{equation}
    P^{\lambda,*}(u,y) = \frac{1}{2}\lvert u - 2\vert^2 + \frac{1}{2}\lvert y\rvert^2 + \frac{\lambda}{2}\lvert 2u-y\rvert^2.\label{eq:penalty_formulation_contour_example}
\end{equation}
Consistent with the notation used above, we denote the exact solution as $(\hat{u}, \hat{y})$ and mark it as a red point in the figure. The solution for the minimization problem with $\lambda_1 = 5$, denoted as $(\ulambdaone, \ylambdaone)$, is marked as a blue point, while the solution for the minimization problem with $\lambda_2 = 0.5$, denoted as $(\ulambdatwo, \ylambdatwo)$, is marked as a green point. The range for $u$ and $y$ is chosen to be $[-0.5, 2]$. We can observe that $(\ulambdaone, \ylambdaone)$ is closer to the exact solution than $(\ulambdatwo, \ylambdatwo)$. 

Additionally, Figure \ref{fig:contour_different_lambda}, parts (a) and (b), plot the contour of each level set of $P^{\lambda,*}$ for these two values of $\lambda$. When $\lambda$ is relatively larger, the contour is more flattened, indicating that the conditioning of the problem is worse \cite{shanno1978matrix}. This results in greater difficulty in finding the optimal value, which aligns with our analysis above.

In a neural network setting, a large $\lambda$ makes the network difficult to train and may not yield the correct solution. Conversely, if $\lambda$ is not sufficiently large, the solution to the corresponding unconstrained problem may not satisfy the constraint adequately, and $u^\lambda$ might not approximate $u^*$ well. Therefore, we seek a practical method that ensures the solution satisfies the constraints while effectively approximating the true solution.

\subsection{Penalty Adversarial Problem}
To resolve the problem mentioned in the end of last subsection, now, we propose a new unconstrained optimization approach based on the penalty method. Instead of using a single penalty parameter \(\lambda\), we simultaneously consider two problems with different penalty parameters. Let \(\lambda_1 > \lambda_2 > 0\), and let \((u^{\lambda_1}, y^{\lambda_1})\) and \((u^{\lambda_2}, y^{\lambda_2})\) denote the solutions that minimize the functionals \(\mathcal{P}^{\lambda_1}(u,y)\) and \(\mathcal{P}^{\lambda_2}(u,y)\), respectively. We assume that the former problem is hard to solve in practice, while the latter is easier to solve. The corresponding objective functionals for \((u^{\lambda_1}, y^{\lambda_1})\) and \((u^{\lambda_2}, y^{\lambda_2})\) can be computed using \eqref{eq:discretized_constrained_problem} and are denoted as \(J(u^{\lambda_1}, y^{\lambda_1})\) and \(J(u^{\lambda_2}, y^{\lambda_2})\). We now consider the following problem:
Find \(u \in \mathbb{R}^n\) and \(y \in \mathbb{R}^m\) to minimize
\begin{equation}
    \label{eq:penalty_adversarial_linearized_problem}
   \Alambdaonelambdatwo(u,y) = \left\{
        \begin{aligned}
            &J(u,y) + \frac{\lambda_1}{2} R(u,y) + \omega \left[J(u, y) - J(u^{\lambda_2}, y^{\lambda_2})\right]^2, &\quad \text{if } (u,y) \in \Omega_1,\\
            &J(u,y) + \frac{\lambda_1}{2} R(u,y), &\quad \text{if } (u,y) \in \Omega_2.
        \end{aligned}
    \right.
\end{equation}
Here, \(\omega > 0\) is a tunable parameter, and 
\begin{equation*}
    \Omega_1 = \{(u,y): J(u, y) > J(u^{\lambda_2}, y^{\lambda_2})\}, \quad \Omega_2 = \{(u,y): J(u, y) \leq J(u^{\lambda_2}, y^{\lambda_2})\}.
\end{equation*}
Clearly, \(\Omega_1\) is an open set and \(\Omega_2\) is a closed set due to the continuity of \(J(u,y)\). We call this minimization problem the penalty adversarial problem (PAP).

Compared to the standard penalty problem, $\Alambdaonelambdatwo(u,y)$ incorporates an additional term to penalize the failure of achieving a sufficiently small objective functional. By designing this problem, we aim to find a solution that adheres closely to the constraint while maintaining feasibility in its solvability. The problem's feasibility is difficult to assert through a mathematical criterion, as it is highly problem-dependent. However, we will demonstrate this advantage through graphical examples in Section \ref{subsubsec:graphical_A_Plambda} and numerical examples in Section \ref{sec:numerics}. On the other hand, a mathematical measure to determine if a solution better adheres to the constraints is straightforward to obtain, which is the value of the remainder function evaluated at the solution. Therefore, mathematically, we aim to show that by choosing an appropriate value of \(\omega\), the minimizer of $\Alambdaonelambdatwo(u,y)$, denoted as $(\ulambdaonelambdatwoomega, \ylambdaonelambdatwoomega)$ hereafter, will satisfy the following:
\begin{equation*}
    R(\ulambdaonelambdatwoomega, \ylambdaonelambdatwoomega) < R(\ulambdatwo, \ylambdatwo).
\end{equation*}
Here $R(u,y)$ is defined in \eqref{eq:discretized_constraint_definition}. As long as this holds, it can be seen that the minimizer for $\Alambdaonelambdatwo(u,y)$ is more accurate than the minimizer for $\Plambdatwo(u,y)$ in terms of satisfying the PDE constraints. Proving this argument will be the main focus of our analysis hereafter.

First, consider the relationship between $J(u,y)$ and $R(u,y)$: For a fixed $\lambda \in \mathbb{R}^+$, $(u^\lambda, y^\lambda)$ minimizes $\Plambda(u,y)$, which is a combination of two parts: one involves the objective functional $J(u,y)$, while the other comprises the remainder function $R(u,y)$. If the value of $J(u,y)$ is fixed, then minimizing $\Plambda(u,y)$ is equivalent to minimizing $R(u,y)$. Therefore, intuitively, the fact that $(u^\lambda, y^\lambda)$ minimizes $\Plambda(u,y)$ indicates that $(u^\lambda, y^\lambda)$ minimizes the remainder function $R(u,y)$ among all $(u,y) \in \mathbb{R}^n \times \mathbb{R}^m$ such that $J(u,y) = J(u^\lambda, y^\lambda)$. The following lemma justifies this.

\begin{lem}
    \label{lem:R_J_relation}
    For a fixed $\lambda > 0$, let $(u^\lambda, y^\lambda) \in \mathbb{R}^n \times \mathbb{R}^m$ minimize $P^\lambda(u,y)$ as given in \eqref{eq:penalty_formulation}. If $(u,y) \in \mathbb{R}^n \times \mathbb{R}^m$, then we have the following assertions:
    \begin{enumerate}
        \item If $J(u,y) \leq J(u^\lambda, y^\lambda)$, then $R(u,y) \geq R(u^\lambda, y^\lambda)$.
        \item If $R(u,y) \leq R(u^\lambda, y^\lambda)$, then $J(u,y) \geq J(u^\lambda, y^\lambda)$.
        \item In each of the previous two assertions, equality can only be attained simultaneously. Specifically, if $J(u,y) < J(u^\lambda, y^\lambda)$, then $R(u,y) > R(u^\lambda, y^\lambda)$. If $R(u,y) < R(u^\lambda, y^\lambda)$, then $J(u,y) > J(u^\lambda, y^\lambda)$.
    \end{enumerate}
\end{lem}
\begin{proof}
    We will only prove the first argument and the proof for the second will follow in the same way. By definition, $(u^\lambda, y^\lambda)$ minimizes $\Plambda(u,y)$. Thus, for any $(u,y) \in \mathbb{R}^n \times \mathbb{R}^m$, we have $\Plambda(u,y) \geq \Plambda(u^\lambda, y^\lambda)$. This implies that
    \begin{equation}\label{eq:J_lambda_quantity_relation}
        J(u,y) + \frac{\lambda}{2} R(u,y) \geq J(u^\lambda, y^\lambda) + \frac{\lambda}{2} R(u^\lambda, y^\lambda).
    \end{equation}
    Since $J(u,y) \leq J(u^\lambda, y^\lambda)$, it follows that
    \begin{equation*}
        R(u,y) \geq R(u^\lambda, y^\lambda),
    \end{equation*}
    as $\lambda > 0$. The third argument also follows immediately from \eqref{eq:J_lambda_quantity_relation} when the inequality relation between $J(u,y)$ and $J(u^\lambda, y^\lambda)$ is strict.
\end{proof}

This result indeed provides a sufficient condition to restrict $(u,y)$ in domain $\Omega_1$. Namely, using the third argument of Lemma \ref{lem:R_J_relation}, as long as $R(u,y) < R(\ulambdatwo, \ylambdatwo)$, then $(u,y) \in \Omega_1$ and the value of $\Alambdaonelambdatwo(u,y)$ will differ from $\Plambdaone(u,y)$ at these points. Another important deduction we can make from here is that the minimum of $\Alambdaonelambdatwo(u,y)$ in the closed region $\Omega_2$ could always be attained at $(\ulambdatwo, \ylambdatwo)$. We conclude this in the following result:

\begin{prop}
    \label{prop:T_lambdaone_lambdatwo_minimum_omega_2}
    If $(u,y) \in \Omega_2$, then 
    \begin{equation*}
        \Alambdaonelambdatwo(u,y) \geq J(\ulambdatwo, \ylambdatwo) + \frac{\lambda_1}{2} R(\ulambdatwo, \ylambdatwo).
    \end{equation*}
    In other words,
    \begin{equation*}
        \min_{(u,y) \in \Omega_2} \Alambdaonelambdatwo(u,y) = J(\ulambdatwo, \ylambdatwo) + \frac{\lambda_1}{2} R(\ulambdatwo, \ylambdatwo).
    \end{equation*}
\end{prop}
\begin{proof}
    $(u,y) \in \Omega_2$ implies that $J(u,y) \leq J(\ulambdatwo, \ylambdatwo)$. According to the first assertion in Lemma \ref{lem:R_J_relation}, we know that
    \begin{equation*}
        R(u,y) \geq R(\ulambdatwo, \ylambdatwo).
    \end{equation*}
    On the other hand, since $(\ulambdatwo, \ylambdatwo)$ minimizes $\Plambdatwo(u,y)$, we have
    \begin{equation*}
        J(u,y) + \frac{\lambda_2}{2} R(u,y) = \Plambdatwo(u,y) \geq \Plambdatwo(\ulambdatwo, \ylambdatwo) = J(\ulambdatwo, \ylambdatwo) + \frac{\lambda_2}{2} R(\ulambdatwo, \ylambdatwo).
    \end{equation*}
    Combining these two inequalities, we get
    \begin{equation*}
        \begin{aligned}
            J(u,y) + \frac{\lambda_1}{2} R(u,y) &= J(u,y) + \frac{\lambda_2}{2} R(u,y) + \frac{\lambda_1 - \lambda_2}{2} R(u,y)\\
            &\geq J(\ulambdatwo, \ylambdatwo) + \frac{\lambda_2}{2} R(\ulambdatwo, \ylambdatwo) + \frac{\lambda_1 - \lambda_2}{2} R(\ulambdatwo, \ylambdatwo)\\
            &= J(\ulambdatwo, \ylambdatwo) + \frac{\lambda_1}{2} R(\ulambdatwo, \ylambdatwo).
        \end{aligned}
    \end{equation*}
    This proves the claim.
\end{proof}

Revealed by this, if for any $(u,y) \in \Omega_1$, $\Alambdaonelambdatwo(u,y)$ is always greater than $\left[J(u^{\lambda_2}, y^{\lambda_2}) + \frac{\lambda_1}{2} R(\ulambdatwo,\ylambdatwo)\right]$, then $(u^{\lambda_2}, y^{\lambda_2})$ can be a solution to the minimization problem related to \eqref{eq:penalty_adversarial_linearized_problem}. In this case, our claim that $(\ulambdaonelambdatwoomega, \ylambdaonelambdatwoomega)$, the minimizer of this problem, will always adhere more closely to the constraint than $(\ulambdatwo, \ylambdatwo)$ would no longer hold. Therefore, to make our claim valid, it is necessary for the minimizer $(\ulambdaonelambdatwoomega, \ylambdaonelambdatwoomega)$ to fall within $\Omega_1$. Enlightened by this, we present the following result:

\begin{lem}\label{lem:T_minimizer_range}
Assuming $\omega > 0$ is a fixed positive constant parameter, the following two arguments are equivalent:

\begin{enumerate}
    \item $(u,y) \in \Omega_1$ and 
    \begin{equation}\label{eq:uy_total_loss_condition}
        \Alambdaonelambdatwo(u,y) < J(u^{\lambda_2}, y^{\lambda_2}) + \frac{\lambda_1}{2} R(\ulambdatwo, \ylambdatwo),
    \end{equation}
    \item $(u,y)$ satisfies $R(u,y) < R(u^{\lambda_2}, y^{\lambda_2})$ and
    \begin{equation}\label{eq:uy_target_bound}
        J(u^{\lambda_2}, y^{\lambda_2}) < J(u,y) < J(u^{\lambda_2}, y^{\lambda_2}) + \frac{\lambda_1 \left[ R(u^{\lambda_2}, y^{\lambda_2}) - R(u,y) \right]}{1 + \sqrt{1 + 2\omega \lambda_1 \left[ R(u^{\lambda_2}, y^{\lambda_2}) - R(u,y) \right]}}.
    \end{equation}
\end{enumerate}
\end{lem}

\begin{proof}
    We first assume the former argument holds and try to deduce the latter one. To start, we define the difference between $J(u,y)$ and $J(u^{\lambda_2}, y^{\lambda_2})$ as $D^{\lambda_1, \lambda_2}(u,y)$, namely,
    \begin{equation*}
        D^{\lambda_1, \lambda_2}(u,y) = J(u,y) - J(u^{\lambda_2}, y^{\lambda_2}).
    \end{equation*}
    As $(u,y) \in \Omega_1$, $(u,y)$ satisfies
    \begin{equation*}
        D^{\lambda_1, \lambda_2}(u,y) > 0.
    \end{equation*}
    Using the definition of $\Alambdaonelambdatwo(u,y)$ given in \eqref{eq:penalty_adversarial_linearized_problem}, we see that \eqref{eq:uy_total_loss_condition} is equivalent to
    \begin{equation}\label{eq:uy_total_loss_condition_equivalence}
        \omega \left[ D^{\lambda_1, \lambda_2}(u,y) \right]^2 + D^{\lambda_1, \lambda_2}(u,y) + \frac{\lambda_1}{2} \left[ R(u,y) - R(u^{\lambda_2}, y^{\lambda_2}) \right] < 0.
    \end{equation}
    Since $\omega > 0$, this inequality will have solutions only when the discriminant of the corresponding quadratic formula is greater than zero, which is equivalent to:
    \begin{equation*}
        1 - 2\omega \lambda_1 \left[ R(u,y) - R(u^{\lambda_2}, y^{\lambda_2}) \right] > 0.
    \end{equation*}
    In this case, $D^{\lambda_1, \lambda_2}(u,y)$ should satisfy
    \begin{equation*}
        \frac{-1 - \sqrt{1 + 2\omega \lambda_1 \left[ R(u^{\lambda_2}, y^{\lambda_2}) - R(u,y) \right]}}{2\omega} < D^{\lambda_1, \lambda_2}(u,y) < \frac{-1 + \sqrt{1 + 2\omega \lambda_1 \left[ R(u^{\lambda_2}, y^{\lambda_2}) - R(u,y) \right]}}{2\omega}.
    \end{equation*}
    Since we require that $D^{\lambda_1, \lambda_2}(u,y) > 0$, it is necessary to have
    \begin{equation*}
        -1 + \sqrt{1 + 2\omega \lambda_1 \left[ R(u^{\lambda_2}, y^{\lambda_2}) - R(u,y) \right]}>0,
    \end{equation*}
    which is equivalent to
    \begin{equation*}
        R(u,y) < R(u^{\lambda_2}, y^{\lambda_2}).
    \end{equation*}
    With this, $J(u,y)$ can be estimated as
    \begin{equation*}
        J(u,y) = J(\ulambdatwo, \ylambdatwo) + D^{\lambda_1, \lambda_2}(u,y) < J(\ulambdatwo, \ylambdatwo) + \frac{\sqrt{1 + 2\omega \lambda_1 \left[ R(u^{\lambda_2}, y^{\lambda_2}) - R(u,y) \right]} - 1}{2\omega},
    \end{equation*}
    which is equivalent to \eqref{eq:uy_target_bound}.
    This proves the latter argument. 
    
    For the other direction, Lemma \ref{lem:R_J_relation} ensures that $R(u,y) < R(\ulambdatwo, \ylambdatwo)$ implies $(u,y) \in \Omega_1$. In addition, \eqref{eq:uy_target_bound} implies that the value of $D^{\lambda_1, \lambda_2}(u,y)$ will ensure \eqref{eq:uy_total_loss_condition_equivalence}, which is equivalent to \eqref{eq:uy_total_loss_condition}, thereby completing the proof.
\end{proof}

Here, we have provided an equivalent condition to characterize the case where the minimum of $\Alambdaonelambdatwo(u,y)$ falls in $\Omega_1$. As long as this condition is met, the minimum point $(\ulambdaonelambdatwoomega, \ylambdaonelambdatwoomega)$ will correspond to a smaller value of $R(u,y)$ compared to $R(\ulambdatwo,\ylambdatwo)$, which is the desired outcome. The remaining task is to justify this condition. To do so, we need to demonstrate that there exists a point in $\Omega_1$ that satisfies \eqref{eq:uy_target_bound}, and then we can invoke the continuity of $J(u,y)$ to conclude the existence of a minimum point. This leads to the main result of this section.
\begin{thm}\label{thm:main_theorem_pan_discretized_problem}
    If $K\ulambdatwo \neq 0$ and there exists $\sigma > 0$ such that
    \begin{equation}\label{eq:minimizer_existence_condition}
        \frac{\lambda_1}{2} R(u^{\lambda_2}, y^{\lambda_2}) > \sigma + J(\hat{u}, \hat{y}) - J(\ulambdatwo, \ylambdatwo),
    \end{equation}
    where $(\hat{u}, \hat{y})$ given in \eqref{eq:analytical_solution} is the exact solution that minimizes $J(u,y)$ subject to $Ku = y$, then there exists $\omega > 0$ such that $(\ulambdaonelambdatwoomega, \ylambdaonelambdatwoomega)$, which minimizes $\Alambdaonelambdatwo(u,y)$ defined in \eqref{eq:penalty_adversarial_linearized_problem}, satisfies
    \begin{equation*}
        R(\ulambdaonelambdatwoomega, \ylambdaonelambdatwoomega) < R(u^{\lambda_2}, y^{\lambda_2}).
    \end{equation*}
\end{thm}

\begin{proof}
    Recalling the definition for $(\hat{u}, \hat{y})$, we will show that $(\hat{u}, \hat{y})$ satisfies 
    \begin{equation} \label{eq:uhat_yhat_smaller_condition}
        \Alambdaonelambdatwo(\hat{u}, \hat{y}) < J(u^{\lambda_2}, y^{\lambda_2}) + \frac{\lambda_1}{2} R(\ulambdatwo, \ylambdatwo).
    \end{equation}
    Using Lemma \ref{lem:T_minimizer_range}, it is equivalent to show $R(\hat{u}, \hat{y}) < R(\ulambdatwo, \ylambdatwo)$ and $(\hat{u}, \hat{y})$ satisfies \eqref{eq:uy_target_bound}.
    
    It is clear that $R(\hat{u}, \hat{y}) = 0$. On the other hand,
    \begin{equation*}
    \begin{aligned}
        R(u^{\lambda_2}, y^{\lambda_2}) &= \|Ku^{\lambda_2} - y^{\lambda_2}\|^2 = \left\|Ku^{\lambda_2} - \frac{\lambda_2}{\rho + \lambda_2} Ku^{\lambda_2}\right\|^2
        = \left\|\frac{\rho}{\rho + \lambda_2} Ku^{\lambda_2}\right\|^2
        > 0,
    \end{aligned}
    \end{equation*}
    according to our assumption. Thus, the relation $R(\hat{u}, \hat{y}) < R(\ulambdatwo, \ylambdatwo)$ holds. 
    
    Meanwhile, using the fact that $R(\hat{u}, \hat{y}) = 0$ again, \eqref{eq:uy_target_bound} in this case reduces to 
    \begin{equation*}
        J(u^{\lambda_2}, y^{\lambda_2}) < J(\hat{u},\hat{y}) < J(u^{\lambda_2}, y^{\lambda_2}) + \frac{\lambda_1 R(u^{\lambda_2}, y^{\lambda_2})}{1 + \sqrt{1 + 2\omega \lambda_1 R(u^{\lambda_2}, y^{\lambda_2})}}.
    \end{equation*}
    Since $\sigma > 0$ and \eqref{eq:minimizer_existence_condition} holds, we can find $\omega > 0$ to satisfy this condition. Therefore, we have shown \eqref{eq:uhat_yhat_smaller_condition}. This indicates that $(\ulambdatwo, \ylambdatwo)$ is not the minimizer of $\Alambdaonelambdatwo(u,y)$ for small enough positive $\omega$, as $\Alambdaonelambdatwo(u,y)$ reaches a smaller value at $(\hat{u}, \hat{y})$ compared to $(\ulambdatwo, \ylambdatwo)$. It remains to show that any minimizer of $\Alambdaonelambdatwo(u, y)$ satisfies $R(\ulambdaonelambdatwoomega, \ylambdaonelambdatwoomega) < R(u^{\lambda_2}, y^{\lambda_2})$. 

    To see this, we consider a subdomain of $\overline{\Omega}_1 = \left\{(u,y) : J(u,y) \geq J(u^{\lambda_2}, y^{\lambda_2})\right\}$, the closure of $\Omega_1$, defined as
    \begin{equation*}
        \Omega^{\lambda_1, \lambda_2} = \left\{(u,y) : J(u^{\lambda_2}, y^{\lambda_2}) \leq J(u,y) \leq J(u^{\lambda_2}, y^{\lambda_2}) + \frac{\lambda_1}{2} R(u^{\lambda_2}, y^{\lambda_2})\right\}.
    \end{equation*}
    Due to the continuity of $J(u,y)$ with respect to $(u,y)$ and the definition of $J(u,y)$, we see that $\Omega^{\lambda_1, \lambda_2}$ is a closed and bounded domain. Hence, it is compact. When $(u,y) \in \overline{\Omega}_1 \setminus \Omega^{\lambda_1, \lambda_2}$,
    \begin{equation*}
        \Alambdaonelambdatwo(u,y) \geq J(u,y) > J(u^{\lambda_2}, y^{\lambda_2}) + \frac{\lambda_1}{2} R(u^{\lambda_2}, y^{\lambda_2}) = \Alambdaonelambdatwo(\ulambdatwo, \ylambdatwo),
    \end{equation*}
    and so it cannot be a minimizer of $\Alambdaonelambdatwo(u,y)$ in $\Omega_1$. On the other hand, as a continuous function defined on a compact set $\Omega^{\lambda_1, \lambda_2}$, $\Alambdaonelambdatwo(u,y)$ attains its minimum at some point in $\Omega^{\lambda_1, \lambda_2}$. For consistent notation, let us denote this point as $(\ulambdaonelambdatwoomega, \ylambdaonelambdatwoomega)$. Therefore, 
    \begin{equation}\label{eq:T_minimizer}
    \begin{aligned}
        J(\ulambdaonelambdatwoomega, \ylambdaonelambdatwoomega) + \frac{\lambda_1}{2} &R(\ulambdaonelambdatwoomega, \ylambdaonelambdatwoomega) + \omega \left[J(\ulambdaonelambdatwoomega, \ylambdaonelambdatwoomega) - J(u^{\lambda_2}, y^{\lambda_2})\right]^2\\ &= \Alambdaonelambdatwo(\ulambdaonelambdatwoomega, \ylambdaonelambdatwoomega) \leq \Alambdaonelambdatwo(\hat{u}, \hat{y}) < J(u^{\lambda_2}, y^{\lambda_2}) + \frac{\lambda_1}{2} R(u^{\lambda_2}, y^{\lambda_2}).
    \end{aligned}
    \end{equation}
    Finally, we will show that $R(\ulambdaonelambdatwoomega, \ylambdaonelambdatwoomega) < R(\ulambdatwo, \ylambdatwo)$. By Lemma \ref{lem:R_J_relation}, it is sufficient to show that $J(\ulambdaonelambdatwoomega, \ylambdaonelambdatwoomega) > J(\ulambdatwo, \ylambdatwo)$. Therefore, we need to prove that $J(\ulambdaonelambdatwoomega, \ylambdaonelambdatwoomega) \neq J(\ulambdatwo, \ylambdatwo)$ since $(\ulambdaonelambdatwoomega, \ylambdaonelambdatwoomega)\in  \Omega^{\lambda_1, \lambda_2} $. In fact, if $J(\ulambdaonelambdatwoomega, \ylambdaonelambdatwoomega) = J(\ulambdatwo, \ylambdatwo)$, then by Lemma \ref{lem:R_J_relation}, we know that $R(\ulambdaonelambdatwoomega, \ylambdaonelambdatwoomega) \geq R(\ulambdatwo, \ylambdatwo)$. From \eqref{eq:T_minimizer}, we have
    \begin{equation*}
    \begin{aligned}
        &\quad\,\quad\,J(\ulambdatwo, \ylambdatwo) + \frac{\lambda_1}{2} R(\ulambdatwo, \ylambdatwo)\\ &\leq J(\ulambdaonelambdatwoomega, \ylambdaonelambdatwoomega) + \frac{\lambda_1}{2} R(\ulambdaonelambdatwoomega, \ylambdaonelambdatwoomega) \\
        &= J(\ulambdaonelambdatwoomega, \ylambdaonelambdatwoomega) + \frac{\lambda_1}{2} R(\ulambdaonelambdatwoomega, \ylambdaonelambdatwoomega) + \omega \left[J(\ulambdaonelambdatwoomega, \ylambdaonelambdatwoomega) - J(u^{\lambda_2}, y^{\lambda_2})\right]^2 \\
        &=\Alambdaonelambdatwo(\ulambdaonelambdatwoomega, \ylambdaonelambdatwoomega)\leq \Alambdaonelambdatwo(\hat{u},\hat{y})< \Alambdaonelambdatwo(\ulambdatwo,\ylambdatwo)= J(u^{\lambda_2}, y^{\lambda_2}) + \frac{\lambda_1}{2} R(u^{\lambda_2}, y^{\lambda_2}).
    \end{aligned}
    \end{equation*}
    This is a contradiction, and so we conclude that $J(\ulambdaonelambdatwoomega, \ylambdaonelambdatwoomega)$ is strictly larger than $J(\ulambdatwo, \ylambdatwo)$ and thus $R(\ulambdaonelambdatwoomega, \ylambdaonelambdatwoomega)$ is strictly smaller than $R(\ulambdatwo, \ylambdatwo)$. This completes the proof.
\end{proof}

To conclude this part, we will make two comments on our strategy of proof for Theorem \ref{thm:main_theorem_pan_discretized_problem}.

\begin{rem}
    \label{rem:minimizer_finding_strategy} The proof for Theorem \ref{thm:main_theorem_pan_discretized_problem} presented here focuses on the existence of a minimizer, which is sufficient to support our claim that $R(\ulambdaonelambdatwoomega, \ylambdaonelambdatwoomega) < R(\ulambdatwo, \ylambdatwo)$. For a deeper study of the properties of $(\ulambdaonelambdatwoomega, \ylambdaonelambdatwoomega)$, one can consider a quantitative analysis of the functional form of $\Alambdaonelambdatwo(u,y)$. 
\end{rem}

\begin{rem}\label{rem:choice_of_reference_point}
    In this context, we use $(\hat{u}, \hat{y})$ as a reference point to demonstrate that there exists at least one point satisfying \eqref{eq:uy_target_bound}. However, this choice is not mandatory. In fact, any point in $\Omega_1$ that can be computed and shown to result in a value smaller than $\Alambdaonelambdatwo(\ulambdatwo, \ylambdatwo)$ could replace $(\hat{u}, \hat{y})$ in our proof.
\end{rem}

\subsection{Choice of $\omega$}

The last part discusses the possibility of designing a penalty adversarial problem to ensure its minimizer adheres to the constraints more closely than $(\ulambdatwo,\ylambdatwo)$. Once $\lambda_1$ and $\lambda_2$ are fixed, the choice of $\omega$ will determine the formulation of the corresponding problem. In this section, we will explore this choice's influence through theoretical analysis and computational examples.

\subsubsection{Upper Bound of $\omega$}

We note that \eqref{eq:minimizer_existence_condition} provides an explicit condition that can be used to determine if the value of $\lambda_1$ is large enough to implement this penalty adversarial strategy independent of $\omega$. Then, after fixing the values for $\lambda_1$ and $\lambda_2$, based on this proof, we can provide an estimate for the upper bound of $\omega$. This conclusion is formalized in the following proposition.

\begin{prop}\label{prop:omega_upper_bound}
    Assuming \eqref{eq:minimizer_existence_condition} holds, then $\omega > 0$ should satisfy the following relation:
\end{prop}
\begin{equation}\label{eq:omega_choice}
    \omega < \frac{\lambda_1 R(u^{\lambda_2}, \ylambdatwo) - 2\left[J(\hat{u}, \hat{y}) - J(\ulambdatwo, \ylambdatwo)\right]}{2\left[J(\hat{u}, \hat{y}) - J(\ulambdatwo, \ylambdatwo)\right]^2}.
\end{equation}
\begin{proof}
    In the proof of Theorem \ref{thm:main_theorem_pan_discretized_problem}, we deduced the following relation:
    \begin{equation*}
        J(\hat{u}, \hat{y}) < J(u^{\lambda_2}, y^{\lambda_2}) + \frac{\lambda_1 R(u^{\lambda_2}, y^{\lambda_2})}{1 + \sqrt{1 + 2\omega \lambda_1 R(u^{\lambda_2}, y^{\lambda_2})}}.
    \end{equation*}
    Rearranging this inequality results in \eqref{eq:omega_choice}.
\end{proof}

We will provide an intuitive understanding of why there exists an upper bound on the choice of $\omega$ in Section \ref{subsubsec:equivalence_auto_tuning}.

\subsubsection{Comparison between $\Alambdaonelambdatwo(u,y)$ and $\Plambda(u,y)$}\label{subsubsec:graphical_A_Plambda}


We have discussed the advantage of minimizing $\Alambdaonelambdatwo(u,y)$ over $\Plambdatwo(u,y)$, specifically its ability to adhere to the constraint more effectively. On the other hand, the main disadvantage of trying to minimize $\Plambdaone(u,y)$, as mentioned earlier, is its poor conditioning, which makes it difficult to solve. Therefore, by choosing to minimize $\Alambdaonelambdatwo(u,y)$ instead of $\Plambdaone(u,y)$, we aim to make the corresponding minimization problem easier to solve.

Here, we will present a graphical example to illustrate this point. Using the same example as in Section \ref{subsec:explicit_solution}, we consider a one-dimensional problem with $m = n = 1$ and aim to find $u, y \in \mathbb{R}$ to minimize

\begin{equation}\label{eq:discretized_constrained_problem_contour_example}
    J^*(u,y) = \frac{1}{2}\lvert u - 2\vert^2 + \frac{1}{2}\lvert y\rvert^2, \quad \text{subject to } 2u = y.
\end{equation}
The corresponding penalty formulation has been stated in \eqref{eq:penalty_formulation_contour_example}, and the corresponding penalty adversarial problem is as follows:
find $u \in \mathbb{R}$ and $y \in \mathbb{R}$ to minimize

\begin{equation}
    \label{eq:penalty_adversarial_linearized_problem_example}
    \Alambdaonelambdatwo(u,y) = \left\{
        \begin{aligned}
            &\frac{1}{2}\lvert u - 2\rvert^2 + \frac{1}{2}\lvert y\rvert^2 + \frac{\lambda_1}{2} \lvert 2u - y\rvert^2 + \omega \left[\frac{1}{2}\lvert u - 2\rvert^2 + \frac{1}{2}\lvert y\rvert^2 - J(u^{\lambda_2}, y^{\lambda_2})\right]^2, &\quad \text{if } (u,y) \in \Omega_1,\\
            &\frac{1}{2}\lvert u - 2\rvert^2 + \frac{1}{2}\lvert y\rvert^2 + \frac{\lambda_1}{2} \lvert 2u - y\rvert^2, &\quad \text{if } (u,y) \in \Omega_2.
        \end{aligned}
    \right.
\end{equation}
where 
\begin{equation*}
    J(\ulambdatwo, \ylambdatwo) = \frac{1}{2}\lvert \ulambdatwo - 2\rvert^2 + \frac{1}{2}\lvert \ylambdatwo\rvert^2
\end{equation*}
is a fixed number that can be calculated using \eqref{eq:analytical_solution_lambda} once $\lambda_2$ is fixed. 

\begin{figure}[htp]
    \centering
    \begin{minipage}{0.36\textwidth}
        \centering
        \includegraphics[width=\textwidth]{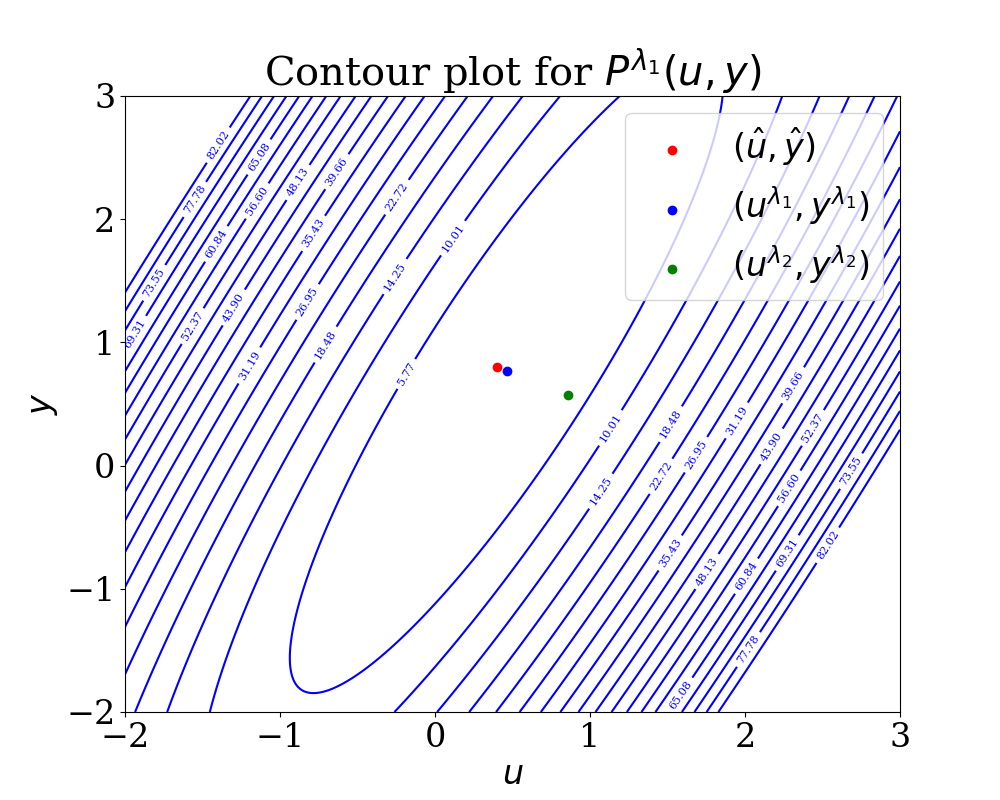}
        \subcaption{Contour plot for $\Plambda(u,y)$ with $\lambda=5$}
        \label{fig:lam1}
    \end{minipage}
    \hspace{0.05\textwidth}
    \begin{minipage}{0.36\textwidth}
        \centering
        \includegraphics[width=\textwidth]{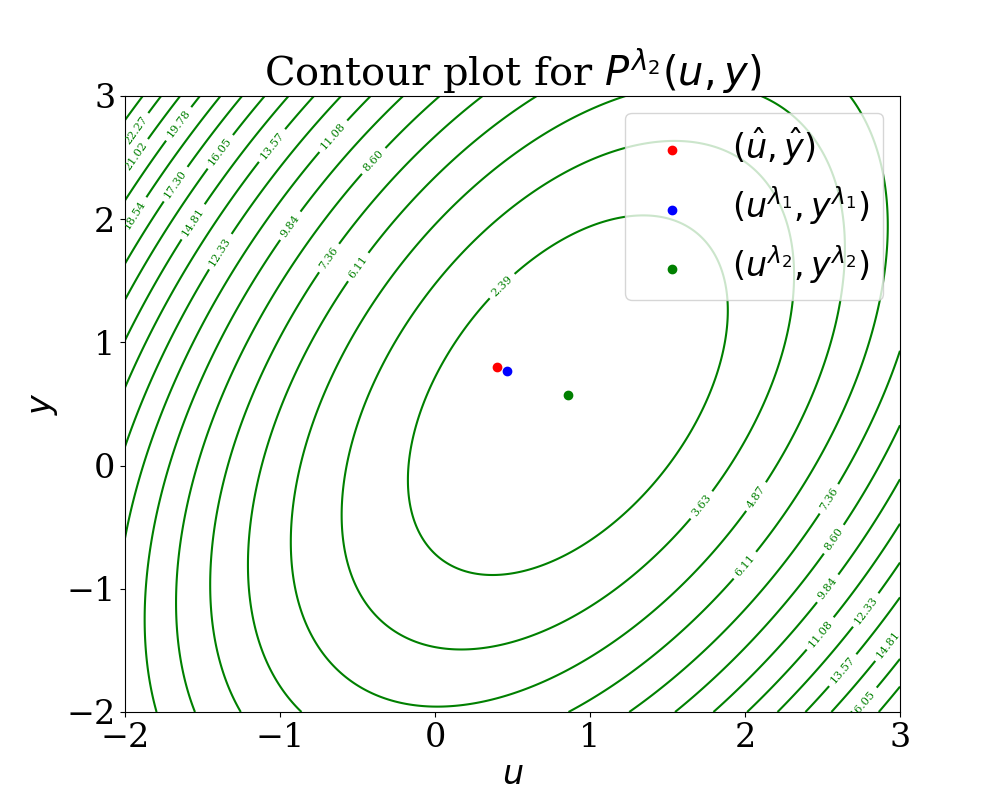}
        \subcaption{Contour plot for $\Plambda(u,y)$ with $\lambda=0.5$}
        \label{fig:lam2}
    \end{minipage}
    
    \vspace{0.5cm} 
    
    \begin{minipage}{0.32\textwidth}
        \centering
        \includegraphics[width=\textwidth]{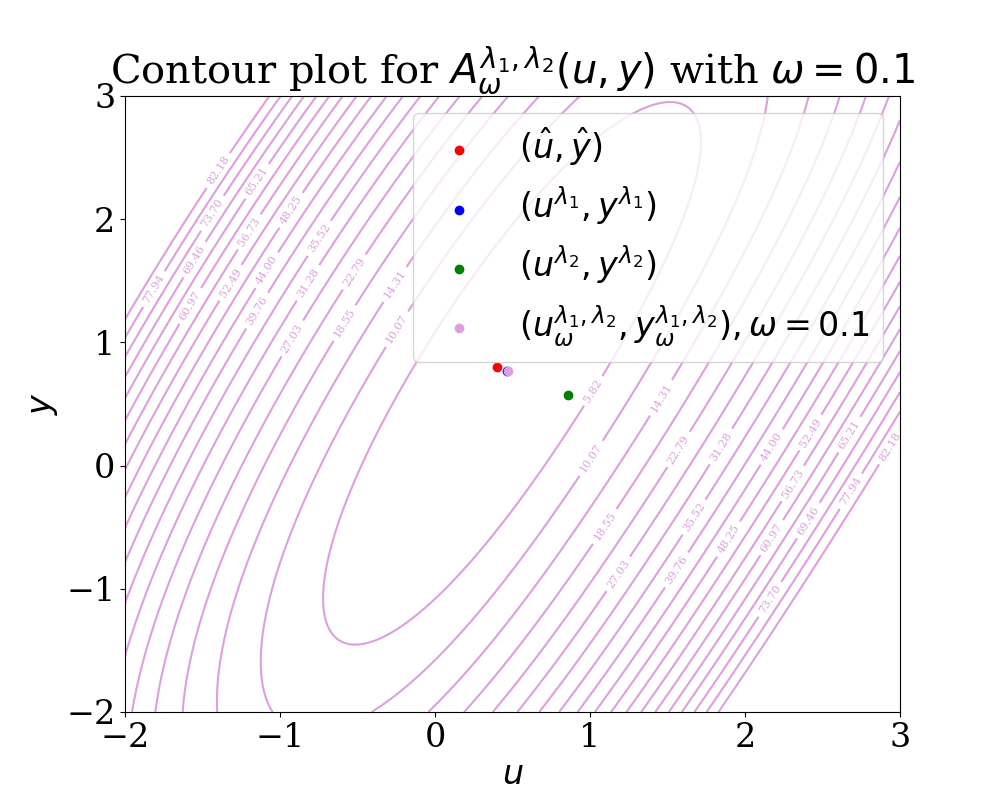}
        \subcaption{Contour plot for $\Alambdaonelambdatwo(u,y)$ with $\lambda_1=5, \lambda_2=0.5, \omega = 0.1$}
        \label{fig:T_omega_0.1}
    \end{minipage}
    \hspace{0.001\textwidth}
    \begin{minipage}{0.32\textwidth}
        \centering
        \includegraphics[width=\textwidth]{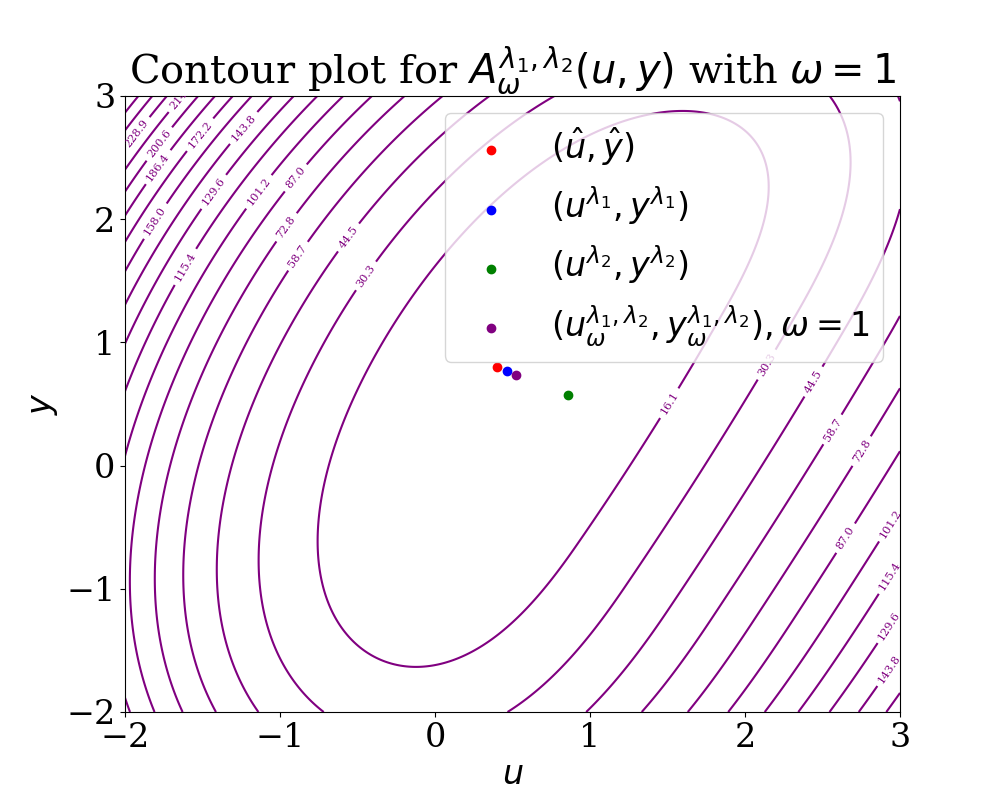}
        \subcaption{Contour plot for $\Alambdaonelambdatwo(u,y)$ with $\lambda_1=5, \lambda_2=0.5, \omega = 1$}
        \label{fig:T_omega_1}
    \end{minipage}
    \hspace{0.001\textwidth}
    \begin{minipage}{0.32\textwidth}
        \centering
        \includegraphics[width=\textwidth]{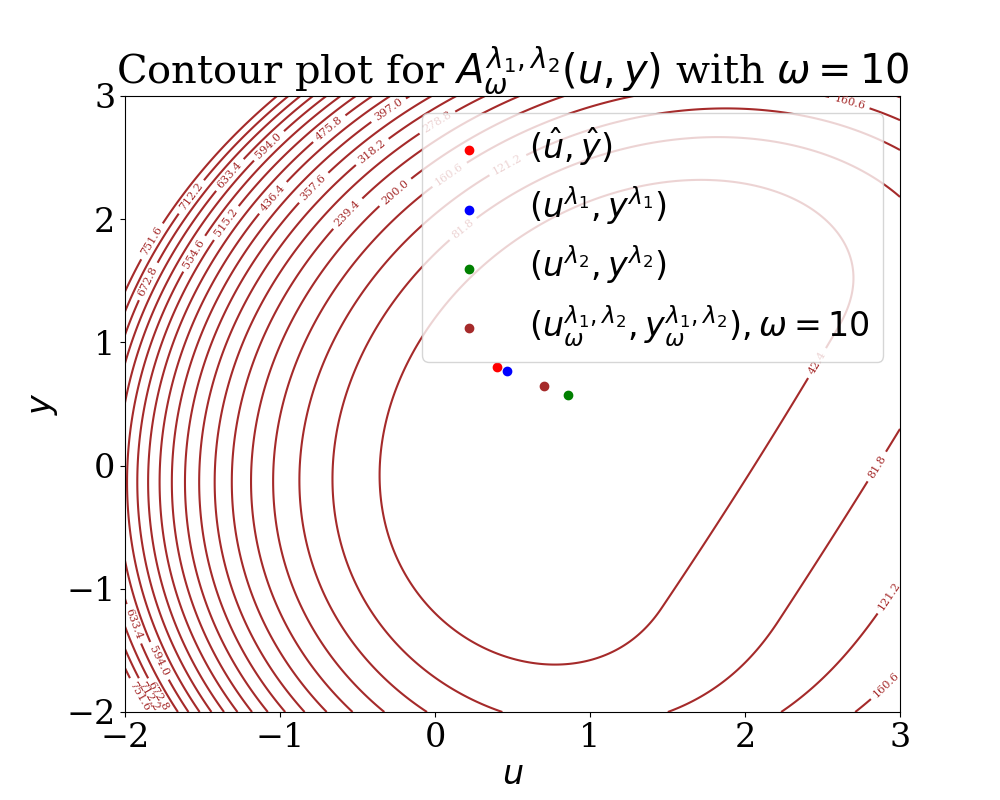}
        \subcaption{Contour plot for $\Alambdaonelambdatwo(u,y)$ with $\lambda_1=5, \lambda_2=0.5, \omega = 10$}
        \label{fig:T_omega_10}
    \end{minipage}

    \caption{Contour plots of the functions: (a) $\Plambda(u,y)$ with $\lambda=5$, (b) $\Plambda(u,y)$ with $\lambda=0.5$, (c) $\Alambdaonelambdatwo(u,y)$ with $\lambda_1=5, \lambda_2=0.5, \omega = 0.1$, (d) $\Alambdaonelambdatwo(u,y)$ with $\lambda_1=5, \lambda_2=0.5, \omega = 1$, (e) $\Alambdaonelambdatwo(u,y)$ with $\lambda_1=5, \lambda_2=0.5, \omega = 10$}
    \label{fig:contour_plots}
\end{figure}

We plot the level sets of $\Plambdaone(u,y)$ and $\Alambdaonelambdatwo(u,y)$ with different parameters in Figure \ref{fig:contour_plots}. To provide more information, we extend the range of each figure to $[-2, 3]$ instead of $[-0.5, 2]$ as in Figure \ref{fig:contour_different_lambda}. As shown in Figure \ref{fig:contour_plots}, the first row displays the contour plots of $\Plambda(u,y)$ with $\lambda_1 = 5$ and $\lambda_2 = 0.5$, while the second row displays the contour plots of $\Alambdaonelambdatwo(u,y)$ with the same values for $\lambda_1$ and $\lambda_2$, and three different choices of $\omega$: 0.1, 1, and 10. The points $(\hat{u}, \hat{y})$, $(\ulambdaone, \ylambdaone)$, and $(\ulambdatwo, \ylambdatwo)$ are marked as red, blue, and green points, respectively, in each figure. In the second row, a fourth point representing the minimizer of $\Alambdaonelambdatwo(u,y)$ for the corresponding value of $\omega$ is also marked.

We make the following observations:

\begin{enumerate}
    \item The point $(\ulambdaonelambdatwoomega, \ylambdaonelambdatwoomega)$ always lies between $(\ulambdaone, \ylambdaone)$ and $(\ulambdatwo, \ylambdatwo)$ and is closer to the analytical solution $(\hat{u}, \hat{y})$ than $(\ulambdatwo, \ylambdatwo)$.
    \item The smaller the value of $\omega$, the closer $(\ulambdaonelambdatwoomega, \ylambdaonelambdatwoomega)$ is to $(\ulambdaone, \ylambdaone)$. In Figure \ref{fig:contour_plots} (c), $(\ulambdaonelambdatwoomega, \ylambdaonelambdatwoomega)$ almost overlaps with $(\ulambdaone, \ylambdaone)$ for $\omega = 0.1$.
    \item As the value of $\omega$ increases, the contour shape changes from elliptical to bowl-shaped, with the direction towards $(\ulambdatwo, \ylambdatwo)$ generally expanding.
\end{enumerate}

The third point illustrates why it is easier to minimize $\Alambdaonelambdatwo(u,y)$ compared to $\Plambdaone(u,y)$. At points away from the minimizer, for a suitable choice of $\omega$, the condition number is more favorable since the contour of the corresponding level set is not as flattened as it is for $\Plambdaone(u,y)$. This behavior enables the effective use of standard iterative methods to find the minimizer.

\subsubsection{Relation to Auto-Tuning Parameter Strategy}\label{subsubsec:equivalence_auto_tuning}
To address the central challenge discussed in this paper regarding the penalty method—specifically, the difficulty of selecting a suitable parameter, as small parameters lead to inaccurate solutions and large parameters lead to ill-conditioned problems—a common strategy is to propose an approach with varying penalty parameters. This strategy is also commonly used when employing neural networks to solve PDE-related problems; for example, see \cite{lu2021physics}. 

To illustrate this, consider using $\lambda^k$ to denote the penalty parameter used in the $k$-th iteration, and as an example given in \cite{lu2021physics}, one possible choice is to take $  \lambda^{k+1}=\beta^k \lambda^k$,
where $\beta^k > 0$ is a constant parameter selected for the $k$-th iteration. The main difficulty with this strategy is that the choices of $\{\beta^k\}_k$ are problem-dependent \cite{bertsekas2014constrained}, often requiring fine-tuning in practice. Such a process is usually challenging or highly technical. It would be favorable if there were an automatic strategy to tune the value of $\lambda^k$. Here, we point out that the penalty adversarial strategy proposed here provides such an automatic mechanism. Specifically, when $J(u, y) \geq J(u^{\lambda_2}, y^{\lambda_2})$, we have
\begin{equation*}
\begin{aligned}
    \Alambdaonelambdatwo(u,y) &= J(u,y) + \frac{\lambda_1}{2} R(u,y) + \omega \left[J(u, y) - J(u^{\lambda_2}, y^{\lambda_2})\right]^2,
\end{aligned}
\end{equation*}
and its gradient can be computed as
\begin{equation}
    \label{eq:gradient_A_lambdaone_lambdatwo}
    \nabla \Alambdaonelambdatwo(u,y)=\left(1+2\omega\left(J(u,y)-J(\ulambdatwo,\ylambdatwo)\right)\right)\nabla J(u,y)+\frac{\lambda_1}{2}\nabla R(u,y).
\end{equation}
On the other hand, for $\Plambda(u,y)$ with a general $\lambda>0$, its gradient can be computed as
\begin{equation}
    \nabla \Plambda(u,y) = \nabla J(u,y)+\frac{\lambda}{2}\nabla R(u,y).
\end{equation}
If one is using an iterative method to find the minimizer and employs an explicit scheme, then minimizing $\Alambdaonelambdatwo(u,y)$ at a fixed stage $(u,y)$ is equivalent to minimizing $P^{\tilde{\lambda}}(u,y)$ in the sense that their gradients have the same direction, with 
\begin{equation}\label{eq:tilde_lambda}
    \tilde{\lambda} = \frac{\lambda_1}{1 + 2\omega \left(J(u,y) - J(u^{\lambda_2}, y^{\lambda_2})\right)} < \lambda_1.
\end{equation}
Thus, the penalty adversarial problem becomes more accessible to minimize under such conditions, as it is equivalent to using a smaller penalty parameter than $\lambda_1$ when $J(u,y)>J(\ulambdatwo,\ylambdatwo)$. In addition, as $J(u,y)$ gets closer to $J(\ulambdatwo,\ylambdatwo)$, $\tilde{\lambda}$ will approach $\lambda_1$.

To conclude, solving the penalty adversarial problem is a strategy to dynamically tune the penalty parameter according to the closeness of the solution to the real solution. When the solution is far from the correct one, the corresponding problem will be equivalent to minimizing $P^\lambda(u,y)$ with a small $\lambda$, allowing it to converge to the actual solution more quickly. As the solution nears the correct solution, the problem transitions to one with a larger penalty term, thereby adhering more closely to the constraints. Compared to other traditional penalty methods with dynamic penalty terms, this approach is more straightforward to implement since no specific strategy is needed to determine the dynamics of the penalty parameters in advance.

To this end, we use this equivalence to provide an intuitive reasoning for the existence of an upper bound for the choice of $\omega$. As seen from \eqref{eq:tilde_lambda}, if $\omega$ is too large, when $J(u,y)>J(\ulambdatwo,\ylambdatwo)$, $\tilde{\lambda}$ can become very small, even smaller than $\lambda_2$. In this case, minimizing $\Alambdaonelambdatwo(u,y)$ behaves similarly to minimizing a penalty problem with a penalty parameter smaller than $\lambda_2$. Consequently, the assertion that $R(\ulambdaonelambdatwoomega,\ylambdaonelambdatwoomega) < R(\ulambdatwo,\ylambdatwo)$ will no longer hold.

\subsection{Alternate Formulations of Penalty Adversarial Problem}

As implemented in \eqref{eq:penalty_adversarial_linearized_problem}, the additional penalty term associated with $\omega$ is chosen as a quadratic form of the difference $\lvert J(u,y) - J(\ulambdatwo, \ylambdatwo) \rvert$. However, this choice is not unique. Other non-negative functions of $\lvert J(u,y) - J(\ulambdatwo, \ylambdatwo) \rvert$ can also be used to add this penalty term. A natural choice is to use a monomial with a power of $k$, and the corresponding function can be defined as:

\begin{equation}
    \label{eq:penalty_adversarial_linearized_problem_k}
    \Alambdaonelambdatwok(u,y) = \left\{
        \begin{aligned}
            &J(u,y) + \frac{\lambda_1}{2} R(u,y) + \omega \left\lvert J(u, y) - J(u^{\lambda_2}, y^{\lambda_2})\right\rvert^k, &\quad \text{if } (u,y) \in \Omega_1,\\
            &J(u,y) + \frac{\lambda_1}{2} R(u,y), &\quad \text{if } (u,y) \in \Omega_2.
        \end{aligned}
    \right.
\end{equation}

\begin{figure}[htp]
    \centering
    \begin{minipage}{0.31\textwidth}
        \centering
        \includegraphics[width=\textwidth]{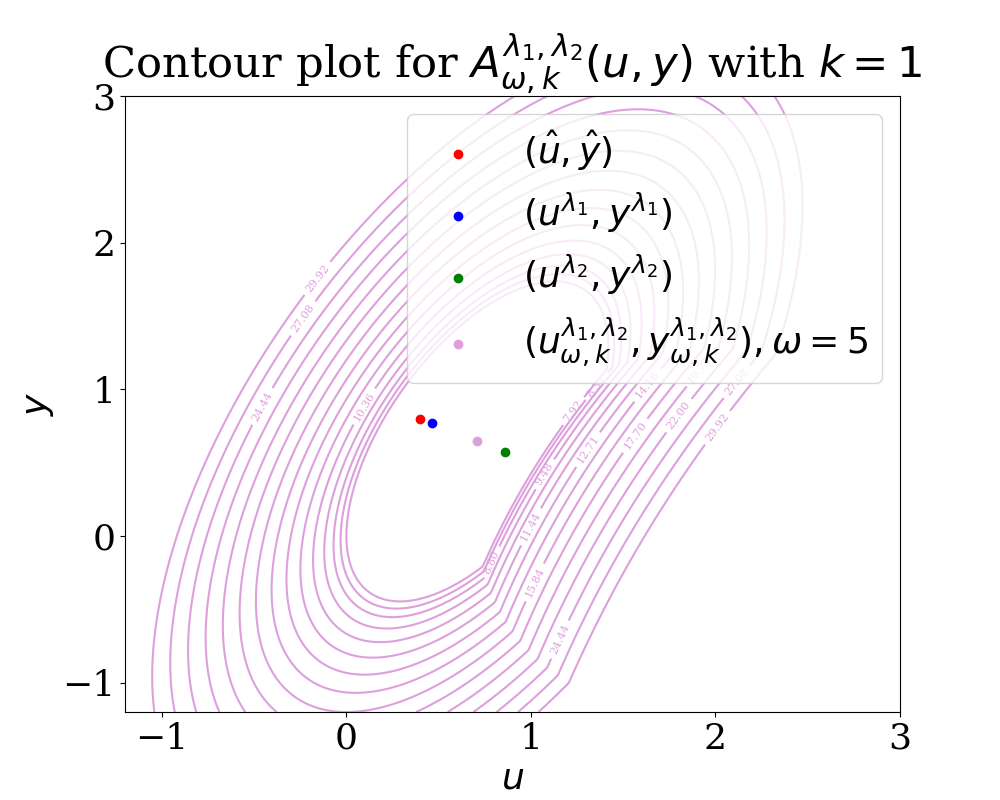}
        \subcaption{Contour plot for $\Alambdaonelambdatwok(u,y)$ with $k=1$}
        \label{fig:T_k_1}
    \end{minipage}
    \hspace{0.02\textwidth}
    \begin{minipage}{0.31\textwidth}
        \centering
        \includegraphics[width=\textwidth]{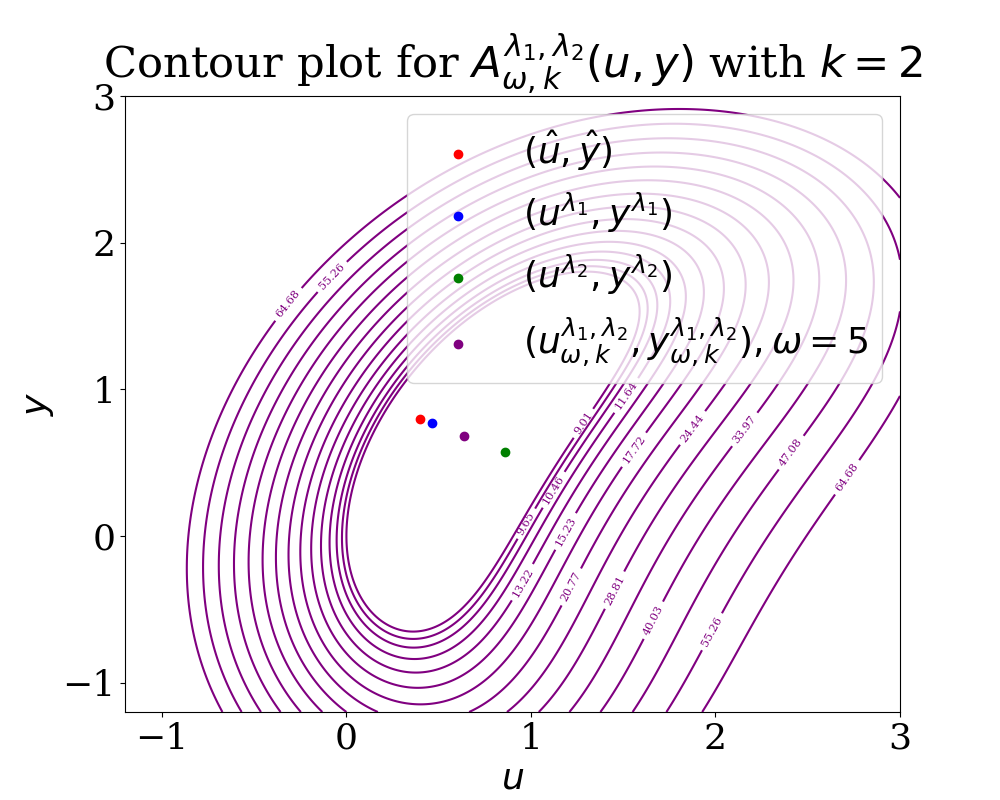}
        \subcaption{Contour plot for $\Alambdaonelambdatwok(u,y)$ with $k=2$}
        \label{fig:T_k_2}
    \end{minipage}
    \hspace{0.02\textwidth}
    \begin{minipage}{0.31\textwidth}
        \centering
        \includegraphics[width=\textwidth]{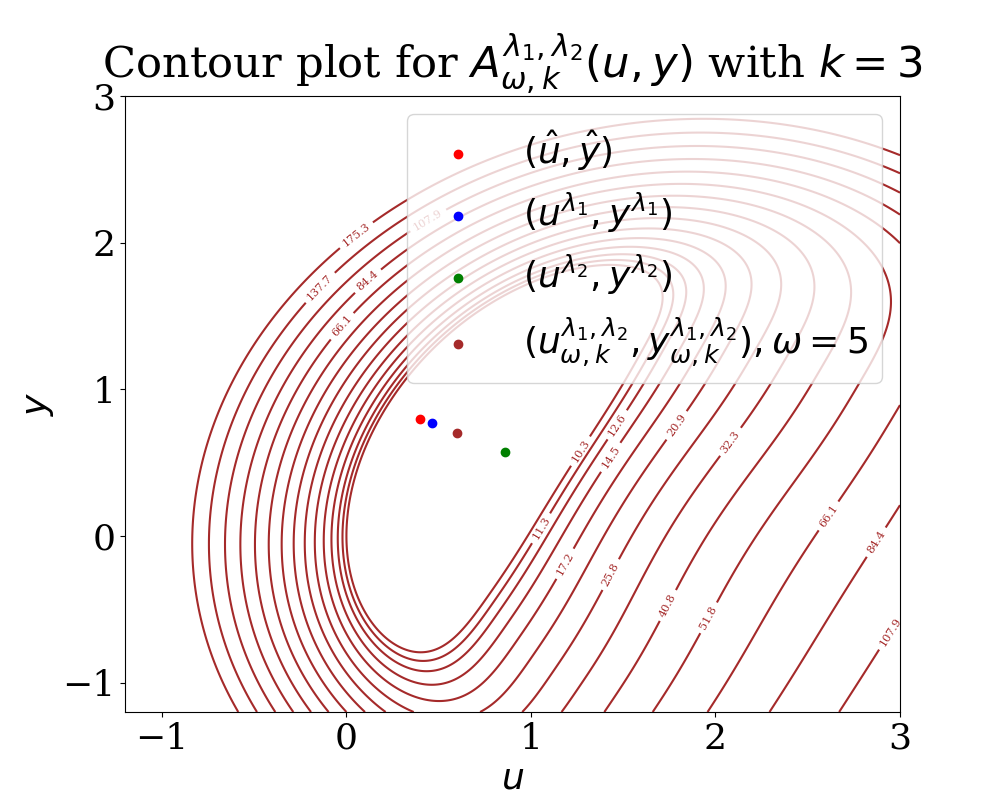}
        \subcaption{Contour plot for $\Alambdaonelambdatwok(u,y)$ with $k=3$}
        \label{fig:T_k_3}
    \end{minipage}

    \vspace{0.5cm} 

    \begin{minipage}{0.31\textwidth}
        \centering
        \includegraphics[width=\textwidth]{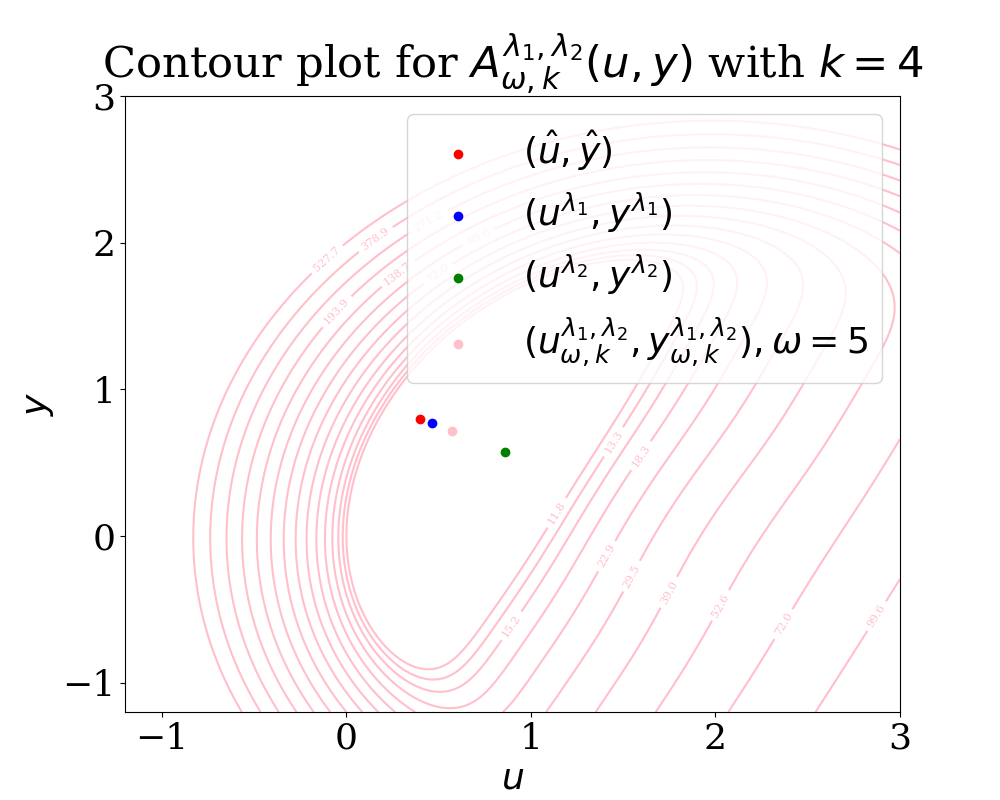}
        \subcaption{Contour plot for $\Alambdaonelambdatwok(u,y)$ with $k=4$}
        \label{fig:T_k_4}
    \end{minipage}
    \hspace{0.02\textwidth}
    \begin{minipage}{0.31\textwidth}
        \centering
        \includegraphics[width=\textwidth]{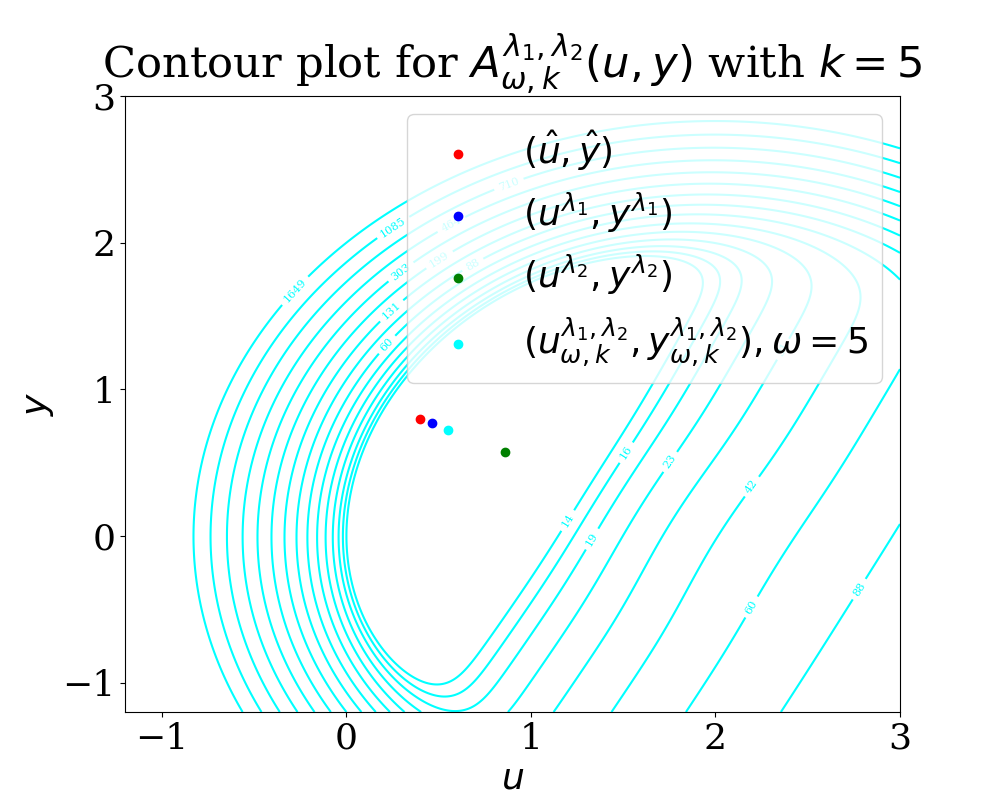}
        \subcaption{Contour plot for $\Alambdaonelambdatwok(u,y)$ with $k=5$}
        \label{fig:T_k_5}
    \end{minipage}
    \hspace{0.02\textwidth}
    \begin{minipage}{0.31\textwidth}
        \centering
        \includegraphics[width=\textwidth]{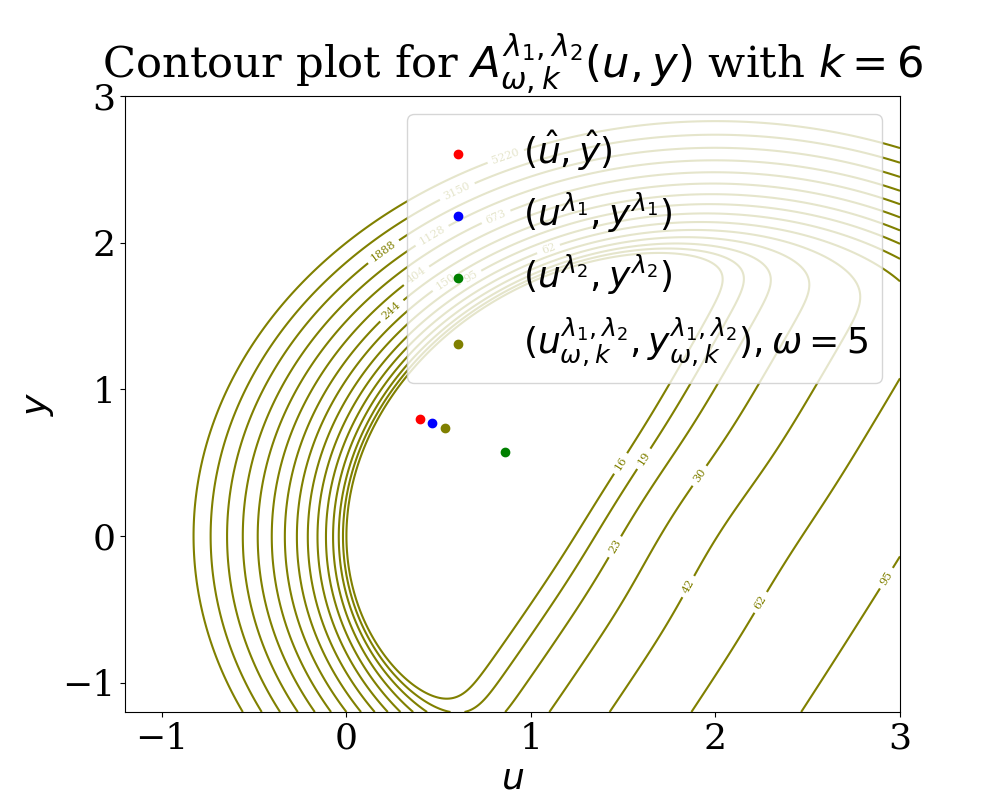}
        \subcaption{Contour plot for $\Alambdaonelambdatwok(u,y)$ with $k=6$}
        \label{fig:T_k_6}
    \end{minipage}

    \vspace{0.5cm} 

    \begin{minipage}{0.31\textwidth}
        \centering
        \includegraphics[width=\textwidth]{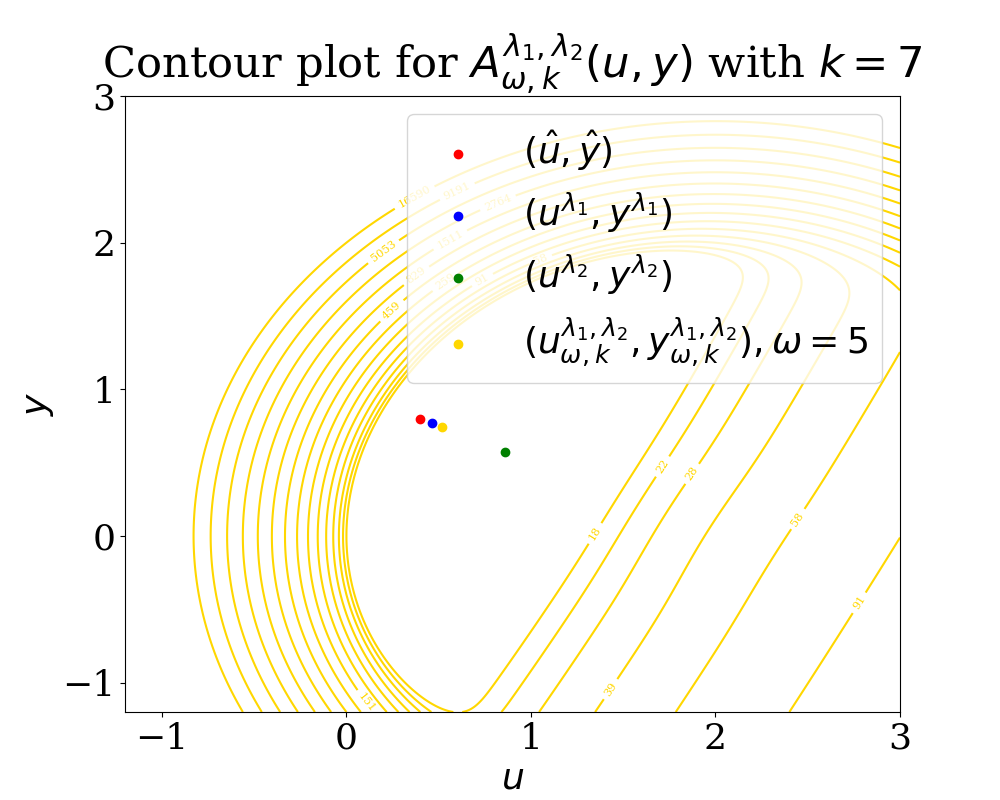}
        \subcaption{Contour plot for $\Alambdaonelambdatwok(u,y)$ with $k=7$}
        \label{fig:T_k_7}
    \end{minipage}
    \hspace{0.01\textwidth}
    \begin{minipage}{0.31\textwidth}
        \centering
        \includegraphics[width=\textwidth]{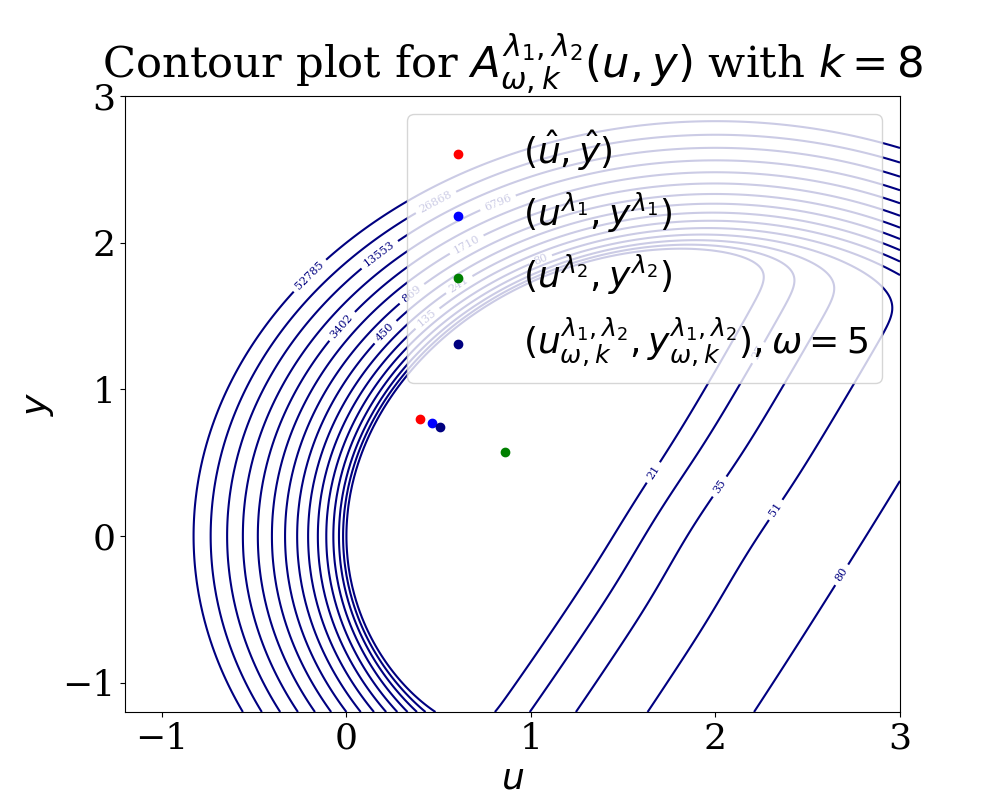}
        \subcaption{Contour plot for $\Alambdaonelambdatwok(u,y)$ with $k=8$}
        \label{fig:T_k_8}
    \end{minipage}
    \hspace{0.01\textwidth}
    \begin{minipage}{0.31\textwidth}
        \centering
        \includegraphics[width=\textwidth]{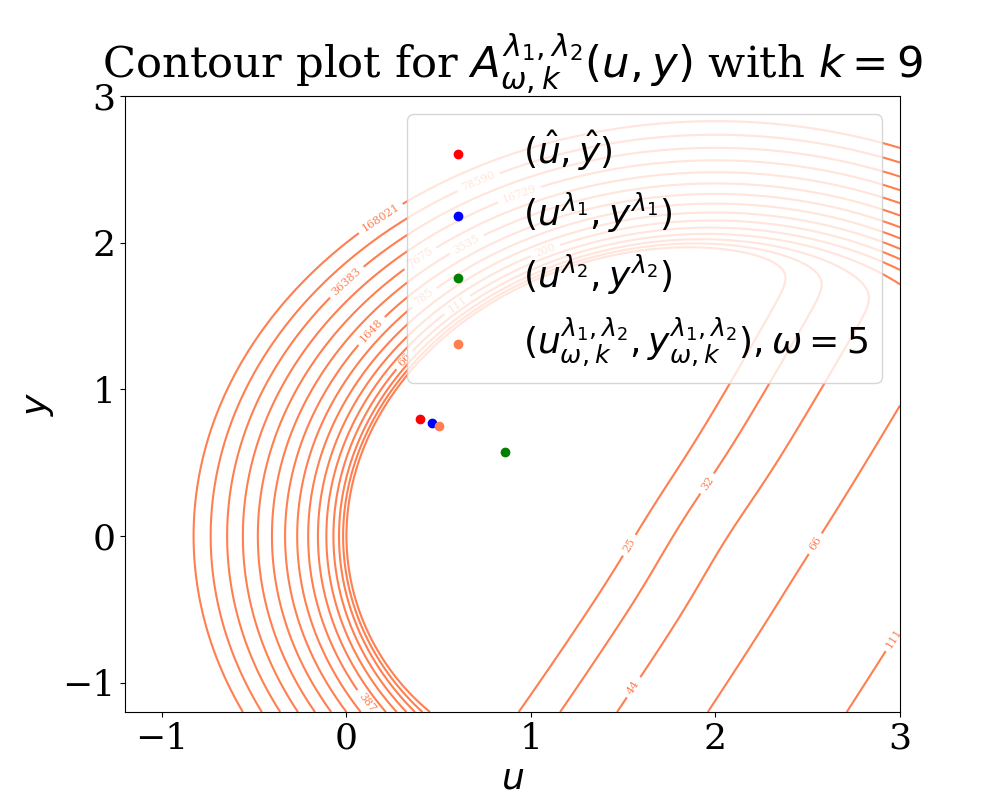}
        \subcaption{Contour plot for $\Alambdaonelambdatwok(u,y)$ with $k=9$}
        \label{fig:T_k_9}
    \end{minipage}
    \caption{Contour plots of the function $\Alambdaonelambdatwok(u,y)$ for different values of $k$: (a) $k=1$, (b) $k=2$, (c) $k=3$, (d) $k=4$, (e) $k=5$, (f) $k=6$, (g) $k=7$, (h) $k=8$, (i) $k=9$}
    \label{fig:contour_T_plots_k}
\end{figure}

The corresponding penalty adversarial problem will be modified to find the minimizer of $\Alambdaonelambdatwok(u,y)$ instead of $\Alambdaonelambdatwo(u,y) = A^{\lambda_1,\lambda_2}_{\omega,2}(u,y)$. Using the same example as presented in \eqref{eq:penalty_adversarial_linearized_problem_example}, we plot the contour of level sets for the corresponding functional $\Alambdaonelambdatwok(u,y)$ with $k$ ranging from 1 to 9, and fixed parameters: $\lambda_1 = 5$, $\lambda_2 = 0.5$, $\omega = 5$. As shown in Figure \ref{fig:contour_T_plots_k}, the corresponding minimizer, denoted as $(\ulambdaonelambdatwoomegak, \ylambdaonelambdatwoomegak)$, moves closer to the exact solution as the penalty power order $k$ increases. Additionally, as $k$ increases, the shape of the contour generally becomes more bowl-like: the curvature of the part near the actual solution remains similar, while the curvature on the other side becomes more flattened, and the curve extends to be wider. This property enhances the solvability of the problem. Moreover, we can observe that $(\ulambdaonelambdatwoomegak,\ylambdaonelambdatwoomegak)$ always lies between the point $(\ulambdaone,\ylambdaone)$ and $(\ulambdatwo,\ylambdatwo)$, and it remains closer to the actual solution than $(\ulambdatwo,\ylambdatwo)$. In conclusion, in the previous subsection, we considered the case of $k = 2$ for simplicity of analysis, but in practice, other values of $k$ could also be valid choices.

\section{Penalty Adversarial Network}\label{sec:PAN}

This section presents the neural network-based algorithm inspired by the penalty adversarial problem discussed in the previous section. Following the standard approach in physics-informed neural networks \cite{lu2021physics, mowlavi2023optimal, zang2020weak}, the input to the neural network consists of the spatial variable $x$ and the temporal variable $t$, while the output is the neural network solution to the corresponding optimal control problem constrained by PDEs, as introduced in \eqref{eq:PDE_constraints}. We denote this neural network solution as $\left(\unn(x,t;\theta), \ynn(x,t;\theta)\right)$, with parameters $\theta$ learned to minimize the corresponding loss function $L[{x}, {t}; \theta]$. By defining the loss function in different ways, various methods can be implemented to solve optimal control problems through neural networks.

As introduced above, a standard approach \cite{mowlavi2023optimal} is to add the constraints as penalty terms to the objective function and then minimize it. We denote the loss function as $L_P[{x}, {t}; \theta]$, which is defined as:

\begin{equation}\label{eq:penalty_loss_nn}
\begin{aligned}
   &\quad\,\, L_P[{x}, {t}; \theta] \\& = \frac{1}{N_J} \sum_{m=1}^{N_J} 
    \underbrace{J\left[\unn(x_m^J, t_m^J; \theta), \ynn(x_m^J, t_m^J; \theta)\right]}_{\text{Objective loss}} + \frac{\lambda_p}{N_p} \sum_{m=1}^{N_p} 
    \underbrace{\mathcal{L}\left[\unn(x_m^p, t_m^p; \theta), \ynn(x_m^p, t_m^p; \theta)\right]^2}_{\text{PDE residual loss}} \\
    &\quad\quad+ \frac{\lambda_b}{N_b} \sum_{j=1}^{N_b} 
    \underbrace{\mathcal{B}\left[\unn(x_j^b, t_j^b; \theta), \ynn(x_j^b, t_j^b; \theta)\right]^2}_{\text{Boundary loss}} + \frac{\lambda_i}{N_i} \sum_{n=1}^{N_i} 
    \underbrace{\mathcal{I}\left[\unn(x_n^i, t_n^i; \theta), \ynn(x_n^i, t_n^i; \theta)\right]^2}_{\text{Initial loss}}.
\end{aligned}
\end{equation}
The variables in the loss function are defined as follows: ${x}$ and ${t}$ represent the vectors of spatial and temporal variables used as input to the neural network, which are collections of sample points $\{x_m^J\}_{m=1}^{N_J}$, $\{x_m^p\}_{m=1}^{N_p}$, $\{x_j^b\}_{j=1}^{N_b}$, $\{x_n^i\}_{n=1}^{N_i}$, and $\{t_m^J\}_{m=1}^{N_J}$, $\{t_m^p\}_{m=1}^{N_p}$, $\{t_j^b\}_{j=1}^{N_b}$, $\{t_n^i\}_{n=1}^{N_i}$, respectively. Here, $N_J$, $N_p$, $N_b$, and $N_i$ denote the number of sample points for the objective functional, PDE residual, boundary conditions, and initial conditions, respectively. The objective loss, evaluated at sample points $(x_m^J, t_m^J)$, is denoted by $J[\unn(x_m^J, t_m^J; \theta), \ynn(x_m^J, t_m^J; \theta)]$. The PDE residual loss, boundary loss, and initial loss, evaluated at their respective sample points, are represented by $\mathcal{L}[\unn(x_m^p, t_m^p; \theta), \ynn(x_m^p, t_m^p; \theta)]$, $\mathcal{B}[\unn(x_j^b, t_j^b; \theta), \ynn(x_j^b, t_j^b; \theta)]$, and $\mathcal{I}[\unn(x_n^i, t_n^i; \theta), \ynn(x_n^i, t_n^i; \theta)]$. We note that the loss function based on the penalty formulation is slightly different from the general penalty formulation we discussed above, as we choose not to divide the penalty parameters by $2$ to keep the form of the loss functions simple. One can, however, use the previous formulations to construct the loss functions here.

In contrast to this penalty approach, we introduce an adversarial network structure enlightened by the penalty adversarial problem. Following the standard terminology used in adversarial network related works \cite{creswell2018generative}, we call the two different network solver network and discriminator network, and denote their solution respectively as $\left(\usnn(x,t;\theta),\ysnn(x,t;\theta)\right)$ and $\left(\udnn(x,t;\theta),\ydnn(x,t;\theta)\right)$. The discriminator network is simply minimizing \eqref{eq:penalty_loss_nn} with some small penalty parameters. Namely, it focuses on minimizing  
\begin{equation}\label{eq:discriminator_loss_nn}
\begin{aligned}
   &\quad\,\, L^d[{x}, {t}; \theta] \\& = \frac{1}{N_J} \sum_{m=1}^{N_J} 
   J\left[\udnn(x_m^J, t_m^J; \theta), \ydnn(x_m^J, t_m^J; \theta)\right] + \frac{\lambda_p^d}{N_p} \sum_{m=1}^{N_p} 
 \mathcal{L}\left[\udnn(x_m^p, t_m^p; \theta), \ydnn(x_m^p, t_m^p; \theta)\right]^2 \\
    &\quad\quad+ \frac{\lambda_b^d}{N_b} \sum_{j=1}^{N_b} 
  \mathcal{B}\left[\udnn(x_j^b, t_j^b; \theta), \ydnn(x_j^b, t_j^b; \theta)\right]^2 + \frac{\lambda_i^d}{N_i} \sum_{n=1}^{N_i} 
   \mathcal{I}\left[\udnn(x_n^i, t_n^i; \theta), \ydnn(x_n^i, t_n^i; \theta)\right]^2,
\end{aligned}
\end{equation}
where $\lambda_p^d, \lambda_b^d, \lambda_i^d>0$ are penalty weights for PDE loss, boundary loss and initial loss of discriminator network, respectively. In practice, we usually choose them to be relatively small to guarantee ease of training. On the other hand, for the solver network, in addition to then standard penalty formulation, we consider an extra term measuring the difference between objectives values of the discriminator network and the solver network. Namely, it minimizes the following loss function:

\begin{equation}\label{eq:solver_loss_nn}
\begin{aligned}
   &\quad\,\, L^s[{x}, {t}; \theta] \\& = \frac{1}{N_J} \sum_{m=1}^{N_J} 
   J\left[\usnn(x_m^J, t_m^J; \theta), \ysnn(x_m^J, t_m^J; \theta)\right] + \frac{\lambda_p^s}{N_p} \sum_{m=1}^{N_p} 
 \mathcal{L}\left[\usnn(x_m^p, t_m^p; \theta), \ysnn(x_m^p, t_m^p; \theta)\right]^2 \\
    &\quad\quad+ \frac{\lambda_b^s}{N_b} \sum_{j=1}^{N_b} 
  \mathcal{B}\left[\usnn(x_j^b, t_j^b; \theta), \ysnn(x_j^b, t_j^b; \theta)\right]^2 + \frac{\lambda_i^s}{N_i} \sum_{n=1}^{N_i} 
   \mathcal{I}\left[\usnn(x_n^i, t_n^i; \theta), \ysnn(x_n^i, t_n^i; \theta)\right]^2\\
   &\quad\quad+\omega \left\lvert \frac{1}{N_J} \sum_{m=1}^{N_J} J\left[\usnn(x_m^J, t_m^J; \theta), \usnn(x_m^J, t_m^J; \theta)\right]-\frac{1}{N_J} \sum_{m=1}^{N_J} J\left[\udnn(x_m^J, t_m^J; \theta), \ydnn(x_m^J, t_m^J; \theta)\right]\right\rvert^2,
\end{aligned}
\end{equation}
where $\omega$ is a tunable hyperparameter which plays the same role as it is in the penalty adversarial problem.

In practice, we initialize both networks and training data, then train these two networks together. Thanks to the fast development of machine learning toolboxes, the training of this problem is relatively standard and can be implemented directly in the TensorFlow framework \cite{abadi2016tensorflow} since it supports automatic differentiation to calculate derivatives of the loss functions with respect to the weights. Backpropagation \cite{rumelhart1986learning} is then applied to update the weights in the network. We summarize the training process in Algorithm \ref{algo:PAN_training}.

\begin{algorithm}[H]
\caption{Training of Penalty Adversarial Network \label{algo:PAN_training}}
\begin{algorithmic} 
\State \textbf{Input:} Parameters $\theta_s$, $\theta_d$, weights $\lambda_p^s$, $\lambda_p^d$, $\lambda_b^s$, $\lambda_b^d$,  $\lambda_i^s$, $\lambda_i^d$, $\omega$, learning rates $\eta_s$, $\eta_d$, maximum epochs $M$
\State \textbf{Initialize:} Solver network $\usnn(x,t;\theta_s)$, Discriminator network $\udnn(x,t;\theta_d)$, training data $\{(x_m^J, t_m^J), (x_m^p, t_m^p), (x_j^b, t_j^b), (x_n^i, t_n^i)\}$
\State Set best solver's objective loss $J^s_{\text{best}} = \infty$, best discriminator loss $L^d_{\text{best}} = \infty$
\For{epoch = 1 to $M$}
    \State \textbf{Step 1: Train Discriminator Network}
    \State Compute discriminator loss $L^d[{x}, {t}; \theta_d]$ as in \eqref{eq:discriminator_loss_nn}
    \State Update $\theta_d \leftarrow \theta_d - \eta_d \nabla_{\theta_d} L^d[{x}, {t}; \theta_d]$
    \State Store best weights if necessary
    \If{$L^d[{x}, {t}; \theta_d] < L^d_{\text{best}}$}
        \State $L^d_{\text{best}} \leftarrow L^d_P[{x}, {t}; \theta_d]$
        \State $\theta_d^{\text{best}} \leftarrow \theta_d$
        \
    \EndIf
    
    \State \textbf{Step 2: Train Solver Network}
    \State Compute solver loss $L^s[{x}, {t}; \theta_s]$ as in \eqref{eq:solver_loss_nn} and $J^s[{x}, {t}; \theta_s]$ according to definition of objective functional
    \State Update $\theta_s \leftarrow \theta_s - \eta_s \nabla_{\theta_s} L^s[{x}, {t}; \theta_s]$
    \State Store best weights if necessary
    \If{$J^s[{x}, {t}; \theta_s] < J^s_{\text{best}}$}
        \State $J^s_{\text{best}} \leftarrow J^s[{x}, {t}; \theta_s]$
        \State $\theta_s^{\text{best}} \leftarrow \theta_s$

    \EndIf
    
\EndFor

\State \textbf{Output:} Best parameters $\theta_s^{\text{best}}$, $\theta_d^{\text{best}}$
\end{algorithmic}
\end{algorithm}

To conclude this section, we note that Algorithm \ref{algo:PAN_training} is not exactly the neural network version of the penalty adversarial problem discussed in Section \ref{sec:linear_analysis}. The penalty adversarial problem uses the exact objective value at $(\ulambdatwo,\ylambdatwo)$ to compute the additional penalty term, which is a fixed number. However, in the neural network algorithm, we use $(\udnn,\ydnn)$, which is not the exact minimum but a solution to the discriminator network trained simultaneously with the solver network. To recover a neural network version of the penalty adversarial problem, it is possible to train the discriminator as a surrogate network in advance and use it to construct the penalty term. At this stage, we are not able to assert which approach is better. We present this work with the current choice for its simplicity in implementation, as it provides an all-at-once solution.

\section{Numerical Results}\label{sec:numerics}

We will present and discuss the performance of the penalty adversarial network when it is applied to solve optimal control problems constrained by different types of equations, including both linear and nonlinear problems. To begin, we will provide additional details to complement Algorithm \ref{algo:PAN_training} for practical implementation.

\subsection{Learning Rate Scheduling}\label{sec:learning_rate_strategy}

Algorithm \ref{algo:PAN_training} outlines a general workflow for simultaneously training the discriminator and solver networks. In practice, to enhance the training process, we employ a learning rate scheduling strategy \cite{darken1992learning}. Specifically, we set a certain number of epochs as the patience parameter $P$. If the loss has not improved after $P$ epochs and the learning rate is still larger than a preset minimum learning rate, the learning rate is reduced to half of its original value. This approach helps automatically fine-tune the training process by decreasing the learning rate when performance stagnates, allowing for more precise adjustments. Additionally, during training, the network often gets stuck in local minima in the first few epochs, so we typically discard the initial epochs and begin recording the best weights only after this initial phase.

\subsection{Numerical Examples}

Here, we will provide a few examples demonstrating the effectiveness of our proposed strategy when applied to various types of optimal control problems constrained by different equations. In each instance, we will present the effect of the direct penalty formulation with a large penalty parameter and then compare the results with the application of our proposed penalty adversarial approach. We will consider three different problems constrained by the 1D Poisson equation, 2D Poisson equation, and 2D Allen-Cahn equation, respectively. 


\subsubsection{Example 1: Boundary Control Problem Constrained by 1D Poisson Equation}

In this example, we consider a boundary control problem constrained by a 1D Poisson equation, following the problem setup described in \eqref{eq:boundary_optimal_control_objective}. We implement Algorithm \ref{algo:PAN_training} and the learning rate adjustment strategy discussed in Section \ref{sec:learning_rate_strategy}. The specific optimal control problem we consider is given by:

\begin{equation}
\begin{aligned}
    &\text{Minimize} && J(u) = \frac{1}{2} \int_0^1 \left[ u(x) - u_d(x) \right]^2 \, dx + \frac{ \rho}{2} \left[ u(0)^2 + u(1)^2 \right], \\
    &\text{subject to} && -\frac{d^2 u}{dx^2} = A \sin(2\pi x), \quad x \in [0,1], 
\end{aligned}
\end{equation}
where the desired state $u_d(x) = \frac{A}{4\pi^2} \sin(2\pi x) + bx + a$, with parameters $a = -10$, $b = 65$, and $A = 8\pi^2$. The control for this problem is the boundary values $u(0)$ and $u(1)$. In fact, if the boundary values are fixed, the solution to the Poisson equation is uniquely determined, which in turn determines the value of $J(u)$. 

This problem has an analytical solution, given by:
\begin{equation}\label{eq:u_analytical}
u_{\text{analytical}}(x) = \frac{A}{4 \pi^2} \sin(2 \pi x) + b^* x + a^*,
\end{equation}
where $a^* = \frac{1}{1 + 2 \rho} a + \frac{2 \rho}{(1 + 2 \rho)(1 + 6 \rho)} b$ and $b^* = \frac{b}{1 + 6 \rho}$. In this example, we take $\rho = 2$, resulting in $a^* = 2$ and $b^* = 5$.

\begin{figure}
    \centering
    \begin{subfigure}[t]{0.32\textwidth}
        \centering
        \includegraphics[width=\textwidth]{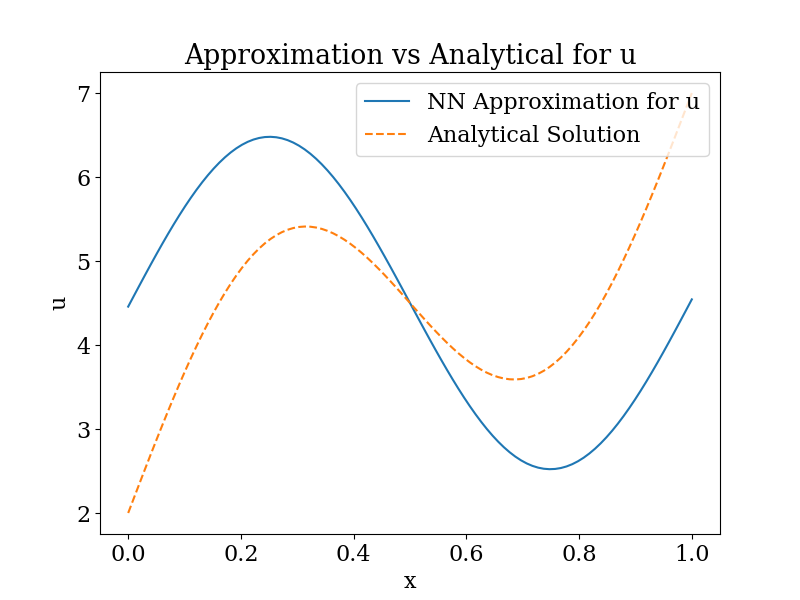}
        \subcaption{approximation vs analytical for $u$}
    \end{subfigure}
    \begin{subfigure}[t]{0.32\textwidth}
        \centering
        \includegraphics[width=\textwidth]{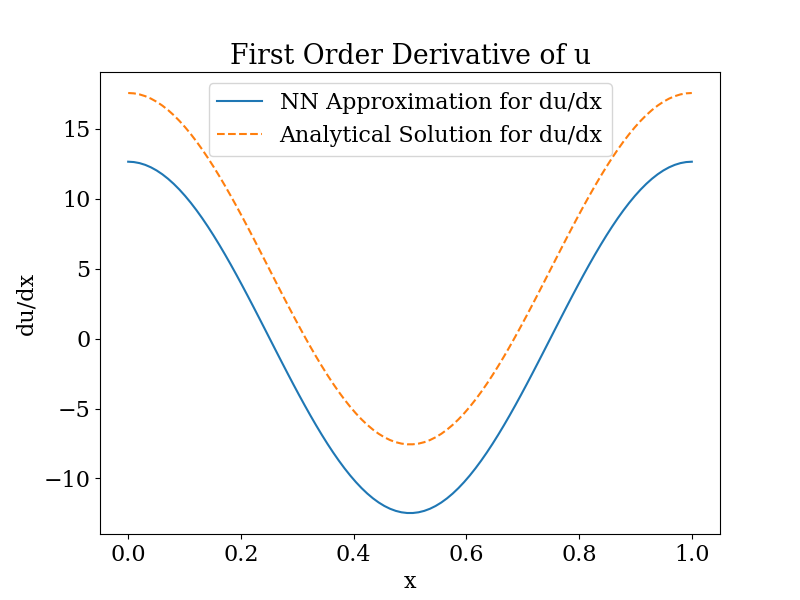}
        \subcaption{first order derivative of $u$}
    \end{subfigure}
    \begin{subfigure}[t]{0.32\textwidth}
        \centering
        \includegraphics[width=\textwidth]{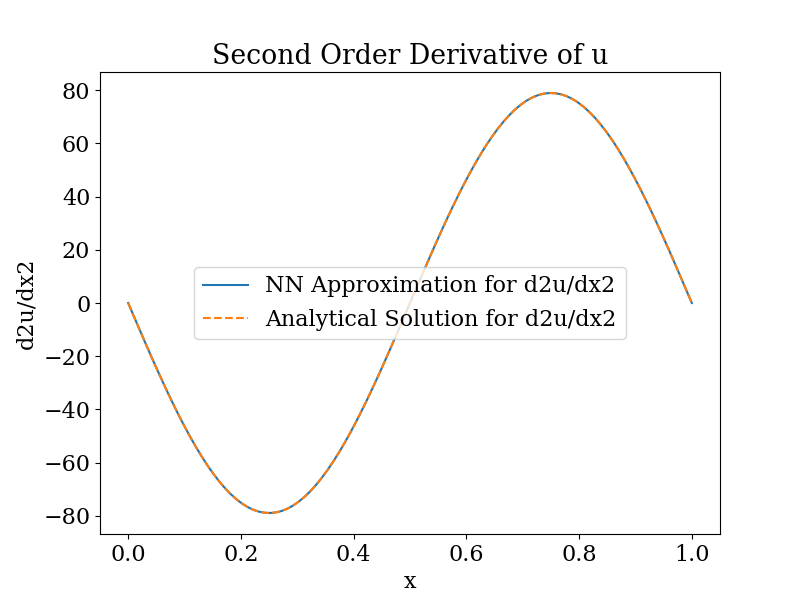}
        \subcaption{second order derivative of $u$}
    \end{subfigure}
    \begin{subfigure}[t]{0.32\textwidth}
        \centering
        \includegraphics[width=\textwidth]{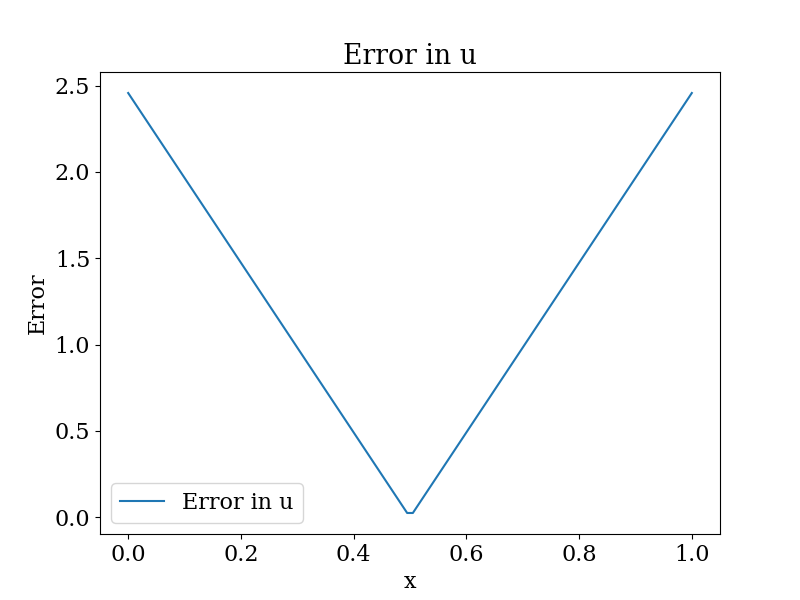}
        \subcaption{error in $u$}
    \end{subfigure}
    \begin{subfigure}[t]{0.32\textwidth}
        \centering
        \includegraphics[width=\textwidth]{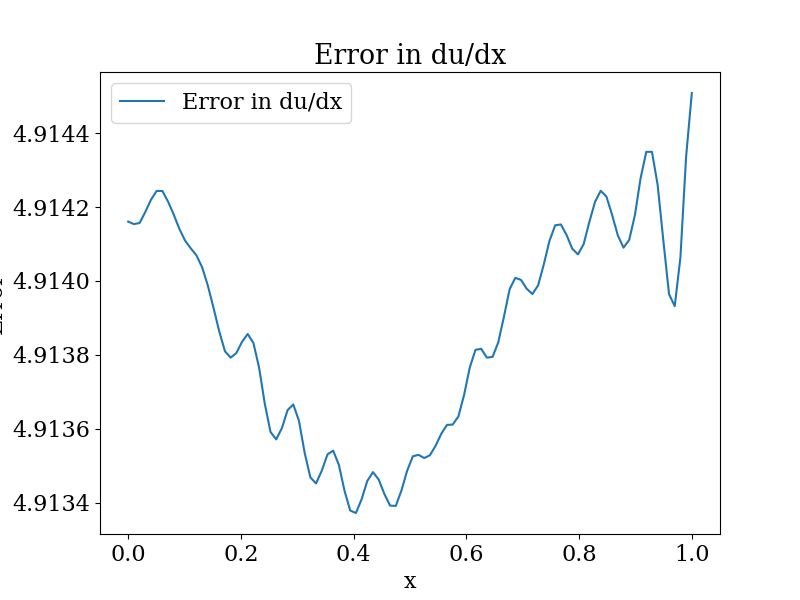}
        \subcaption{error in $\frac{du}{dx}$}
    \end{subfigure}
    \begin{subfigure}[t]{0.32\textwidth}
        \centering
        \includegraphics[width=\textwidth]{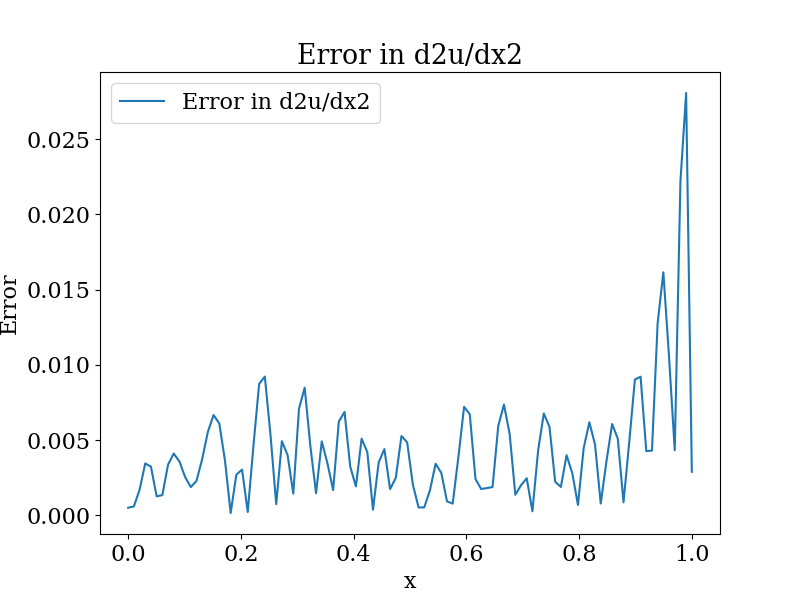}
        \subcaption{error in $\frac{d^2u}{dx^2}$}
    \end{subfigure}
    \vfill
    \begin{subfigure}[t]{0.5\textwidth}
        \centering
        \includegraphics[width=\textwidth]{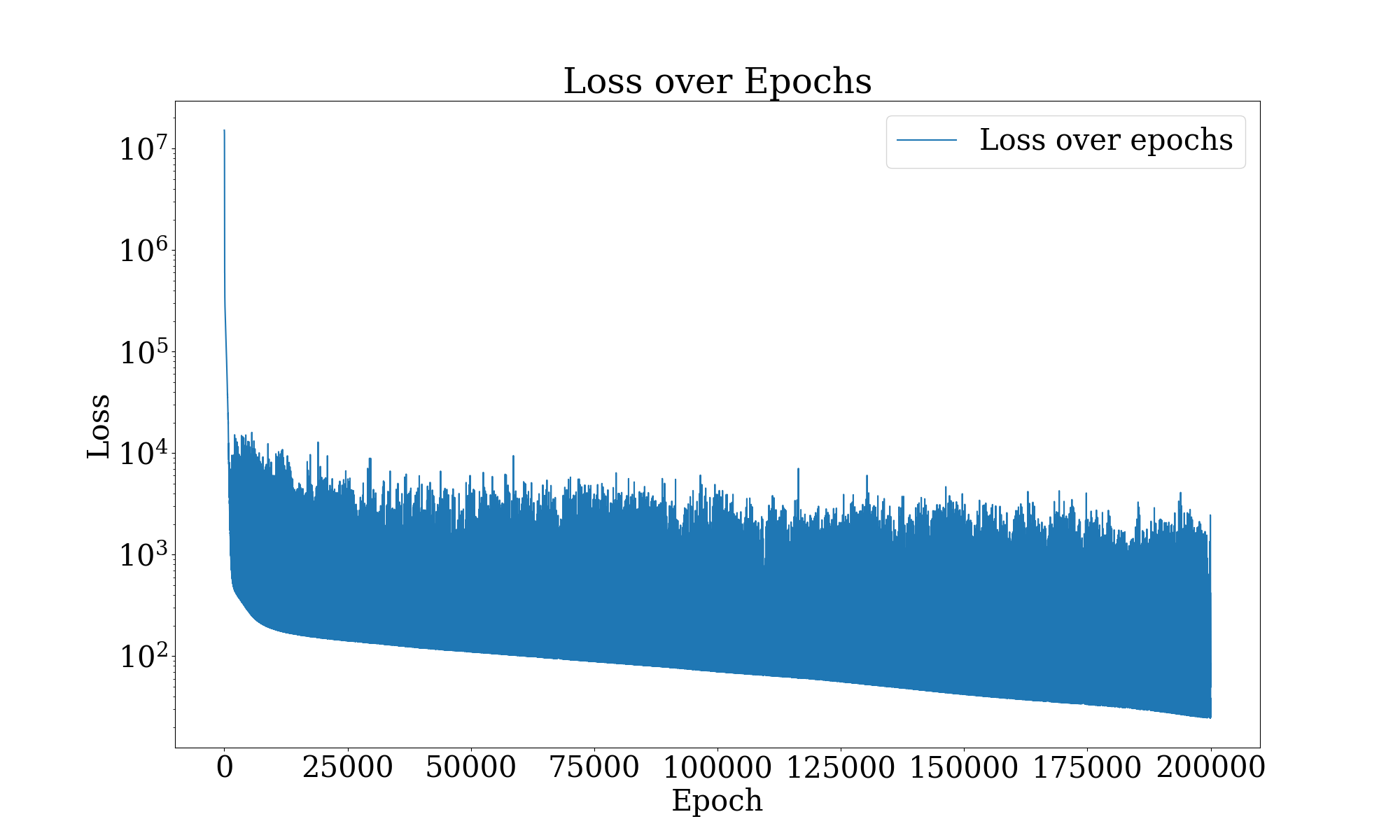}
        \subcaption{loss over epochs}
    \end{subfigure}
    \caption{Results for minimizing loss function \eqref{eq:1D_Poisson_penalty_loss_nn} with $N=32$, $\lambda_p = 5000$ and $\rho =2$. Subfigure (a) shows the approximation vs analytical solution for $u$, (b) shows the first order derivative $\frac{du}{dx}$, (c) shows the second order derivative $\frac{d^2u}{dx^2}$, (d) shows the error in $u$, (e) shows the error in $\frac{du}{dx}$, and (f) shows the error in $\frac{d^2u}{dx^2}$. Subfigure (g) presents the loss over epochs.}
        \label{fig:1D_Poisson_model}
\end{figure}

Firstly, we consider using a standard penalty formulation, as given in \eqref{eq:penalty_loss_nn}, with a large penalty parameter to solve the problem. Namely, the corresponding neural network aims to output $(\unn,\ynn)$ to minimize the following loss function:
\begin{equation}\label{eq:1D_Poisson_penalty_loss_nn}
\begin{aligned}
   L[{x}; \theta] = \frac{ \rho}{2} \left[ \unn(0;\theta)^2 + \unn(1;\theta)^2 \right]&+\frac{1}{2N} \sum_{m=1}^{N}   \left[ \unn(x_m;\theta) - u_d(x) \right]^2\\& + \frac{\lambda_p}{N} 
   \sum_{m=1}^{N} 
   \left[\frac{d^2 \unn}{dx^2}(x_m;\theta) + A \sin(2\pi x_m)\right]^2.
\end{aligned}
\end{equation}
The training data $\{x_m\}_m$ are $N = 32$ points uniformly chosen from $[0,1]$. The hyperparameters for this example chosen in implementation are as follows. We set the neural network depth to $4$, width to $40$, the penalty parameter $\lambda_p=5000$. The initial learning rate is $0.001$, which can be reduced to a minimum learning rate of $0.0001$ based on the learning rate scheduling strategy with patience set to be $3000$. The training is conducted over a maximum of $200000$ epochs. The numerical findings are presented in Figure \ref{fig:1D_Poisson_model}.

We observe that after 200000 epochs, the loss function continues to decrease. However, the numerical solution computed via the neural network remains significantly distant from the analytical solution. The derivative of the numerical solution exhibits a nearly constant error compared to the analytical derivative. In contrast, the second derivative of the numerical solution, which corresponds to the constraint in this problem, remains close to its analytical value. This outcome aligns with our expectations, as we enforce a large penalty parameter on the PDE constraint. While this enforcement ensures the solution adheres well to the constraint, it also makes it difficult for the network to be trained effectively and to reach the actual solution.

It is important to note that the training process here is standard without applying specific techniques. While refining the training procedure or increasing the number of training epochs can bring the numerical solution closer to the analytical solution, the results demonstrate that a basic neural network may struggle to find the optimal solution when a large penalty parameter is employed. In comparison, we now focus on solving the same problem using our proposed penalty adversarial strategy.

Instead of using one single network, to utilize penalty adversarial network approach, we create and train the solver and the discriminator network together. Namely, the discriminator network aims to output $(\udnn,\ydnn)$ to minimize the following loss function
\begin{equation}\label{eq:1D_Poisson_penalty_loss_nn_discriminator}
\begin{aligned}
   L^d[{x}; \theta] = \frac{ \rho}{2} \left[ \udnn(0;\theta)^2 + \udnn(1;\theta)^2 \right]&+\frac{1}{2N} \sum_{m=1}^{N}   \left[ \udnn(x_m;\theta) - u_d(x) \right]^2\\& + \frac{\lambda_p^d}{N} 
   \sum_{m=1}^{N} 
   \left[\frac{d^2 \udnn}{dx^2}(x_m;\theta) + A \sin(2\pi x_m)\right]^2,
\end{aligned}
\end{equation}
and the solver network aims to output $(\usnn,\ysnn)$ to minimize the following loss function
\begin{equation}\label{eq:1D_Poisson_penalty_loss_nn_solver}
\begin{aligned}
   L^s[{x}; \theta] =& \frac{ \rho}{2} \left[ \usnn(0;\theta)^2 + \usnn(1;\theta)^2 \right]+\frac{1}{2N} \sum_{m=1}^{N}   \left[ \usnn(x_m;\theta) - u_d(x) \right]^2\\
   &+ \frac{\lambda_p^s}{N} 
   \sum_{m=1}^{N} 
   \left[\frac{d^2 \usnn}{dx^2}(x_m;\theta) + A \sin(2\pi x_m)\right]^2\\
    &+\omega\, \bigg{\{}\frac{ \rho}{2} \left[ \usnn(0;\theta)^2 + \usnn(1;\theta)^2 \right]-\frac{ \rho}{2} \left[ \udnn(0;\theta)^2+\udnn(1;\theta)^2 \right] \\
    &\quad\,\quad\,\quad+\frac{1}{2N} \sum_{m=1}^{N}   \left[ \usnn(x_m;\theta) - u_d(x) \right]^2-\frac{1}{2N} \sum_{m=1}^{N}   \left[ \udnn(x_m;\theta) - u_d(x) \right]^2\bigg{\}}^2.
\end{aligned}
\end{equation}
The parameters here are chosen as $\lambda_p^d = 1$, $\lambda_p^s= 5000$ and $\omega = 1$.

Other hyperparameters for the construction of the networks remain the same as in the previous setting. The neural network's depth and width are set to 4 and 40, respectively. The training data consists of $N=32$ points uniformly drawn from $[0,1]$. The initial learning rate is set to $0.001$, which can be reduced to a minimum learning rate of $0.0001$ based on the learning rate scheduling strategy, with patience set to 3000. The training is conducted over a maximum of 200000 epochs.

\begin{figure}
    \centering
    \begin{subfigure}[t]{0.32\textwidth}
        \centering
        \includegraphics[width=\textwidth]{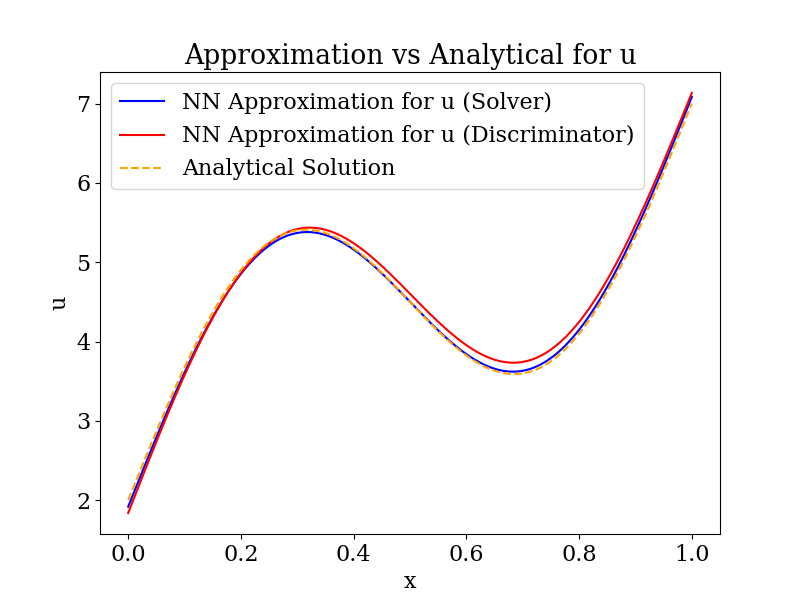}
        \subcaption{approximation vs analytical for $u$}
    \end{subfigure}
    \begin{subfigure}[t]{0.32\textwidth}
        \centering
        \includegraphics[width=\textwidth]{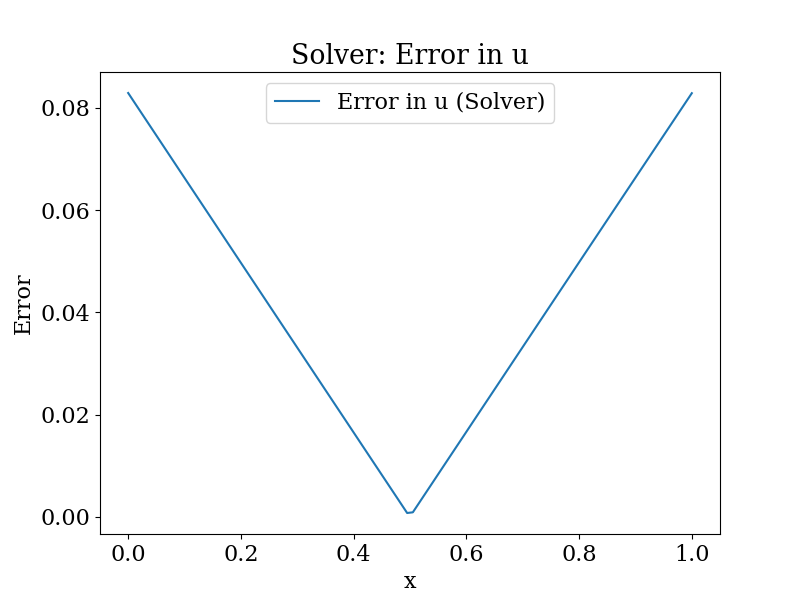}
        \subcaption{error in $u$ (solver)}
    \end{subfigure}
    \begin{subfigure}[t]{0.32\textwidth}
        \centering
        \includegraphics[width=\textwidth]{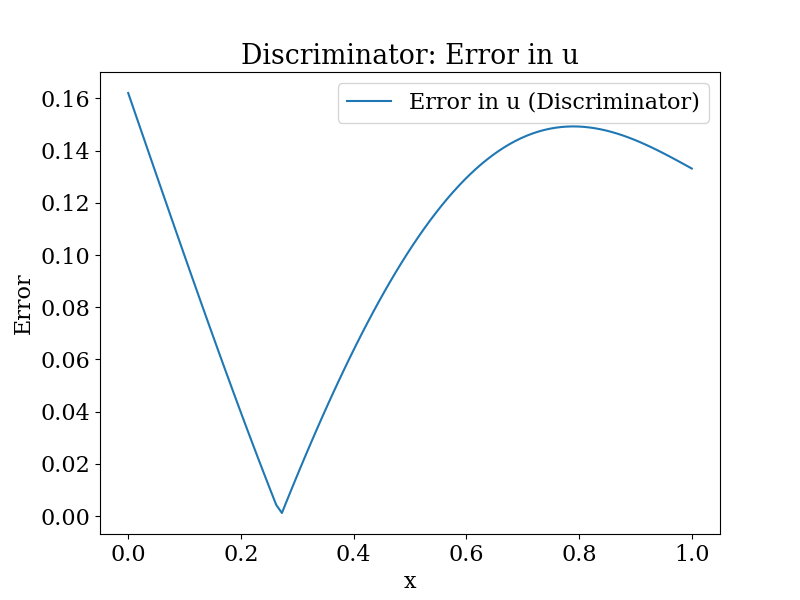}
        \subcaption{error in $u$ (discriminator)}
    \end{subfigure}
    \begin{subfigure}[t]{0.32\textwidth}
        \centering
        \includegraphics[width=\textwidth]{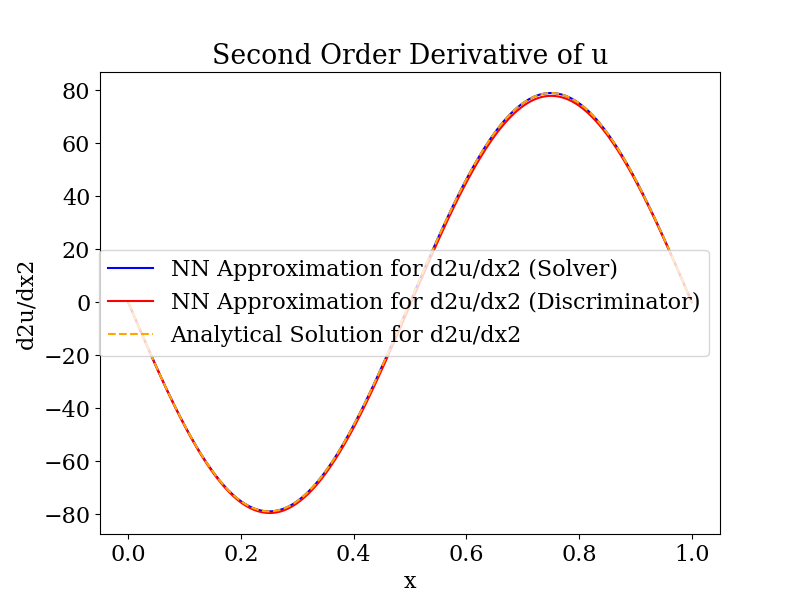}
        \subcaption{approximation vs analytical for $\frac{d^2u}{dx^2}$ }
    \end{subfigure}
    \begin{subfigure}[t]{0.32\textwidth}
        \centering
        \includegraphics[width=\textwidth]{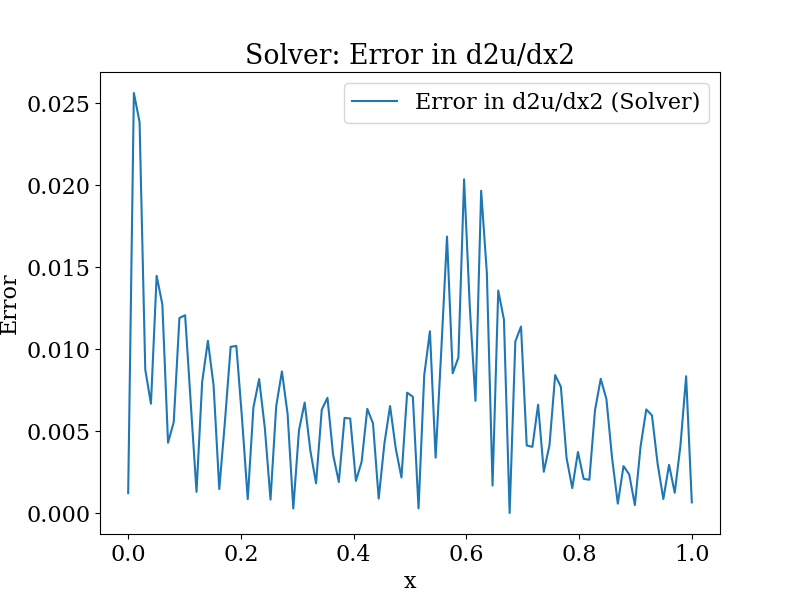}
        \subcaption{error in $\frac{d^2u}{dx^2}$ (solver)}
    \end{subfigure}
    \begin{subfigure}[t]{0.32\textwidth}
        \centering
        \includegraphics[width=\textwidth]{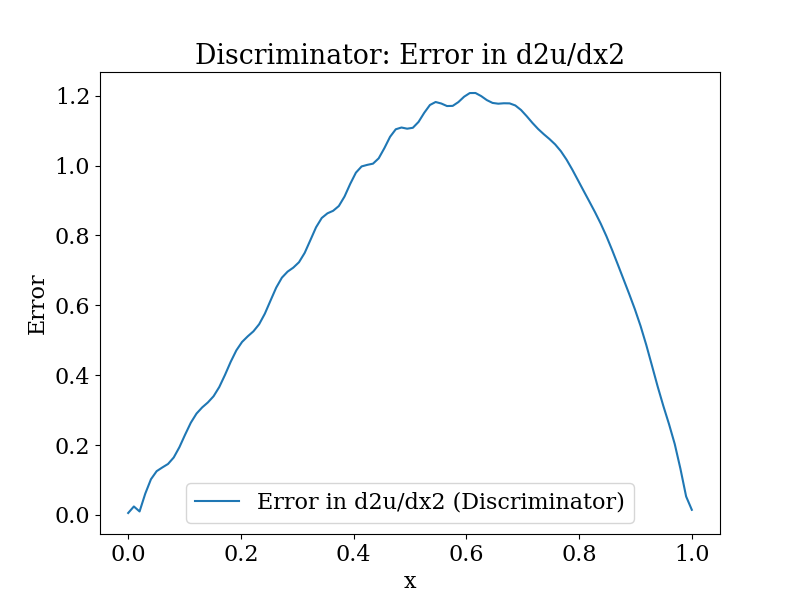}
        \subcaption{error in $\frac{d^2u}{dx^2}$ (discriminator)}
    \end{subfigure}
    \vfill
    \begin{subfigure}[t]{0.5\textwidth}
        \centering
        \includegraphics[width=\textwidth]{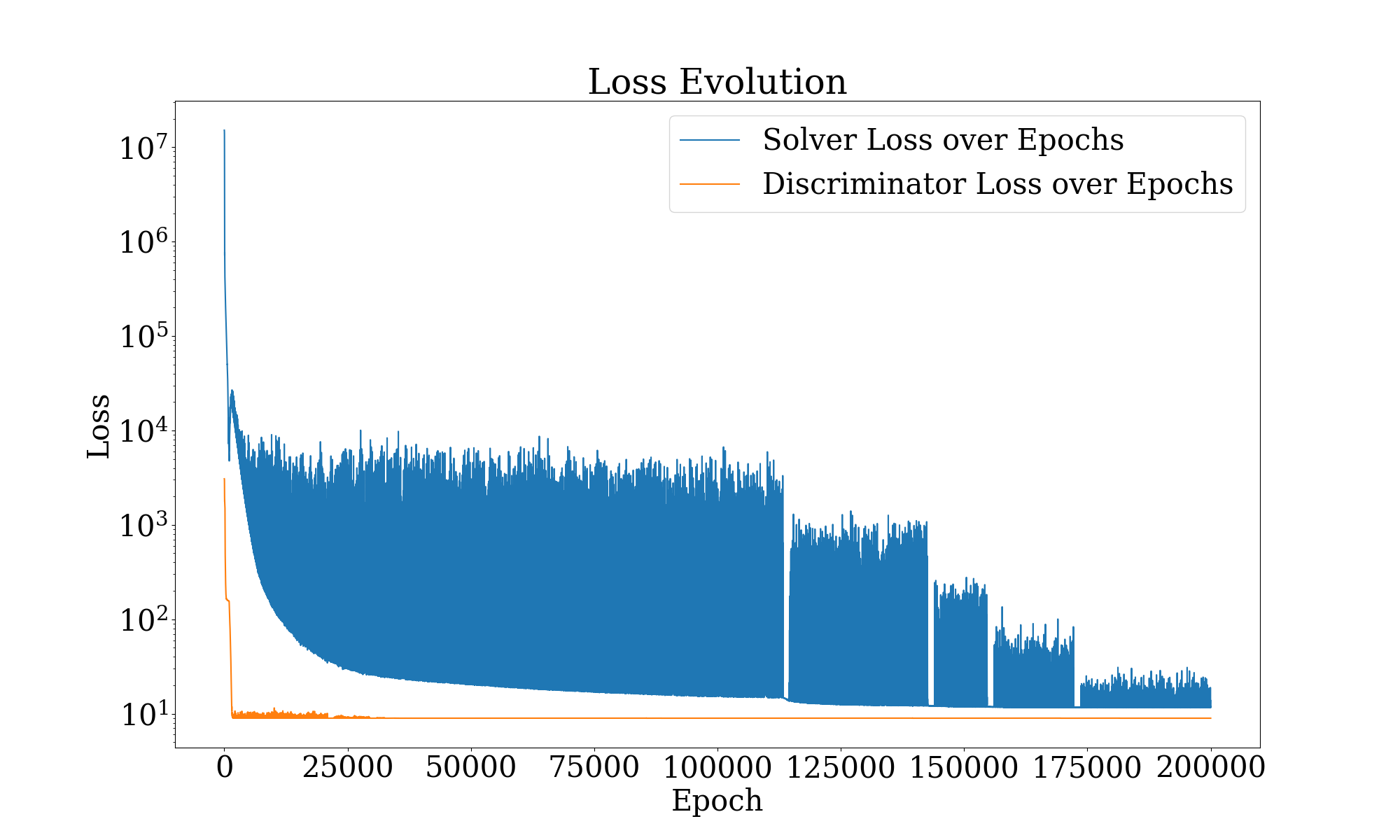}
        \subcaption{loss evolution (solver and discriminator)}
    \end{subfigure}
   \caption{Results for the penalty adversarial network approach with $\lambda_p^d=1, \lambda_p^s = 5000$, $\omega = 1$ and $\rho = 2$. The first row displays results for the solver and discriminator concerning $u$, with (a) showing the approximation versus analytical solution for $u$ (solver), (b) the error in $u$ (solver), and (c) the error in $u$ (discriminator). The second row presents the second order derivative and its error, with (d) the second order derivative of $u$ (solver), (e) error in $\frac{d^2u}{dx^2}$ (solver), and (f) error in $\frac{d^2u}{dx^2}$ (discriminator). The final image (g) details the loss evolution for both the solver and discriminator.}
    \label{fig:1D_Poisson_adversarial}
\end{figure}

The numerical findings are presented in Figure \ref{fig:1D_Poisson_adversarial}. The first row compares the neural network approximations from both the solver and discriminator networks with the analytical solution for the value of \(u\), explicitly showing the magnitude of the difference between these two numerical solutions and the analytical solution. Both numerical solutions approximate the analytical solution well. However, the solver network's solution is closer to the real solution, with a maximum error of $0.08$, while the maximum error for the discriminator network reaches $0.16$. 

A more notable difference is observed in the second row, which compares the second-order derivatives. The accuracy of the second-order derivative reflects the degree of constraint satisfaction. Here, the solver network's solution maintains the error of \(\frac{d^2u}{dx^2}\) under $0.025$, which is quite small compared to the magnitude of the second-order derivative. However, the discriminator network exhibits a much larger error in this term, reaching a maximum error of $1.2$.

This example demonstrates the advantage of using an adversarial network compared to a standard penalty formulation with both large and small penalty parameters. As shown in Figure \ref{fig:1D_Poisson_model}, using a large penalty parameter makes the neural network harder to train, resulting in a numerical solution that deviates from the real solution after a certain number of epochs. The adversarial approach produces a much more accurate approximation while maintaining a large penalty parameter to guarantee satisfying the underlying PDE. Conversely, using a small penalty parameter results in a formulation whose theoretical solution is away from the real analytical solution $u_{\text{analytical}}$, as given in \eqref{eq:u_analytical}. This fact means that even if the discriminator neural network is easier to train and converges, the obtained solution will not be accurate enough, particularly not satisfying the PDE constraint well, as shown in Figure \ref{fig:1D_Poisson_adversarial}. The adversarial approach successfully achieves a solution closer to the real analytical solution without requiring complicated training techniques beyond introducing the adversarial structure. This example highlights the effectiveness of applying the penalty adversarial approach.

\subsubsection{Example 2: Distributed Control Problem Constrained by 2D Poisson Equation}

Here we will present that our approach will also work for multi-dimensional problems. We consider a distributed control problem constrained by a 2D Poisson equation. The same problem is investigated in \cite{ghasemi2019artificial,pearson2013radial} as well.

\begin{figure}
    \centering
    
    \begin{subfigure}[t]{0.45\textwidth}
        \centering
        \includegraphics[width=\textwidth]{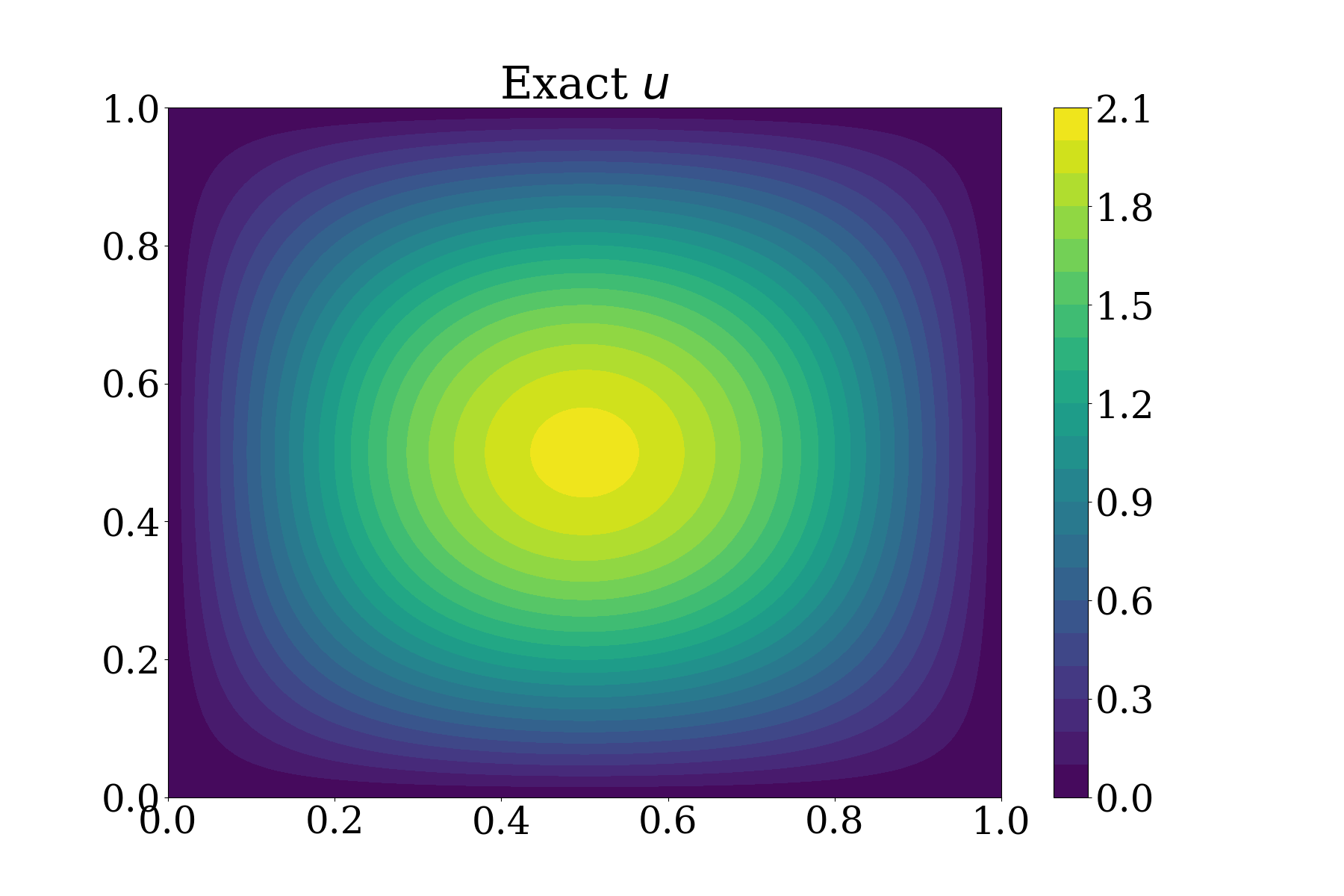}
        \subcaption{Desired state $u_d$}
    \end{subfigure}
    \begin{subfigure}[t]{0.45\textwidth}
        \centering
        \includegraphics[width=\textwidth]{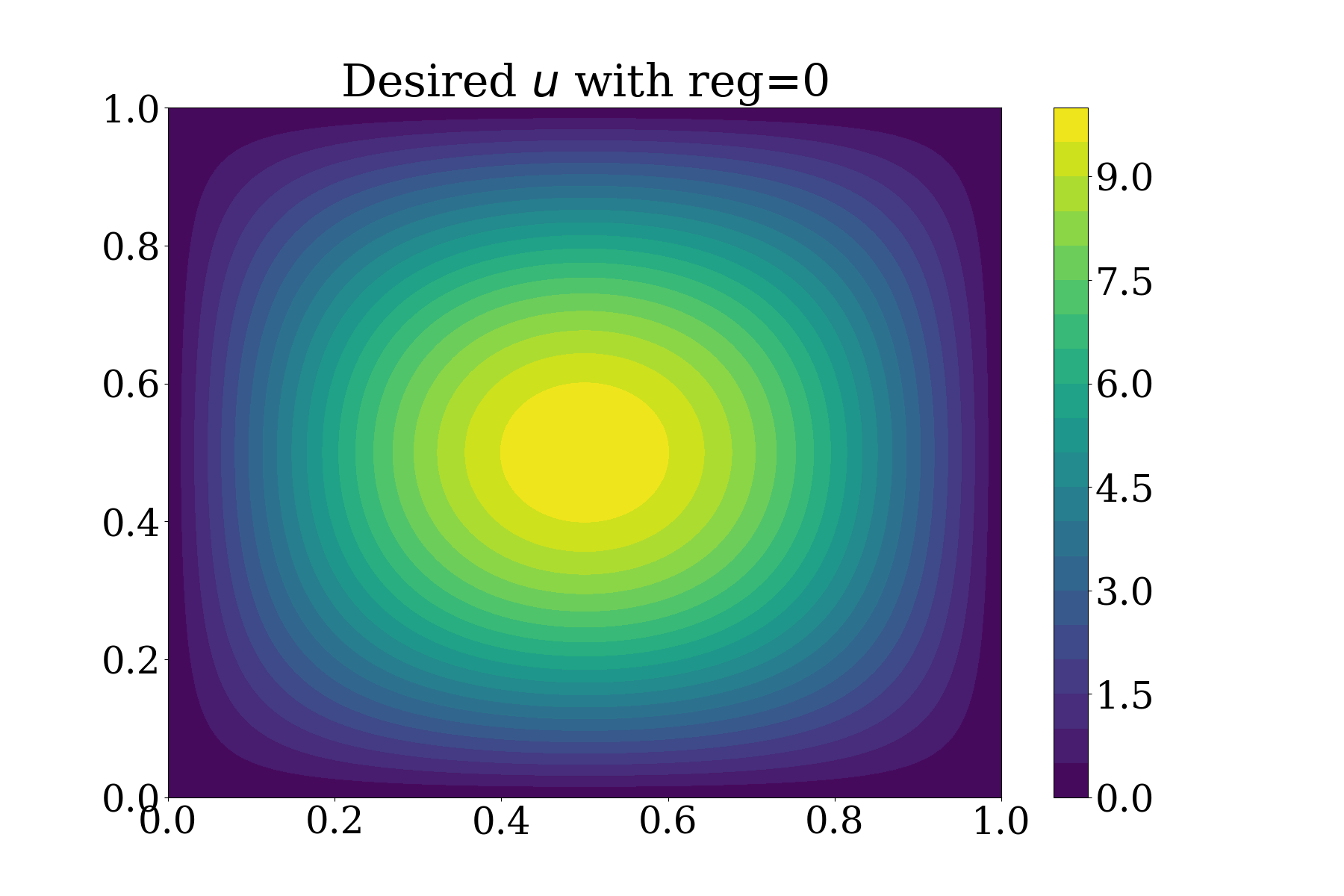}
         \subcaption{Analytical solution $u_{\text{analytical}}$}
    \end{subfigure}
    
    \caption{The desired state and analytical solution for problem \eqref{eq:2D_Poisson_control_example}}.
    \label{fig:2D_Poisson_desired_analytical}
\end{figure}

The specific optimal control problem we consider is given by:
\begin{equation}\label{eq:2D_Poisson_control_example}
\begin{aligned}
    &\text{Minimize} && J(u, f) = \frac{1}{2} \int_0^1 \int_0^1 \left[ u(x, y) - u_d(x, y) \right]^2 \, dx \, dy + \frac{ \rho}{2} \int_0^1 \int_0^1 f(x, y)^2 \, dx \, dy, \\
    &\text{subject to} && -\Delta  u(x,y) = f(x, y), \quad (x, y) \in [0,1] \times [0,1],
\end{aligned}
\end{equation}
with $u(x,y)=0$ on the boundary and the desired state 
$$u_d(x, y) = A\sin(\pi x) \sin(\pi y),$$ with parameter $A = 10$. The control for this problem is the distributed control $f(x, y)$ over the domain. This problem also has an analytical solution, given by:
\begin{equation}\label{eq:u_analytical_2D}
u_{\text{analytical}}(x, y) = \frac{A}{1 + 4 \rho \pi^4} \sin(\pi x) \sin(\pi y),
\end{equation}
and the exact control is:
\begin{equation}\label{eq:f_analytical_2D}
f_{\text{analytical}}(x, y) = \frac{2 \pi^2 A}{1 + 4 \rho \pi^4} \sin(\pi x) \sin(\pi y).
\end{equation}
In this example, we take $\rho = 0.01$. We present the desired state $u_d$ and the analytical solution $u_{\text{analytical}}$ in Figure \ref{fig:2D_Poisson_desired_analytical}. Under this choice of $\rho$, there is an apparent difference between the desired state and the analytical solution to the control problem. We deliberately choose such an example instead of a problem with a much smaller regularization parameter (for example, $\rho = 0.0001$) to ensure that the neural network genuinely finds the solution to the optimal control problem rather than merely solving an approximation problem that approximates the desired state as if learning a given function.

\begin{figure}
    \centering
    
    \begin{subfigure}[t]{0.42\textwidth}
        \centering
        \includegraphics[width=\textwidth]{exact_u_depth_4_width_60_w_e_2000_reg_0.01_2D_Poisson.png}
        \subcaption{Exact $u$}
    \end{subfigure}
    \begin{subfigure}[t]{0.42\textwidth}
        \centering
        \includegraphics[width=\textwidth]{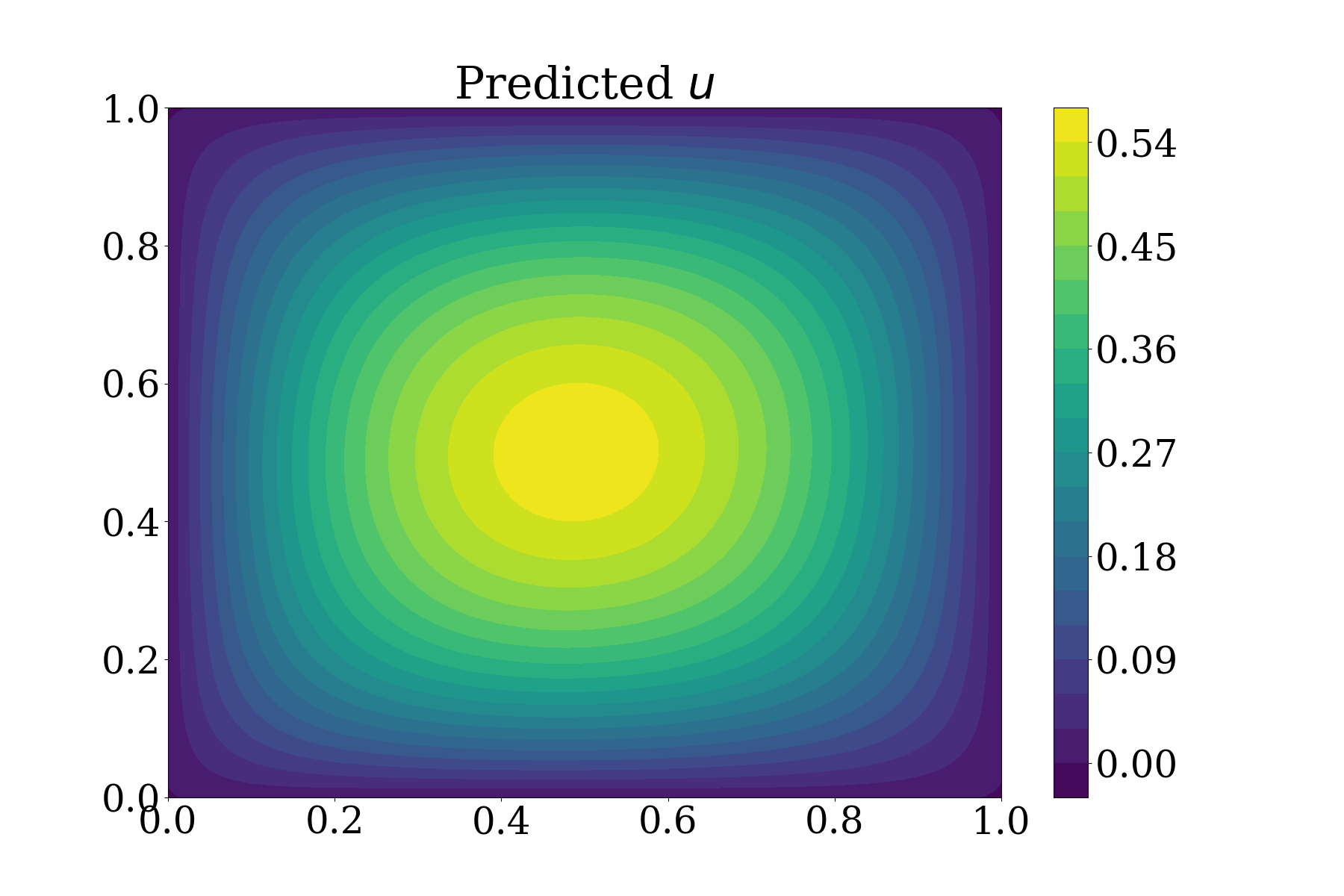}
        \subcaption{Numerical $u$}
    \end{subfigure}

    \begin{subfigure}[t]{0.5\textwidth}
        \centering
        \includegraphics[width=\textwidth]{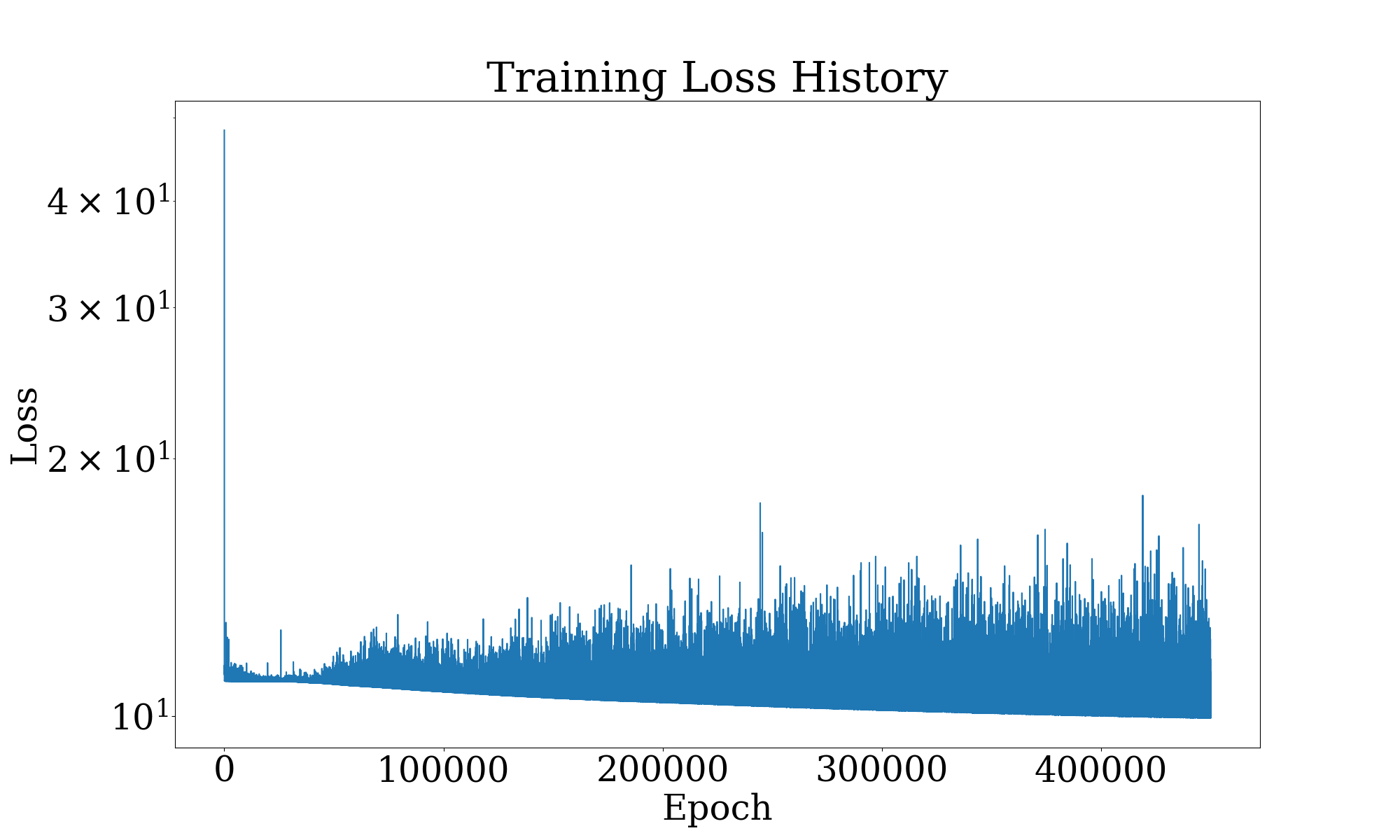}
        \subcaption{Loss over epochs}
    \end{subfigure}
    \caption{Results for minimizing the loss function with $\lambda_p = 2000$, $\lambda_b = 2000$, and $\rho = 0.01$. Subfigures (a) and (b) show the exact and numerical solutions for $u$, respectively. 
    Subfigure (c) presents the loss over epochs.}
    \label{fig:2D_Poisson_model}
\end{figure}

Firstly, as a naive approach, we aim to solve this problem using a penalty formulation. We consider using a neural network with a depth of $4$ layers and width of $60$ neurons to minimize the following loss function:
\begin{equation}\label{eq:2D_Poisson_penalty_loss_nn}
\begin{aligned}
   L[{x, y}; \theta] = 
   \frac{1}{2N^2} \sum_{m=1}^{N} \sum_{n=1}^N&\left[ \unn(x_m, y_n; \theta) - u_d(x_m, y_n) \right]^2 + \frac{\rho}{2N^2} \sum_{m=1}^{N} \sum_{n=1}^N \fnn (x_m, y_n; \theta)^2  \\
   & + \frac{\lambda_p}{N^2} \sum_{m=1}^{N} \sum_{n=1}^N\left[ \frac{\partial^2 \unn}{\partial x^2}(x_m, y_n; \theta) + \frac{\partial^2 \unn}{\partial y^2}(x_m, y_n; \theta) + \fnn(x_m, y_n; \theta) \right]^2\\
   &\quad\quad\quad+ \frac{\lambda_b}{(N_b)} \sum_{k=1}^{N_b}\left[ \unn(x_k^b, y_k^b; \theta) \right]^2,
\end{aligned}
\end{equation}
The weights for the equation and boundary loss terms are set to $\lambda_e = \lambda_b = 2000$. The training data used to compute the objective function and equation loss are selected as $\{x_m\}_{m=1}^N, \{y_n\}_{n=1}^N$ with $N = 16$, resulting in a total of 256 interior points. The training data used to compute the boundary loss are selected as $\{(x_k^b,y_k^b)\}_{k=1}^{N_b}$ with $N_b = 32$. Both the interior points and boundary data are uniformly sampled from their corresponding domains. For the training details, the initial learning rate is set to $0.001$, which can be reduced to a minimum learning rate of $0.0001$ based on the learning rate scheduling strategy, with patience set to 3000 epochs. The training is conducted over a maximum of 450000 epochs. The numerical findings are presented in Figure \ref{fig:2D_Poisson_model}.

The results indicate that the numerical solution does not approximate the analytical solution accurately. The exact value for $u$ ranges from 0 to approximately $2.04$, while the numerical value for $u$ remains below $0.6$, indicating a significant discrepancy. Although the training loss continues to decrease, it is evident that the model does not converge to the correct solution even after 450000 epochs. While we cannot exclude the possibility that the network might eventually converge to the actual solution with more training epochs, this convergence has not been observed at this point after such a large number of epochs. This phenomenon showcases the limitations of a direct approach using penalty formulation with large penalty parameters to implement the PDE constraints.

In contrast, let us apply the penalty adversarial network to solve the same problem. The loss function that the discriminator network aims to minimize is given by:
\begin{equation}\label{eq:2D_Poisson_penalty_loss_nn_discriminator}
\begin{aligned}
   L^d[{x, y}; \theta] = 
   \frac{1}{2N^2} &\sum_{m=1}^{N} \sum_{n=1}^N\left[ \udnn(x_m, y_n; \theta) - u_d(x_m, y_n) \right]^2 + \frac{\rho}{2N^2} \sum_{m=1}^{N} \sum_{n=1}^N \fdnn(x_m, y_n; \theta)^2  \\
   & + \frac{\lambda_p^d}{N^2} \sum_{m=1}^{N} \sum_{n=1}^N\left[ \frac{\partial^2 \udnn}{\partial x^2}(x_m, y_n; \theta) + \frac{\partial^2 \udnn}{\partial y^2}(x_m, y_n; \theta) + \fdnn(x_m, y_n; \theta) \right]^2\\
   &\quad\quad\quad+ \frac{\lambda^d_b}{(N_b)} \sum_{k=1}^{N_b}\left[ \udnn(x_k^b, y_k^b; \theta) \right]^2.
\end{aligned}
\end{equation}
On the other hand, the solver network aims to minimize:
\begin{equation}\label{eq:2D_Poisson_penalty_loss_nn_solver}
\begin{aligned}
   L^s[{x, y}; \theta] &= 
   \frac{1}{2N^2} \sum_{m=1}^{N} \sum_{n=1}^N\left[ \usnn(x_m, y_n; \theta) - u_d(x_m, y_n) \right]^2 + \frac{\rho}{2N^2} \sum_{m=1}^{N} \sum_{n=1}^N \fsnn(x_m, y_n; \theta)^2  \\
   &\quad\quad + \frac{\lambda_p^s}{N^2} \sum_{m=1}^{N} \sum_{n=1}^N\left[ \frac{\partial^2 \usnn}{\partial x^2}(x_m, y_n; \theta) + \frac{\partial^2 \usnn}{\partial y^2}(x_m, y_n; \theta) + \fsnn(x_m, y_n; \theta) \right]^2\\
&\quad\quad+ \frac{\lambda^s_b}{(N_b)} \sum_{k=1}^{N_b}\left[ \usnn(x_k^b, y_k^b; \theta) \right]^2\\
    &\quad\quad+\omega\, \bigg{\{}\frac{\rho}{2N} \sum_{m=1}^{N} \sum_{n=1}^N \fsnn(x_m, y_n; \theta)^2-\frac{\rho}{2N} \sum_{m=1}^{N} \sum_{n=1}^N \fdnn(x_m, y_n; \theta)^2\\
    &\,\quad\quad\quad\,\quad\,\quad+  \frac{1}{2N^2} \sum_{m=1}^{N} \sum_{n=1}^N\left[ \usnn(x_m, y_n; \theta) - u_d(x_m, y_n) \right]^2  \\
  &\quad\,\quad\,\quad\quad \quad\,\quad\,\quad\quad -   \frac{1}{2N^2} \sum_{m=1}^{N} \sum_{n=1}^N\left[ \udnn(x_m, y_n; \theta) - u_d(x_m, y_n) \right]^2 \bigg{\}}^2.
\end{aligned}
\end{equation}
The parameters are chosen as $\lambda_p^s=\lambda_b^s = 2000$, $\lambda_p^d=\lambda_b^d = 10$ and $\omega = 100$. All other parameters related to the structure and the training data remain the same. 

\begin{figure}
    \centering
    \begin{subfigure}[t]{0.32\textwidth}
        \centering
        \includegraphics[width=\textwidth]{exact_u_depth_4_width_60_w_e_2000_reg_0.01_2D_Poisson.png}
        \subcaption{Exact $u$}
    \end{subfigure}
    \begin{subfigure}[t]{0.32\textwidth}
        \centering
        \includegraphics[width=\textwidth]{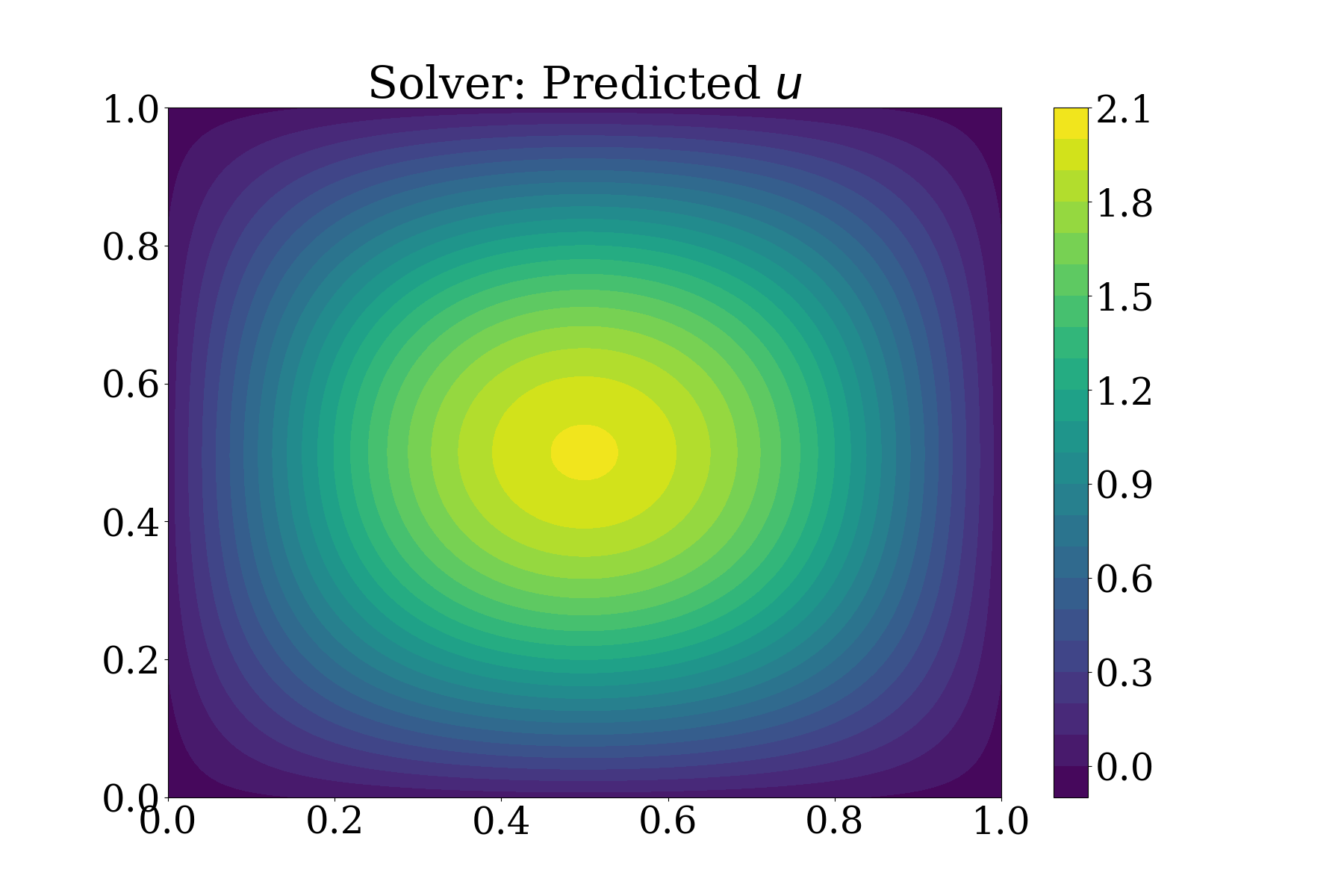}
        \subcaption{Solver: Predicted $u$}
    \end{subfigure}
    \begin{subfigure}[t]{0.32\textwidth}
        \centering
        \includegraphics[width=\textwidth]{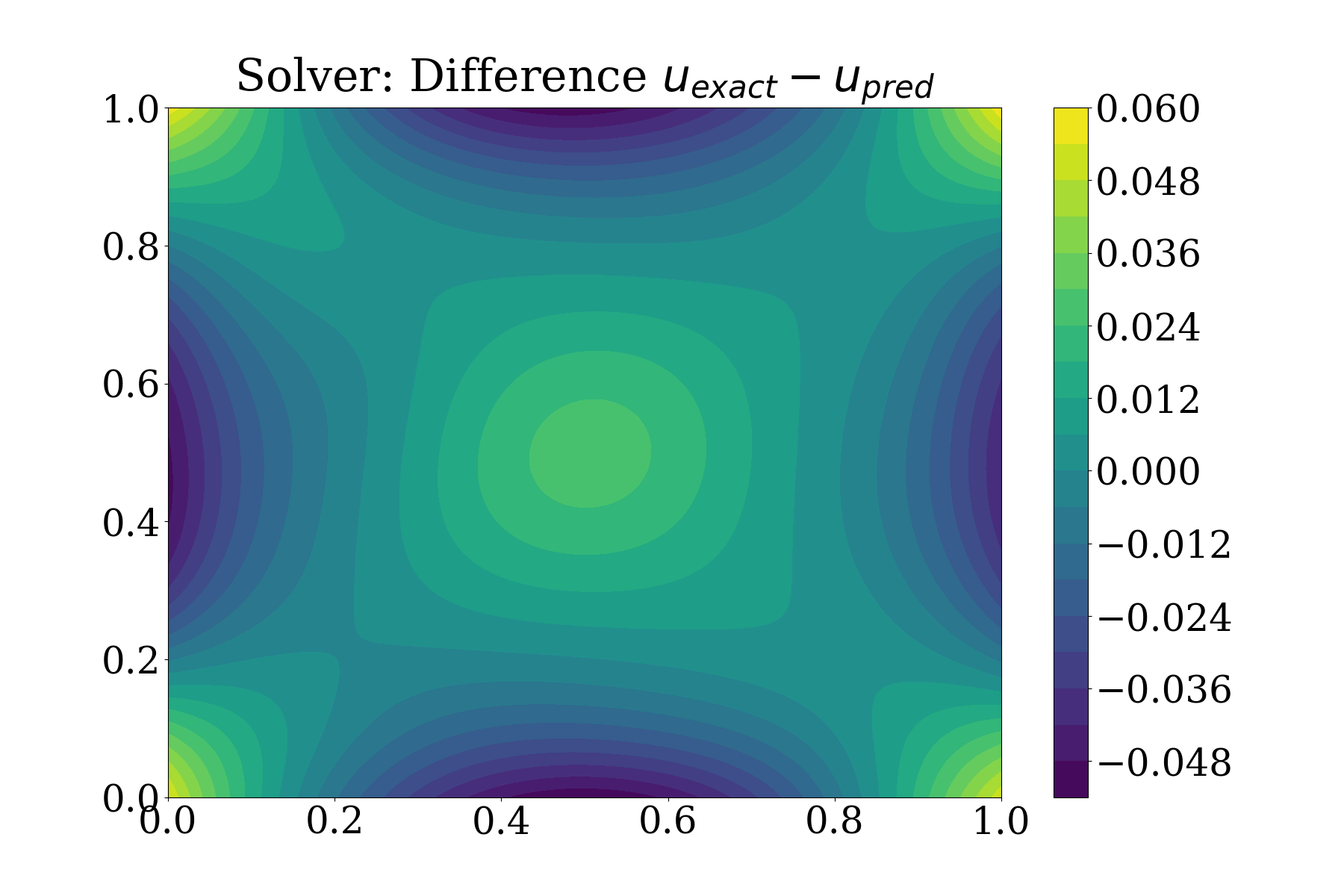}
        \subcaption{Solver: Difference $u$}
    \end{subfigure}
    
    \begin{subfigure}[t]{0.32\textwidth}
        \centering
        \includegraphics[width=\textwidth]{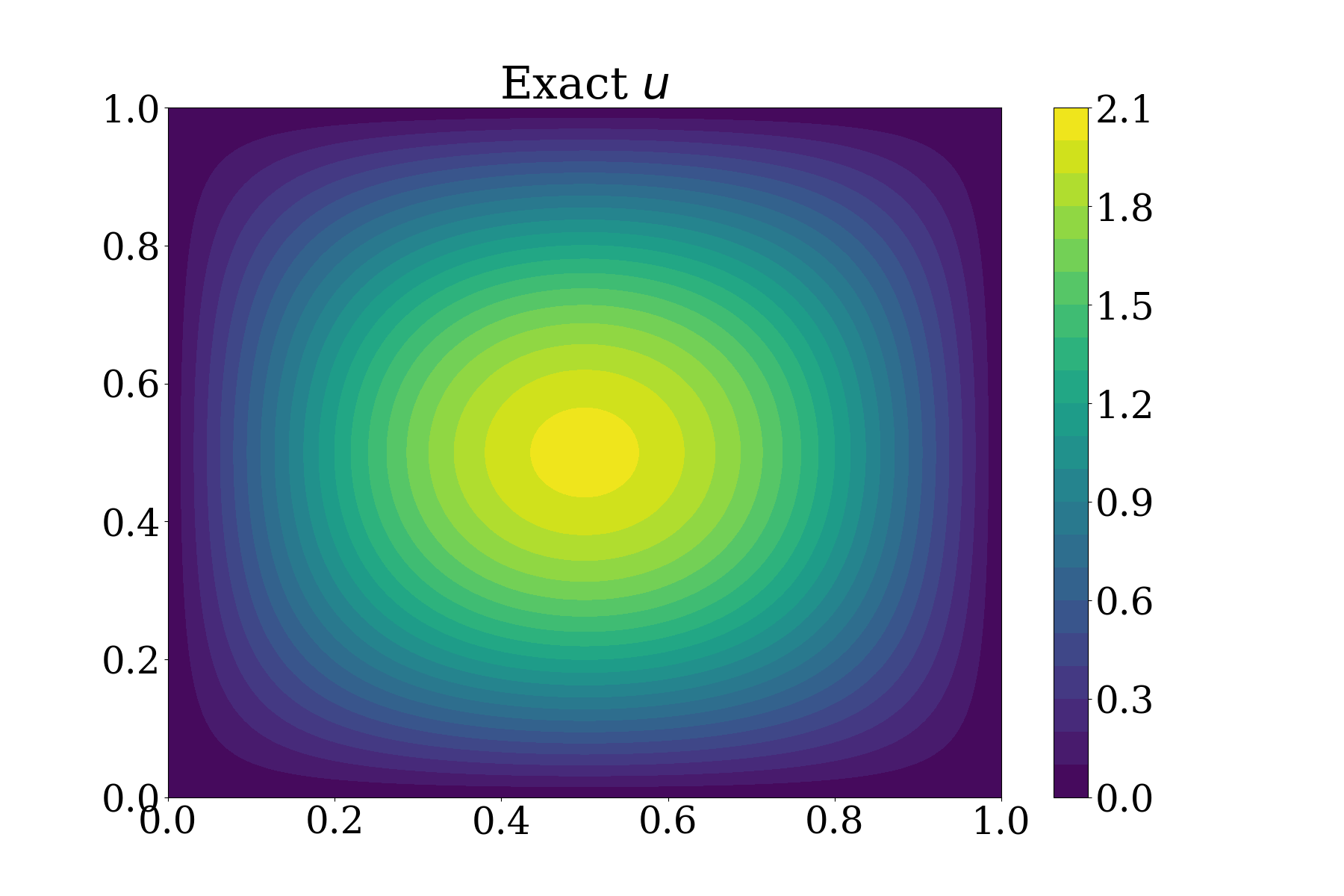}
        \subcaption{Exact $u$}
    \end{subfigure}
    \begin{subfigure}[t]{0.32\textwidth}
        \centering
        \includegraphics[width=\textwidth]{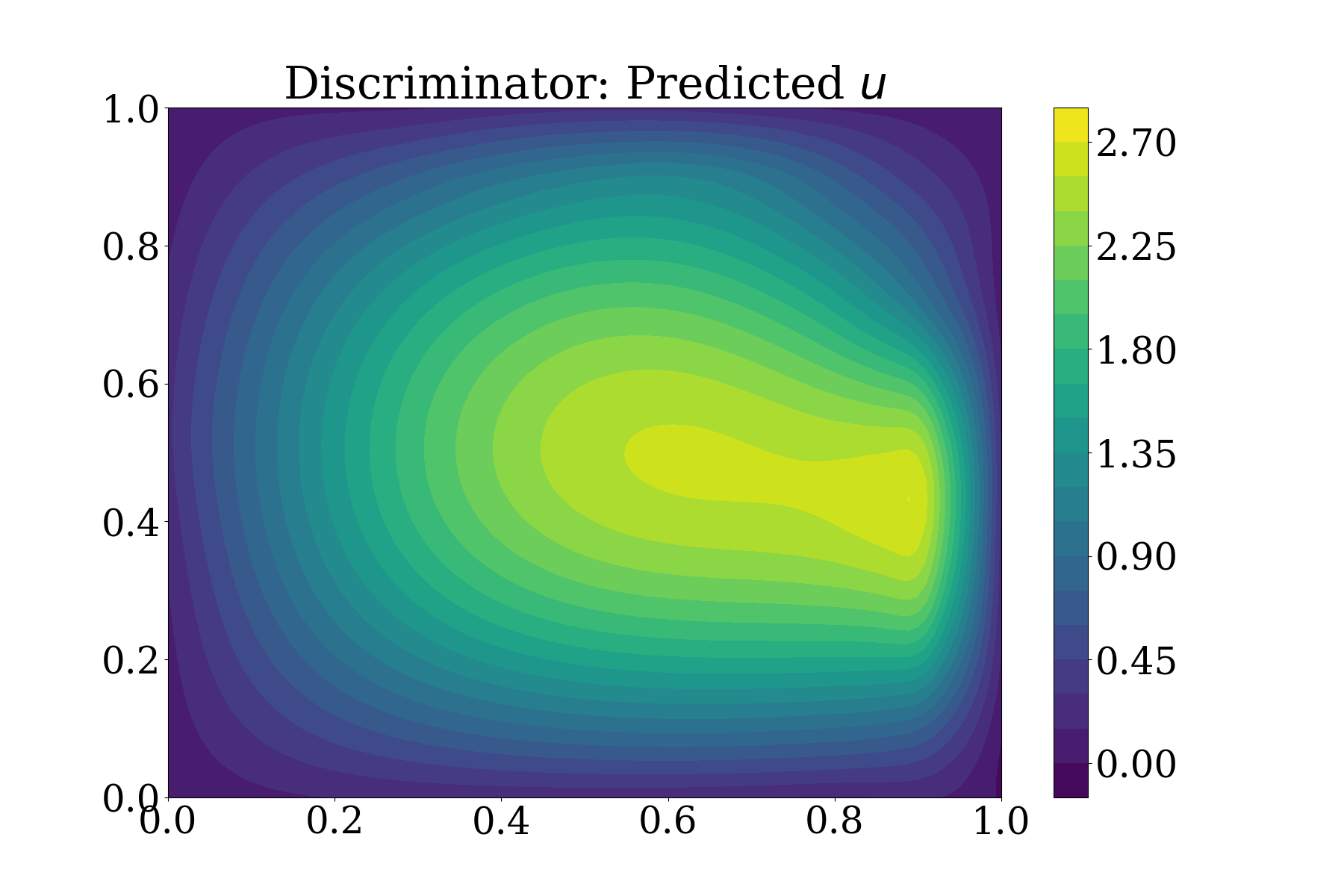}
        \subcaption{Discriminator: Predicted $u$}
    \end{subfigure}
    \begin{subfigure}[t]{0.32\textwidth}
        \centering
        \includegraphics[width=\textwidth]{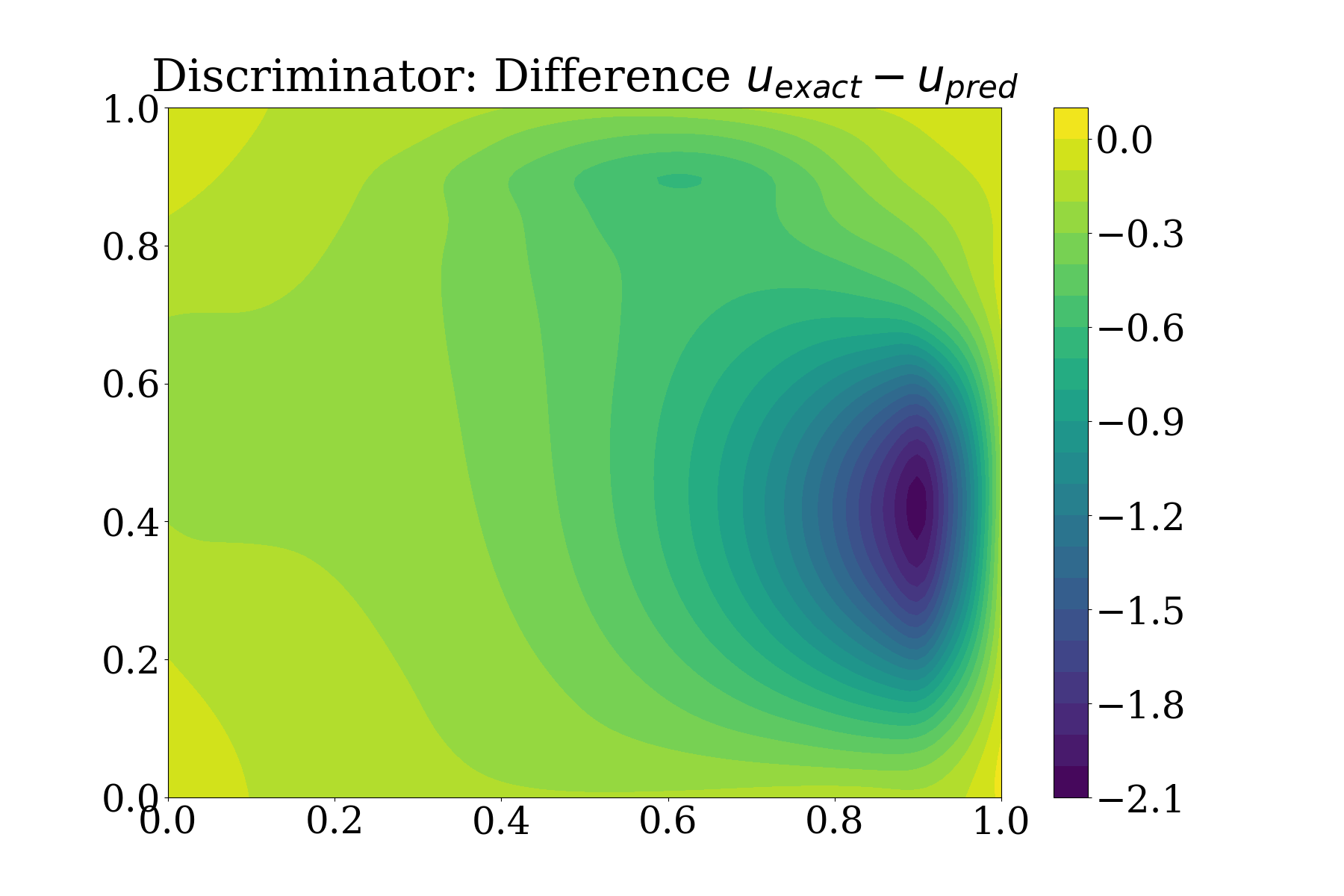}
        \subcaption{Discriminator: Difference $u$}
    \end{subfigure}

    \begin{subfigure}[t]{0.4\textwidth}
        \centering
        \includegraphics[width=\textwidth]{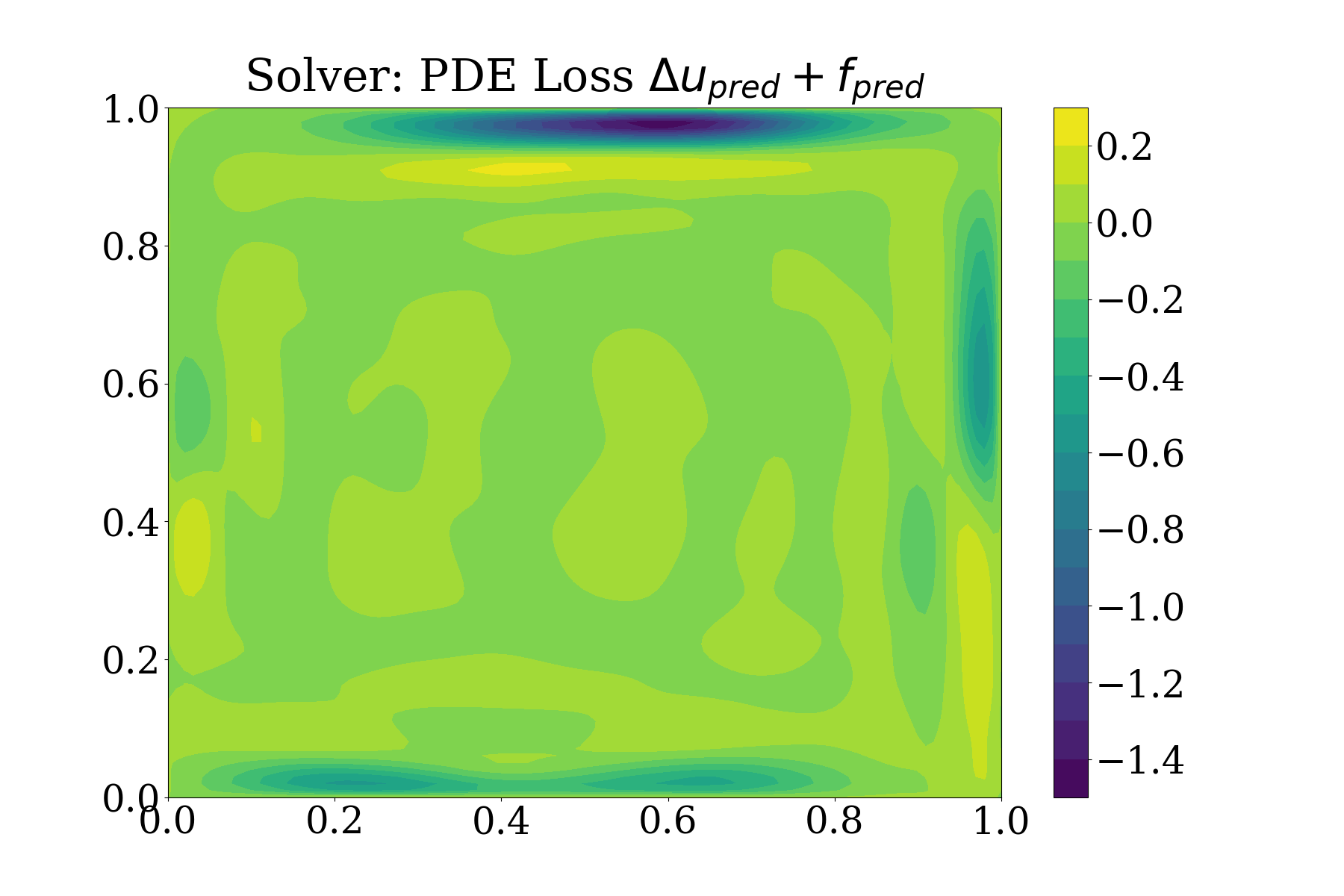}
        \subcaption{Solver: PDE Loss $\Delta u_{pred} + f_{pred}$}
    \end{subfigure}
    \begin{subfigure}[t]{0.4\textwidth}
        \centering
        \includegraphics[width=\textwidth]{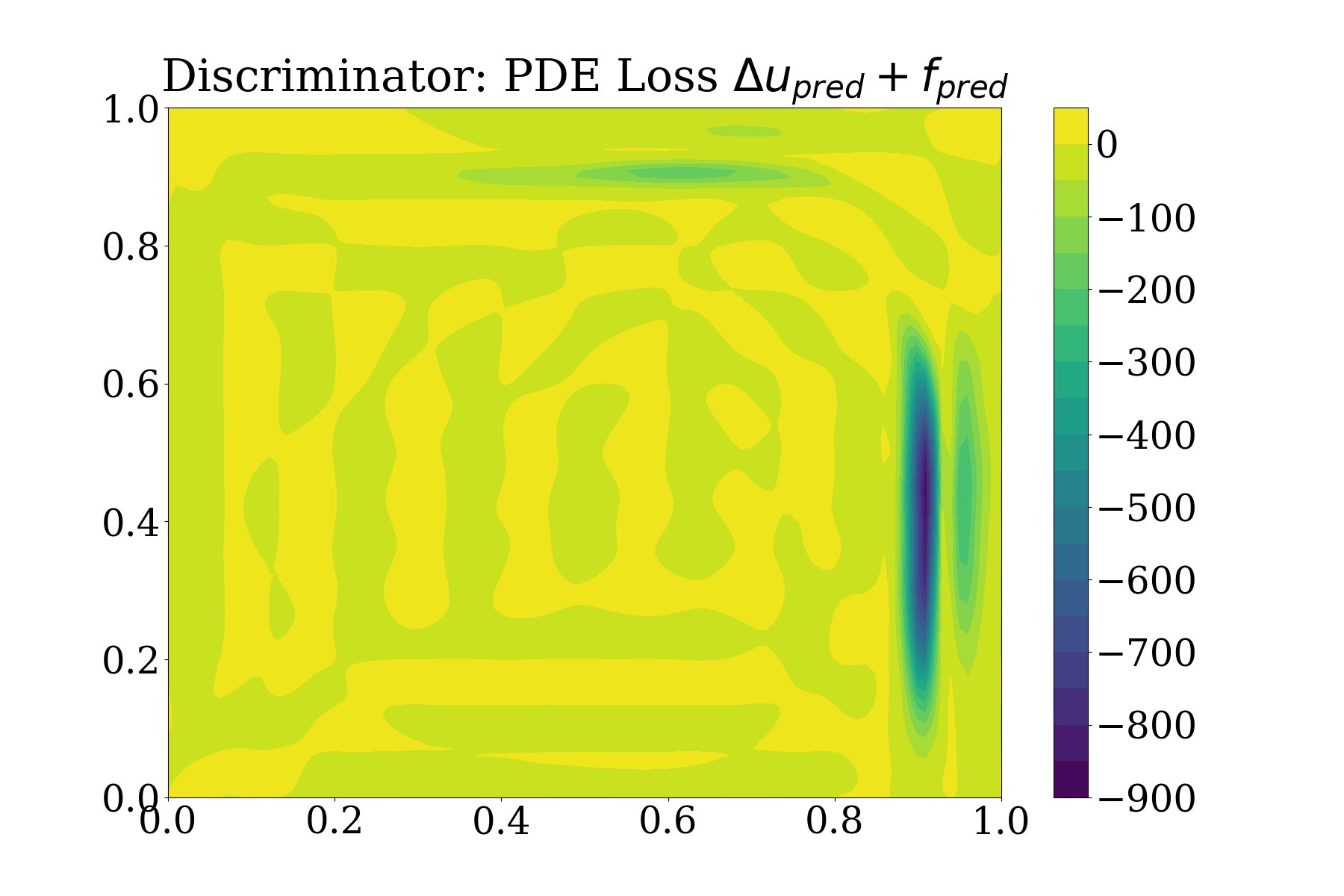}
        \subcaption{Discriminator: PDE Loss $\Delta u_{pred} + f_{pred}$}
    \end{subfigure}

    \begin{subfigure}[t]{0.5\textwidth}
        \centering
        \includegraphics[width=\textwidth]{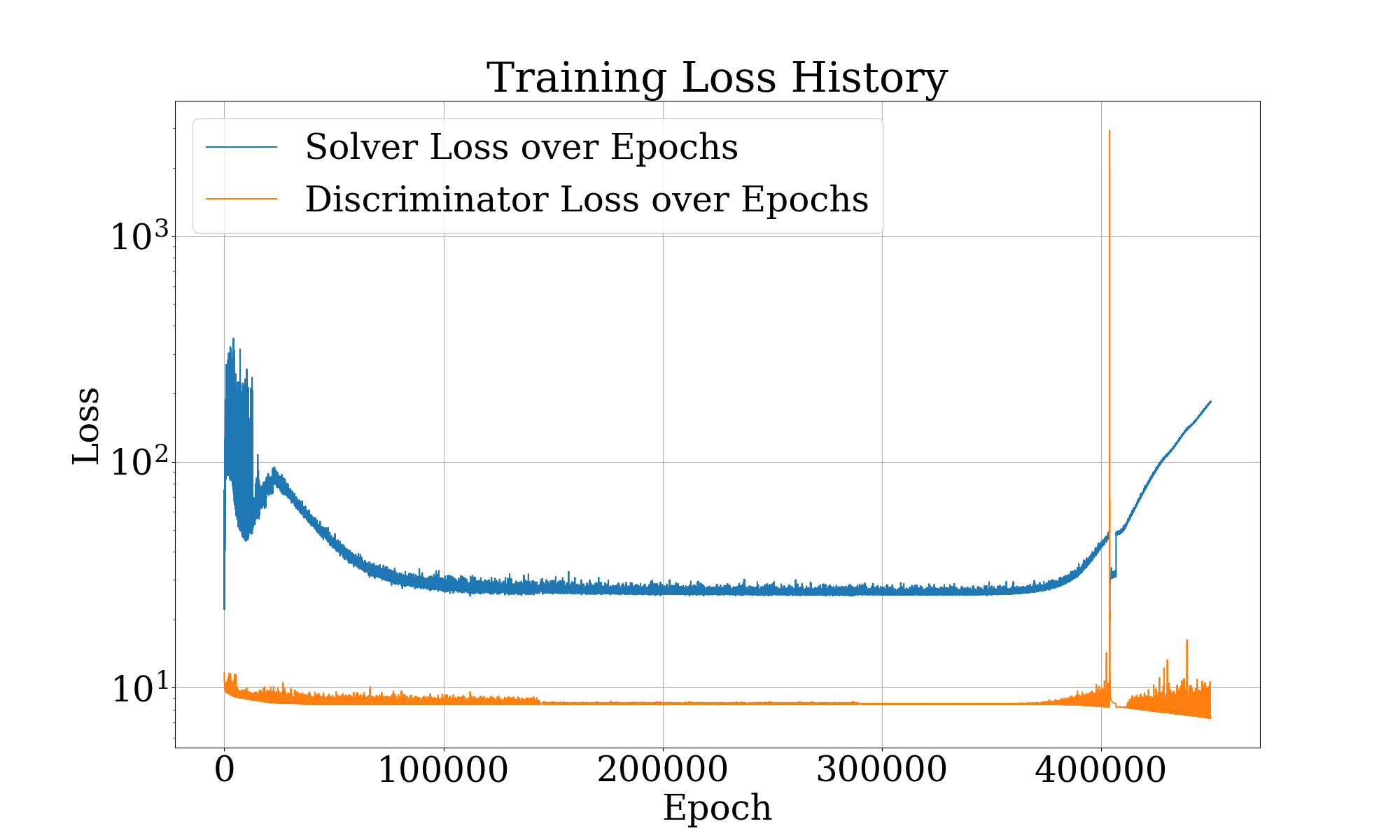}
        \subcaption{Loss over epochs}
    \end{subfigure}
    \caption{Results for minimizing the loss function \eqref{eq:2D_Poisson_penalty_loss_nn_discriminator} and \eqref{eq:2D_Poisson_penalty_loss_nn_solver} with $\lambda_p^d=\lambda_b^d = 10$, $\lambda_p^s=\lambda_b^s = 2000$,  $\omega = 100$ and $\rho = 0.01$. Subfigures (a) and (b) show the exact and numerical solutions for $u$ for the solver, respectively. Subfigure (c) shows the difference between exact and predicted $u$ for the solver. Subfigures (d) and (e) show the exact and numerical solutions for $u$ for the discriminator, respectively. Subfigure (f) shows the difference between exact and predicted $u$ for the discriminator. Subfigure (g) presents the PDE loss for the solver, and subfigure (h) presents the PDE loss for the discriminator. Subfigure (i) presents the loss over epochs.}
    \label{fig:2D_Poisson_adversarial}
\end{figure}

The numerical results of applying the penalty adversarial network to solve problem \eqref{eq:2D_Poisson_control_example} are presented in Figure \ref{fig:2D_Poisson_adversarial}. We observe that the solver successfully finds the real solution. Subplot (c) in Figure \ref{fig:2D_Poisson_adversarial} shows that the difference between the exact solution, shown in Subplot (a), and the predicted result from the solver, shown in Subplot (b), is minimal. However, the discriminator fails to reach the exact solution, exhibiting a significant error in the bottom-right area compared to other regions, as shown in Subplot (f).

Focusing on the ability of the networks to adhere to the constraints, as reflected by the PDE losses, we see that the solver's prediction for $u$ and $f$ generally follows the Poisson equation quite well in most areas, except for some deviation at the top boundary. In comparison, while the discriminator's prediction also generally follows the Poisson equation, it deviates strongly in the bottom-right area, corresponding to the same region where the discriminator makes a significant error in predicting the values of $u$.

Examining the evolution of the loss function over epochs, we observe an interesting phenomenon. As seen in Subplot (i) of Figure \ref{fig:2D_Poisson_adversarial}, the loss for the discriminator decreases slowly and remains almost flat during the first 400000 epochs. However, after 400000 epochs, the loss decreases much more rapidly. At the same time, the loss for the solver increases quickly due to the difference in the objective value between the solver and the discriminator, which causes the additional penalty term introduced by the penalty adversarial approach to increase rapidly.

Though we do not fully understand this behavior at present, we propose one potential explanation: both the solver and the discriminator gradually converge toward the analytical solution at first. However, since the discriminator uses a small penalty parameter, theoretically, the minimum that the discriminator network can converge to deviates a certain distance from the analytical solution $u_{\text{analytical}}$. This deviation arises from a tradeoff between disobeying the constraints and decreasing the objective function, leading to a situation where the corresponding numerical solution might attempt to decrease the objective function by strongly disobeying the PDE constraints in certain areas. This explanation aligns with the observed phenomenon in the discriminator's result, which significantly deviates from satisfying the Poisson equation in a specific region. Comparatively, although the loss function for the solver increases after 400000 epochs, the minimum that the solver network reaches does not change significantly and remains close to the real analytical solution.

\subsubsection{Example 3: Distributed Control Problem Constrained by 2D Allen-Cahn Equation}

As the last example presented here, we will showcase that the penalty adversarial network can also be applied to solve nonlinear problems, albeit more complicated. We consider a distributed control problem constrained by a 2D Allen-Cahn equation, which has received considerable research interest from the computational community \cite{guillen2014second,lwight2021solving,yang2017linearly}.

\begin{figure}
    \centering
    
    \begin{subfigure}[t]{0.45\textwidth}
        \centering
        \includegraphics[width=\textwidth]{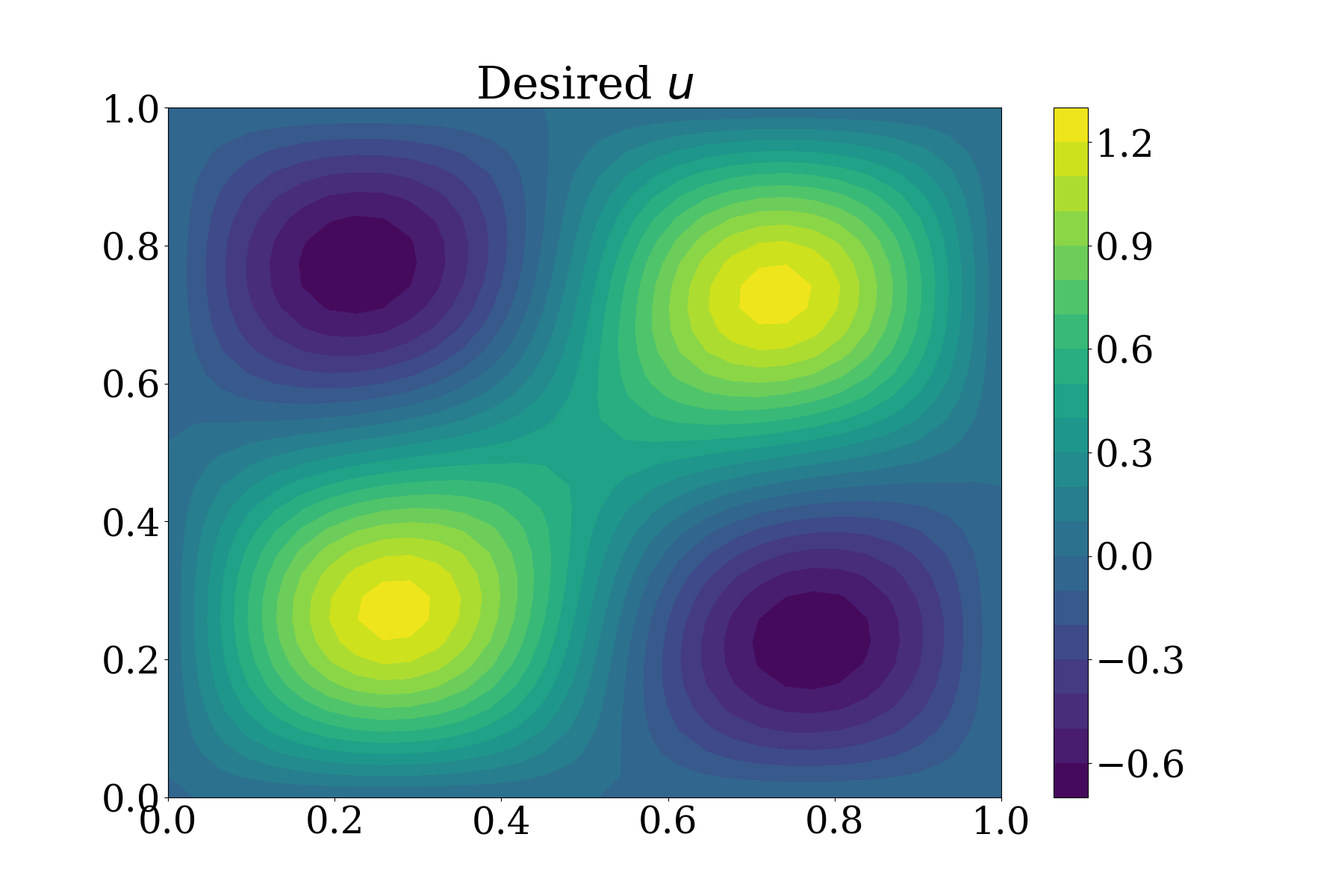}
        \subcaption{Desired state $u_d$}
    \end{subfigure}
    \begin{subfigure}[t]{0.45\textwidth}
        \centering
\includegraphics[width=\textwidth]{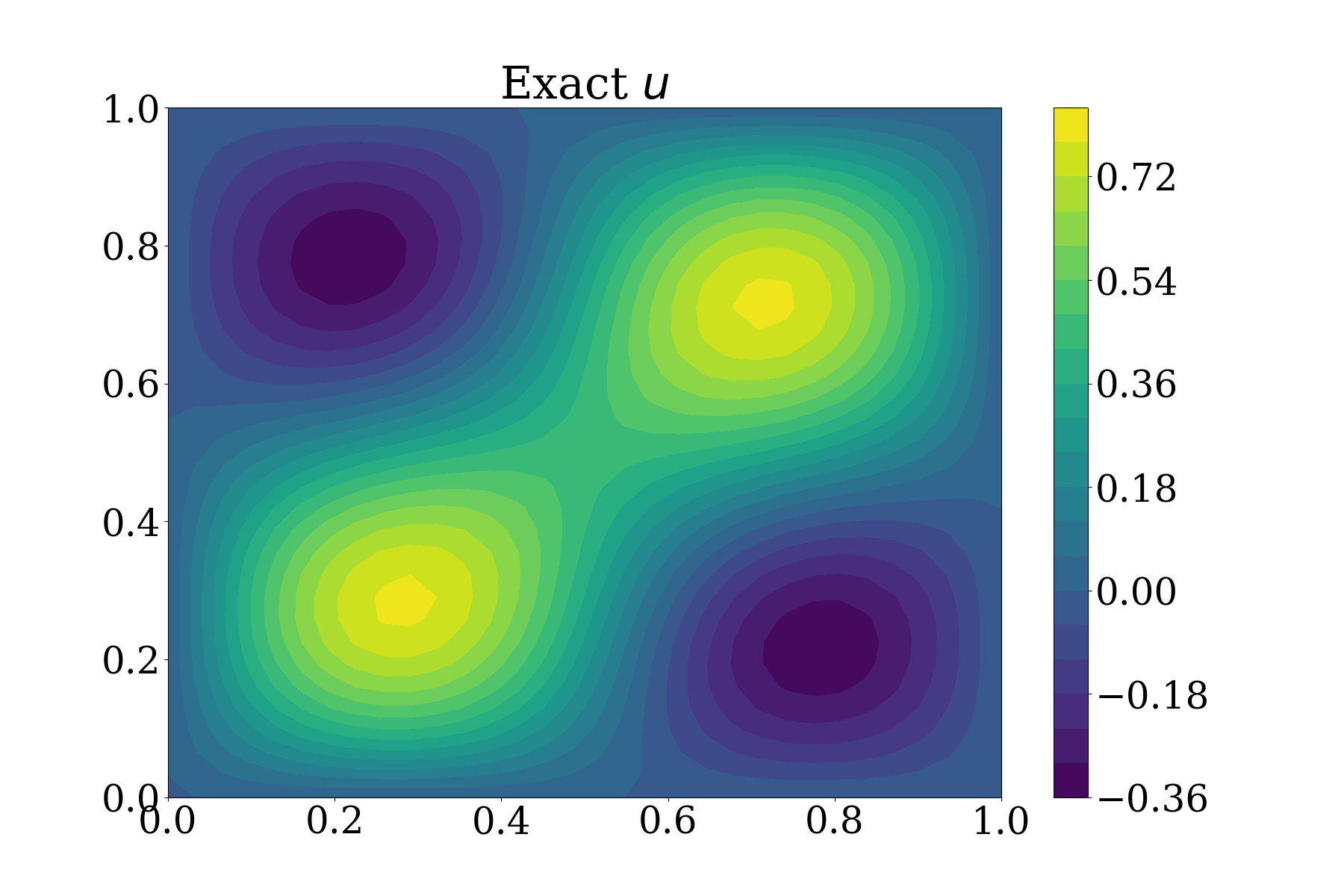}
         \subcaption{Analytical solution $u_{\text{analytical}}$}
    \end{subfigure}
    
    \caption{The desired state and analytical solution for problem \eqref{eq:2D_Allen_Cahn_control_example}}.
    \label{fig:2D_Allen_Cahn_desired_analytical}
\end{figure}

The specific optimal control problem has the same formulation for the objective functional as the previous example \eqref{eq:2D_Poisson_control_example} but is constrained by a different equation, given by:
\begin{equation}\label{eq:2D_Allen_Cahn_control_example}
\begin{aligned}
    &\text{Minimize} && J(u, f) = \frac{1}{2} \int_0^1 \int_0^1 \left[ u(x, y) - u_d(x, y) \right]^2 \, dx \, dy + \frac{ \rho}{2} \int_0^1 \int_0^1 f(x, y)^2 \, dx \, dy, \\
    &\text{subject to} && -\Delta  u(x,y) + \frac{1}{\epsilon^2}\left[u(x,y)^3-u(x,y)\right] = f(x, y), \quad (x, y) \in [0,1] \times [0,1],
\end{aligned}
\end{equation}
with $u(x,y) = 0$ on the boundary and a chosen desired state. The parameter $\epsilon > 0$ is related to the PDE constraint itself. Using an adjoint formulation \cite{bertsekas2014constrained}, we can find the following pair of functions that serve as a solution to problem \eqref{eq:2D_Allen_Cahn_control_example}: We set

\begin{equation}
    u_{\text{analytical}}(x, y) = \alpha \sin(\pi x) \sin(\pi y) + \beta \sin(2\pi x) \sin(2\pi y),
\end{equation}
and
\begin{equation}
\begin{aligned}
   & \quad\,\, f_{\text{analytical}}(x, y) \\&=  2\pi^2 \left[\alpha \sin(\pi x) \sin(\pi y) + 4\beta \sin(2\pi x) \sin(2\pi y)\right]\\
     &\quad\quad+ \frac{1}{\epsilon^2}\left[\left(\alpha \sin(\pi x) \sin(\pi y) + \beta \sin(2\pi x) \sin(2\pi y)\right)^3 - \alpha \sin(\pi x) \sin(\pi y) - \beta \sin(2\pi x) \sin(2\pi y)\right],
\end{aligned}
\end{equation}
while the desired state is
\begin{equation}
    \begin{aligned}
u_{d}(x, y) &= u_{\text{analytical}}(x, y) + \rho\, \Delta^2 u_{\text{analytical}}(x, y) -  \frac{3\rho}{\epsilon^2}u_{\text{analytical}}^2(x, y) \Delta u_{\text{analytical}}(x, y) \\
&\quad - \frac{\rho}{\epsilon^2} \left[  6\,  u_{\text{analytical}}(x, y) \left\lvert \nabla u_{\text{analytical}}(x, y)\right\rvert^2 + \Delta u_{\text{analytical}}(x, y) \right] \\
&\quad - \frac{\rho}{\epsilon^2} \left[  \Delta u_{\text{analytical}}(x, y) - \frac{1}{\epsilon^2} \left( u_{\text{analytical}}^3(x, y) - u_{\text{analytical}}(x, y) \right) \right] \left(3 u_{\text{analytical}}^2(x, y) - 1\right).
\end{aligned}
\end{equation}
Here, $\alpha, \beta\in\mathbb{R}$ are parameters. Note that if $\rho = 0$, then $u_d = u_{\text{analytical}}$, which coincides with the expectation that, in this case, the minimization problem \eqref{eq:2D_Allen_Cahn_control_example} is simply equivalent to learning a known function. However, when $\rho > 0$, $u_d$ can have obvious differences from $u_{\text{analytical}}$. In the example presented here, we choose $\epsilon = 0.4$, $\alpha = 0.45$, $\beta = 0.55$, and $\rho = 0.0001$. Under these choices, the corresponding desired state $u_d$ and $u_{\text{analytical}}$ are plotted in Figure \ref{fig:2D_Allen_Cahn_desired_analytical}. Obvious differences in values at the same points can be observed from this comparison, preventing the network from simply learning a function.

For the sake of brevity of our presentation, we will skip the experiment of using a naive approach as a penalty formulation with a large penalty parameter. As expected, it will fail to work. We will only present the results of applying a penalty adversarial network. For this approach, we define the corresponding loss function for the discriminator and solver network, respectively, as
\begin{equation}\label{eq:2D_Allen_Cahn_penalty_loss_nn_discriminator}
\begin{aligned}
 L^d[{x, y}; \theta] = 
   \frac{1}{2N^2} &\sum_{m=1}^{N} \sum_{n=1}^N\left[ \udnn(x_m, y_n; \theta) - u_d(x_m, y_n) \right]^2 + \frac{\rho}{2N^2} \sum_{m=1}^{N} \sum_{n=1}^N \fdnn(x_m, y_n; \theta)^2  \\
   & + \frac{\lambda_p^d}{N^2} \sum_{m=1}^{N} \sum_{n=1}^N\Bigg[ \frac{\partial^2 \udnn}{\partial x^2}(x_m, y_n; \theta) + \frac{\partial^2 \udnn}{\partial y^2}(x_m, y_n; \theta) \\
   &\quad\quad\quad\quad\quad\quad\quad\quad -\frac{1}{\epsilon^2}\left(\udnn(x_m, y_n; \theta)^3-\udnn(x_m, y_n; \theta)\right) +\fdnn(x_m, y_n; \theta) \Bigg]^2\\
   &+ \frac{\lambda^d_b}{(N_b)} \sum_{k=1}^{N_b}\left[ \udnn(x_k^b, y_k^b; \theta) \right]^2,
\end{aligned}
\end{equation}
and
\begin{equation}\label{eq:2D_Allen_Cahn_penalty_loss_nn_solver}
\begin{aligned}
   L^s[{x, y}; \theta] &= 
   \frac{1}{2N^2} \sum_{m=1}^{N} \sum_{n=1}^N\left[ \usnn(x_m, y_n; \theta) - u_d(x_m, y_n) \right]^2 + \frac{\rho}{2N^2} \sum_{m=1}^{N} \sum_{n=1}^N \fsnn(x_m, y_n; \theta)^2  \\
    &\quad\quad + \frac{\lambda_p^d}{N^2} \sum_{m=1}^{N} \sum_{n=1}^N\Bigg[ \frac{\partial^2 \udnn}{\partial x^2}(x_m, y_n; \theta) + \frac{\partial^2 \udnn}{\partial y^2}(x_m, y_n; \theta) \\
   &\quad\quad\quad\quad\quad\quad\quad\quad -\frac{1}{\epsilon^2}\left(\udnn(x_m, y_n; \theta)^3-\udnn(x_m, y_n; \theta)\right) +\fdnn(x_m, y_n; \theta) \Bigg]^2\\
&\quad\quad+ \frac{\lambda^s_b}{(N_b)} \sum_{k=1}^{N_b}\left[ \usnn(x_k^b, y_k^b; \theta) \right]^2\\
    &\quad\quad+\omega\, \bigg{\{}\frac{\rho}{2N^2} \sum_{m=1}^{N} \sum_{n=1}^N \fsnn(x_m, y_n; \theta)^2-\frac{\rho}{2N^2} \sum_{m=1}^{N} \sum_{n=1}^N \fdnn(x_m, y_n; \theta)^2\\
    &\,\quad\quad\quad\,\quad\,\quad+  \frac{1}{2N^2} \sum_{m=1}^{N} \sum_{n=1}^N\left[ \usnn(x_m, y_n; \theta) - u_d(x_m, y_n) \right]^2  \\
  &\quad\,\quad\,\quad\quad \quad\,\quad\,\quad\quad -   \frac{1}{2N^2} \sum_{m=1}^{N} \sum_{n=1}^N\left[ \udnn(x_m, y_n; \theta) - u_d(x_m, y_n) \right]^2 \bigg{\}}^2.
\end{aligned}
\end{equation}

The hyperparameters for this experiment are configured as follows: The training data for both the objective function and equation loss are sampled as $\{x_m\}_{m=1}^N$ and $\{y_n\}_{n=1}^N$, with $N = 32$. The boundary training data consist of $\{(x_k^b, y_k^b)\}_{b=1}^{N_b}$, with $N_b = 32$. Both the interior and boundary points are uniformly sampled from their respective domains. For the training specifics, the initial learning rate is set at $0.001$, potentially decreasing to a minimum of $0.0001$ using a learning rate scheduling strategy with a patience parameter of 10000 epochs. The learning rate adjustment does not begin until after 300000 epochs to avoid early-stage inaccuracies. The total number of epochs is extended to 1.5 million to capture the full training dynamics, as the solver network, as seen in the experiment results, continues to improve throughout. This extended training period also highlights the common challenges of addressing nonlinear equations in PINN applications \cite{mowlavi2023optimal}.

\begin{figure}
    \centering
    \begin{subfigure}[t]{0.32\textwidth}
        \centering
        \includegraphics[width=\textwidth]{exact_u_depth_4_width_60_w_e_s_1000_w_e_d_0.2_w_d_20000_reg_0.0001_2D_Allen_Cahn.png}
        \subcaption{Exact $u$}
    \end{subfigure}
    \begin{subfigure}[t]{0.32\textwidth}
        \centering
        \includegraphics[width=\textwidth]{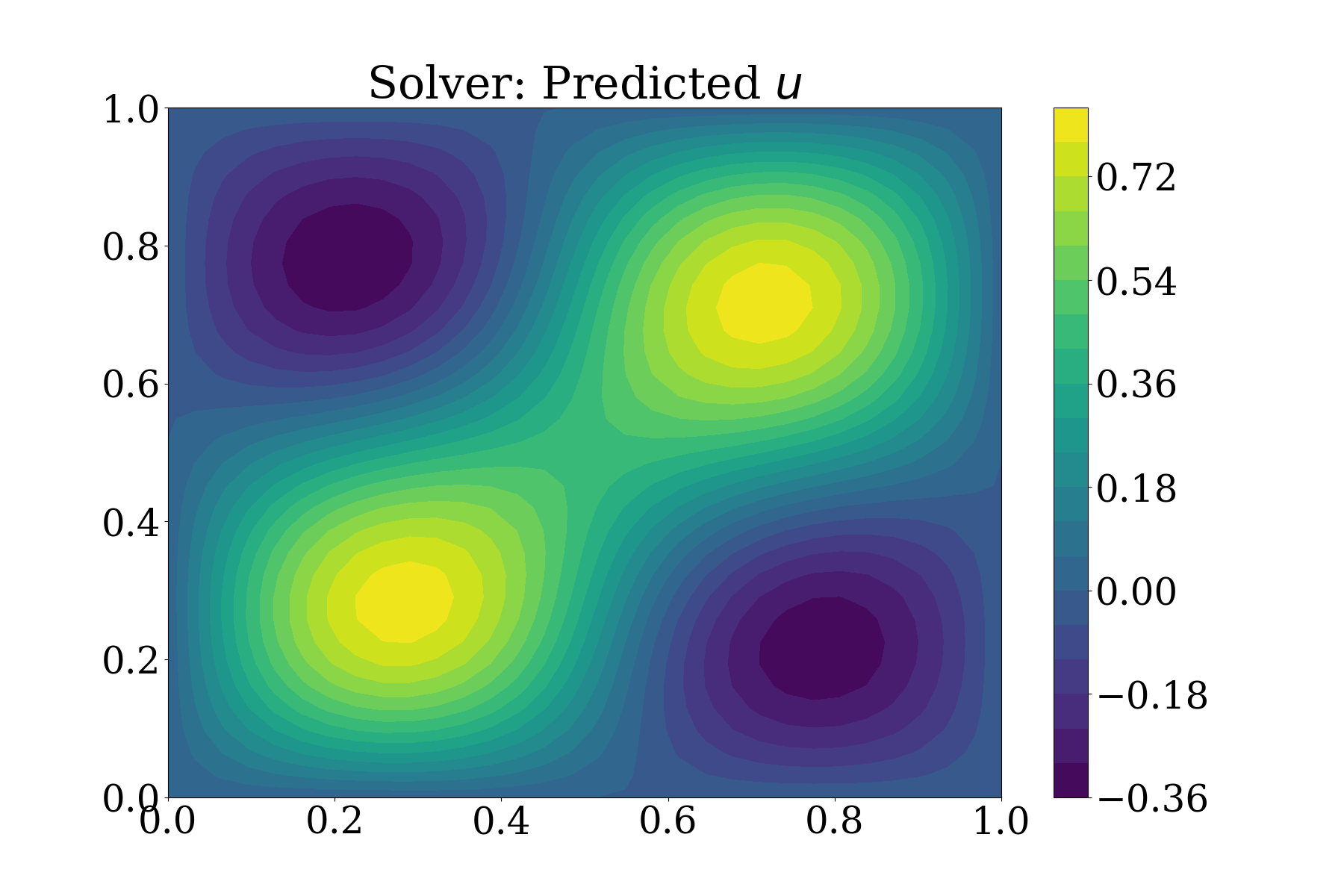}
        \subcaption{Solver: Predicted $u$}
    \end{subfigure}
    \begin{subfigure}[t]{0.32\textwidth}
        \centering
        \includegraphics[width=\textwidth]{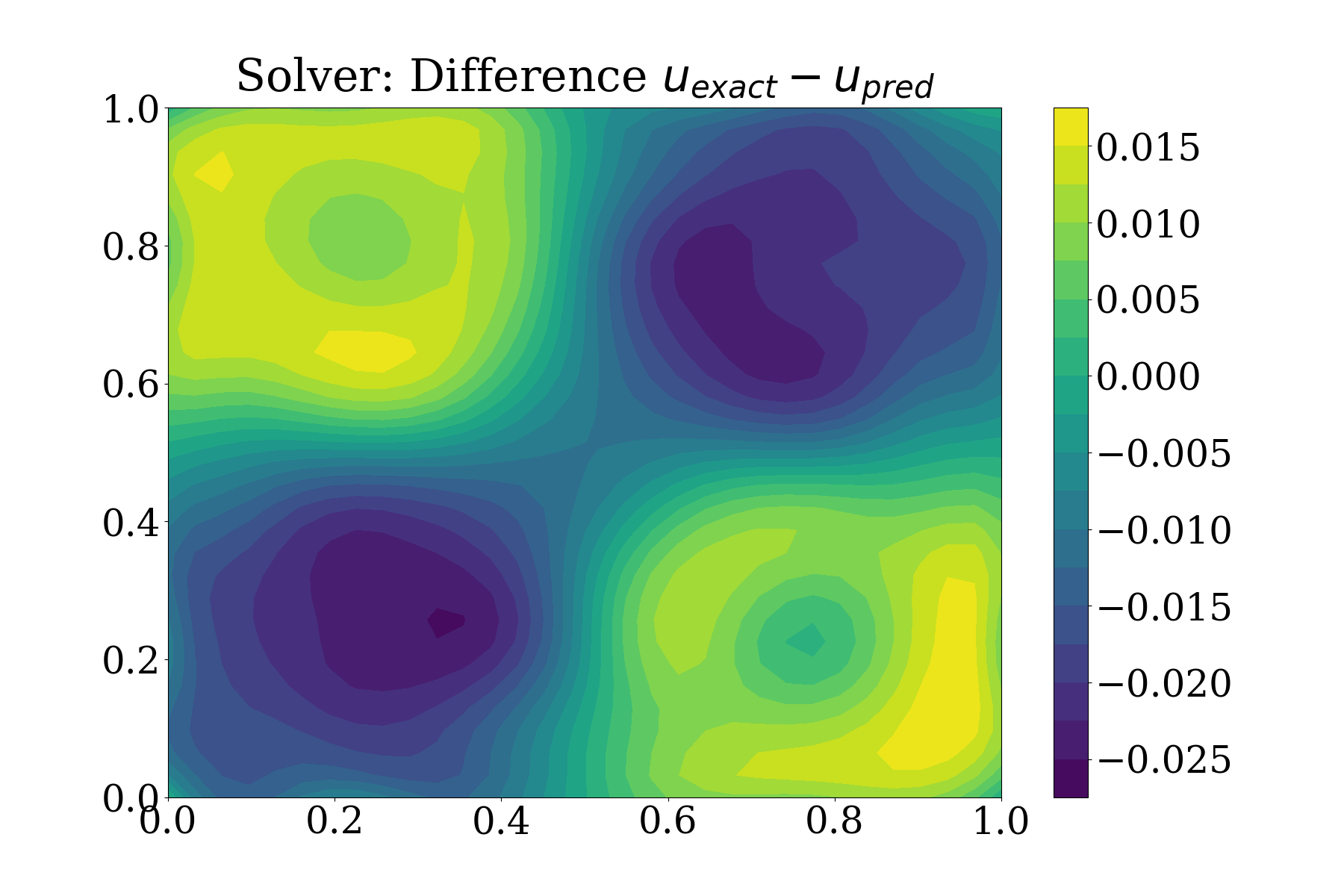}
        \subcaption{Solver: Difference $u$}
    \end{subfigure}
    
    \begin{subfigure}[t]{0.32\textwidth}
        \centering
        \includegraphics[width=\textwidth]{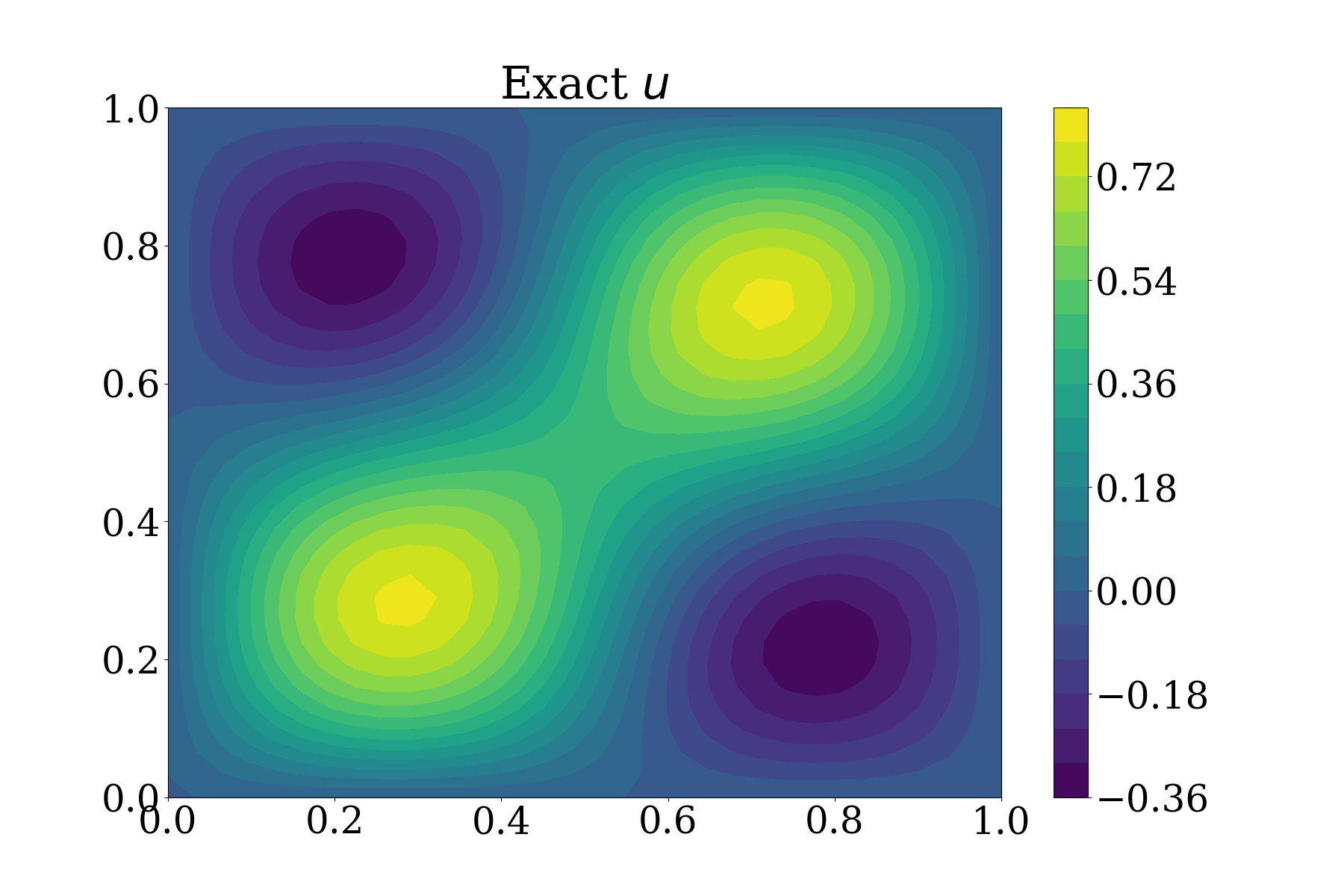}
        \subcaption{Exact $u$}
    \end{subfigure}
    \begin{subfigure}[t]{0.32\textwidth}
        \centering
        \includegraphics[width=\textwidth]{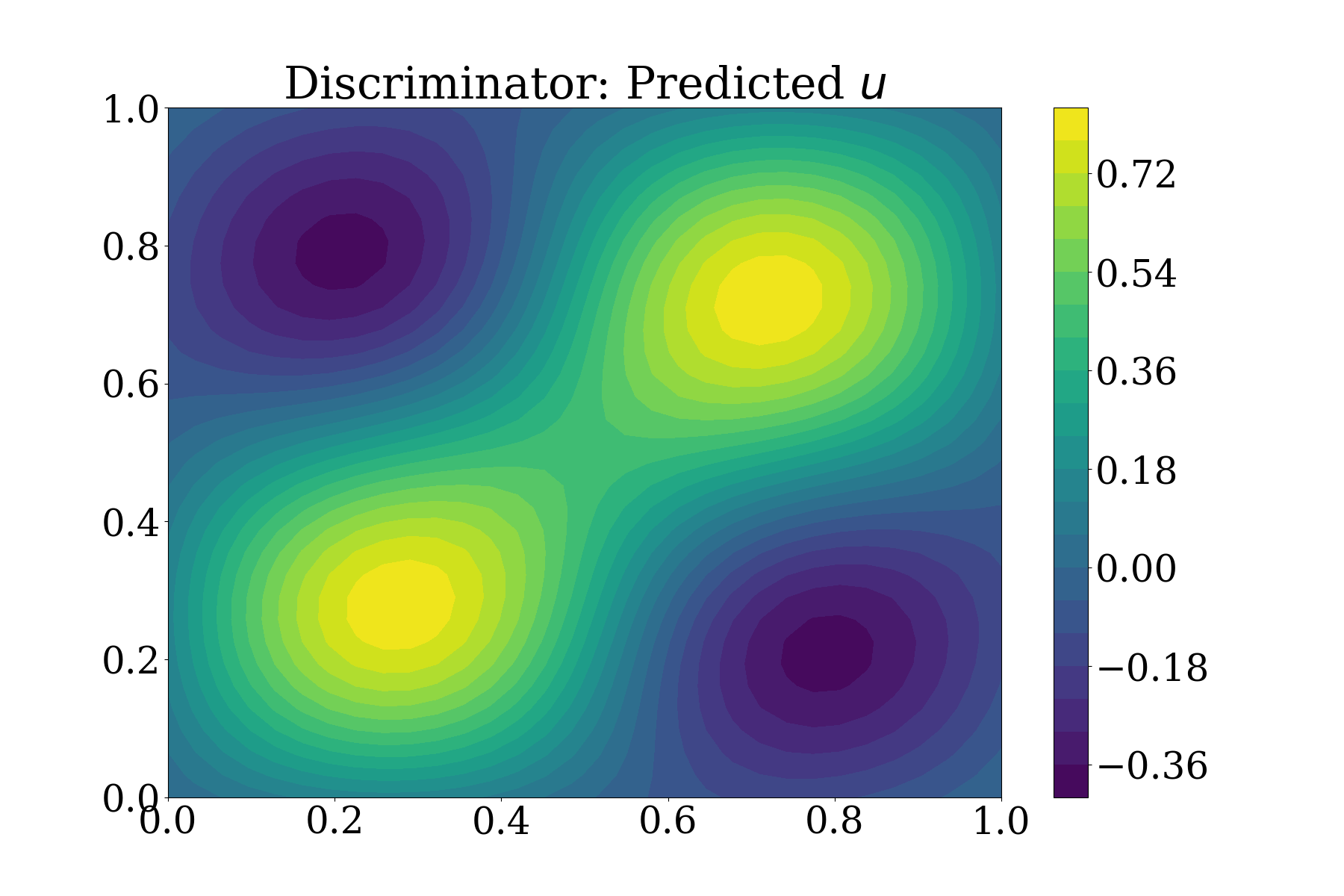}
        \subcaption{Discriminator: Predicted $u$}
    \end{subfigure}
    \begin{subfigure}[t]{0.32\textwidth}
        \centering
        \includegraphics[width=\textwidth]{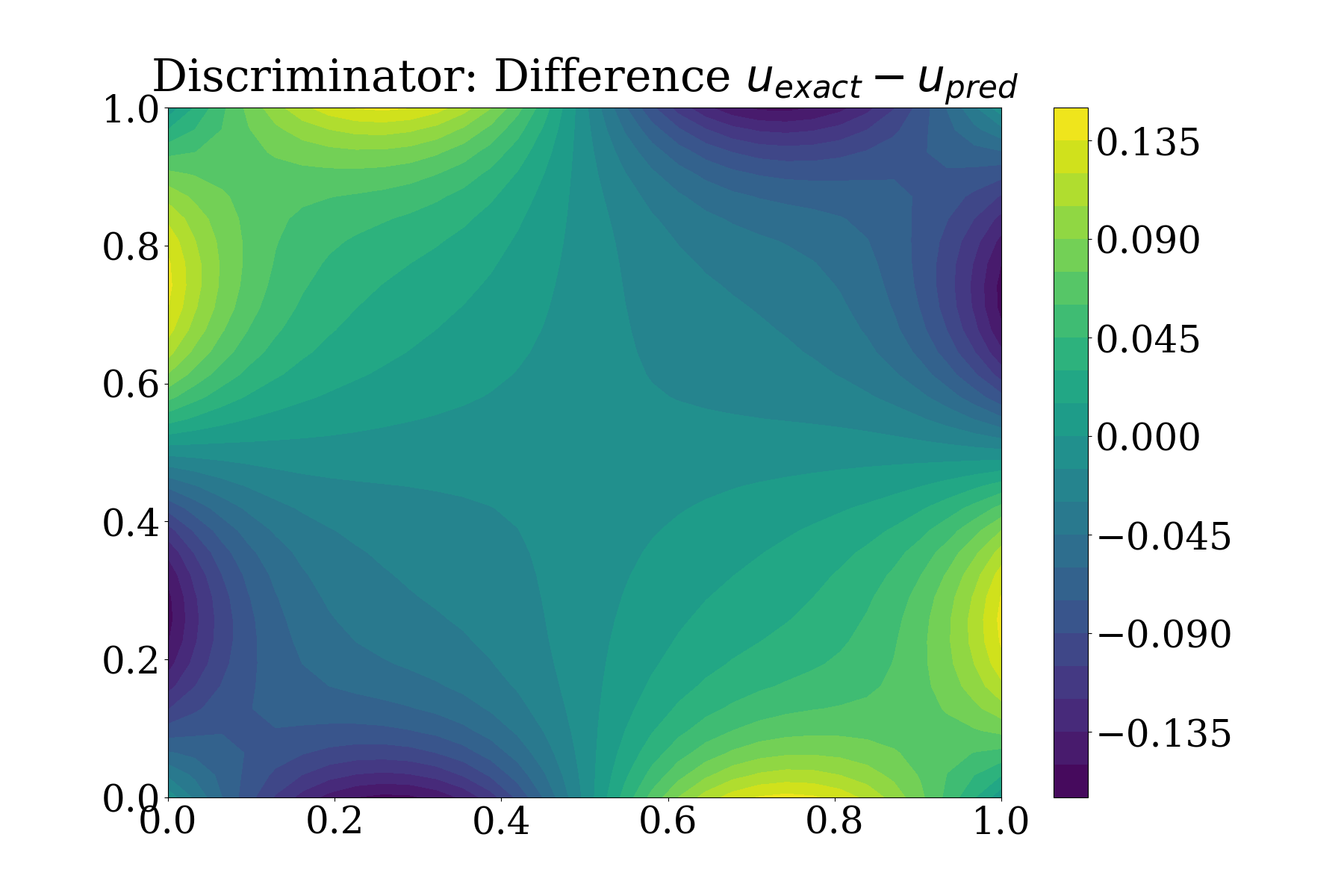}
        \subcaption{Discriminator: Difference $u$}
    \end{subfigure}

    \begin{subfigure}[t]{0.4\textwidth}
        \centering
        \includegraphics[width=\textwidth]{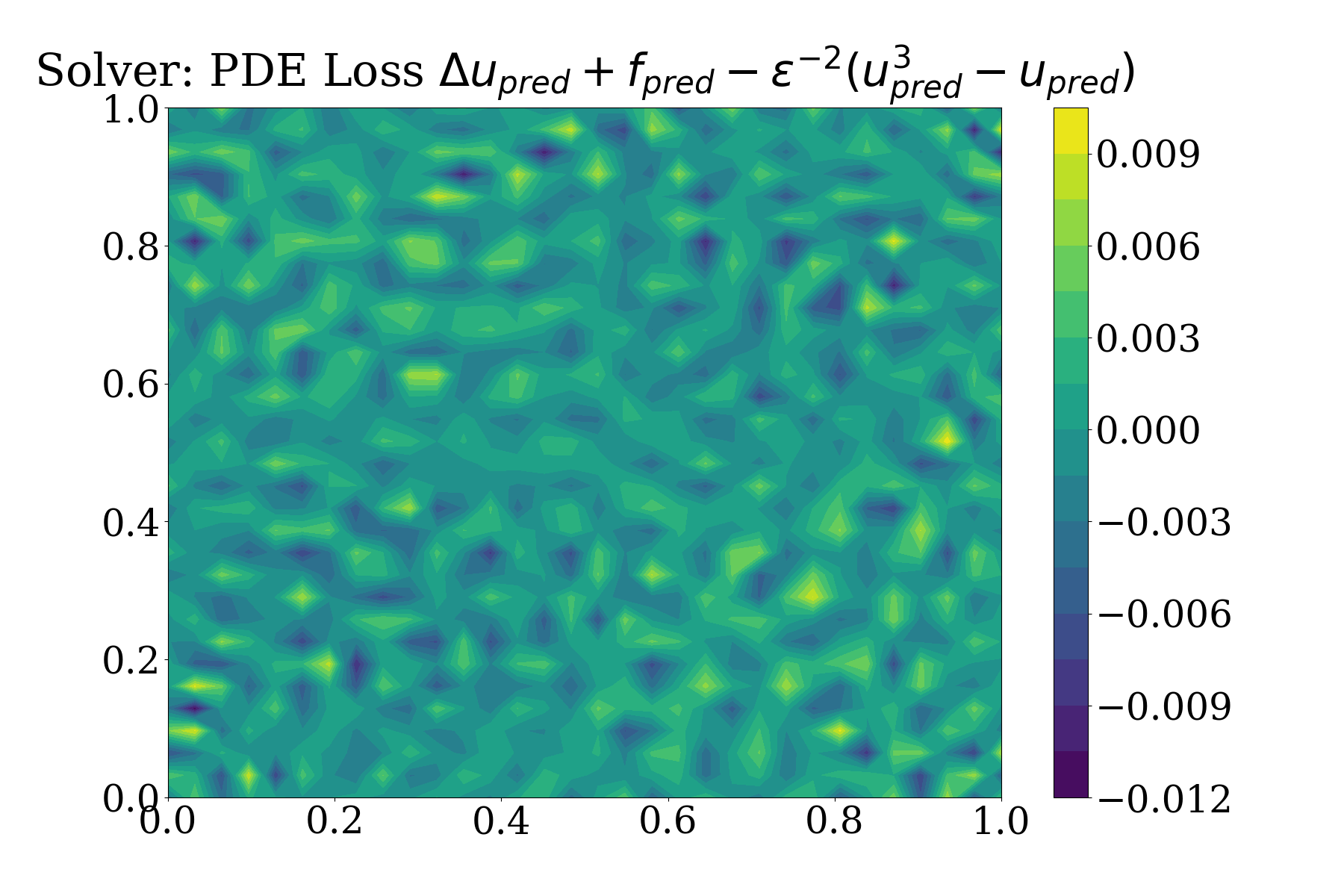}
        \subcaption{Solver: PDE Loss $\Delta u_{pred} + f_{pred}-\frac{1}{\epsilon^2}(u_{pred}^3-u_{pred})$}
    \end{subfigure}
    \begin{subfigure}[t]{0.4\textwidth}
        \centering
        \includegraphics[width=\textwidth]{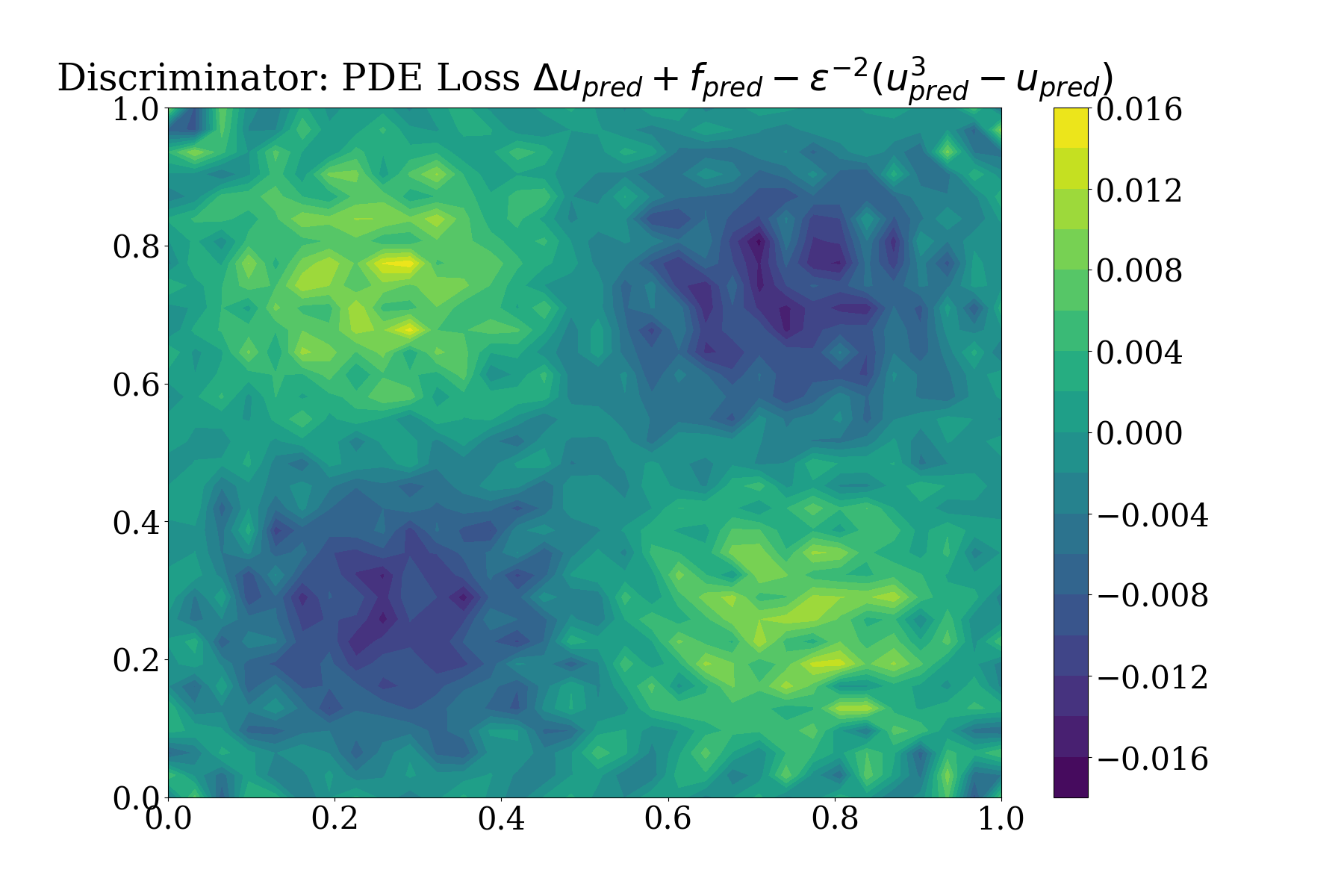}
        \subcaption{Discriminator: PDE Loss $\Delta u_{pred} + f_{pred}-\frac{1}{\epsilon^2}(u_{pred}^3-u_{pred})$}
    \end{subfigure}

    \begin{subfigure}[t]{0.5\textwidth}
        \centering
        \includegraphics[width=\textwidth]{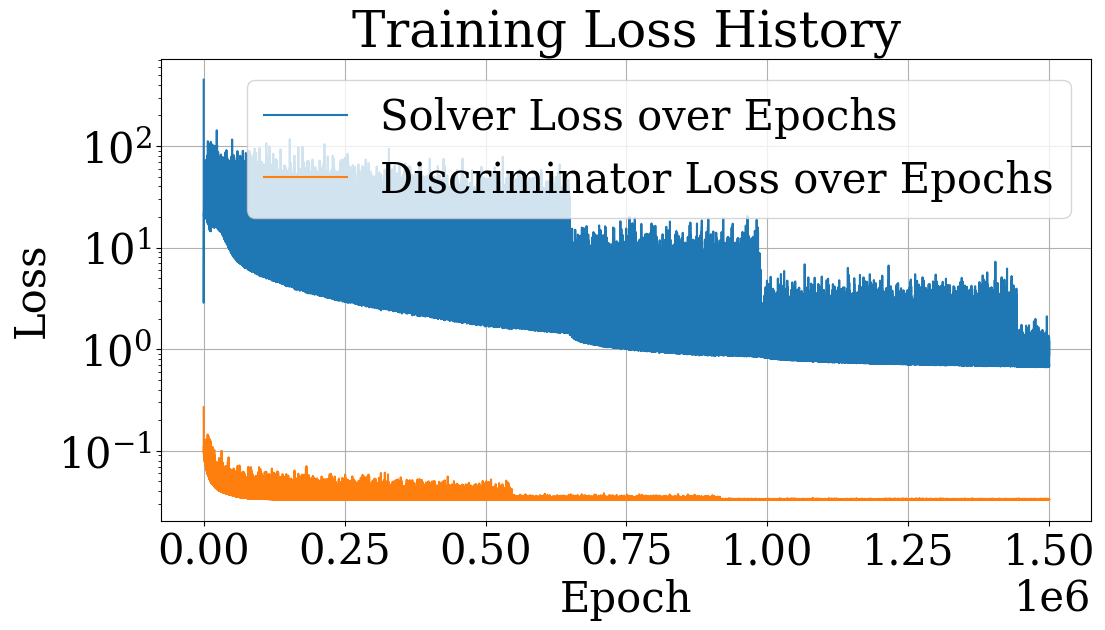}
        \subcaption{Loss over epochs}
    \end{subfigure}
    \caption{Results for minimizing the loss function \eqref{eq:2D_Allen_Cahn_penalty_loss_nn_discriminator} and \eqref{eq:2D_Allen_Cahn_penalty_loss_nn_solver} with $\lambda_p^d=\lambda_b^d = 0.2$, $\lambda_p^s=\lambda_b^s = 1000$,  $\omega = 20000$ and $\rho = 0.0001$. Subfigures (a) and (b) show the exact and numerical solutions for $u$ for the solver, respectively. Subfigure (c) shows the difference between exact and predicted $u$ for the solver. Subfigures (d) and (e) show the exact and numerical solutions for $u$ for the discriminator, respectively. Subfigure (f) shows the difference between exact and predicted $u$ for the discriminator. Subfigure (g) presents the PDE loss for the solver, and subfigure (h) presents the PDE loss for the discriminator. Subfigure (i) presents the loss over epochs.}
    \label{fig:2D_Allen_Cahn_adversarial}
\end{figure}

The penalty parameters selected for this example are as follows: $\lambda_p^s = \lambda_b^s = 1000$, $\lambda_p^d = \lambda_b^d = 0.2$, and $\omega = 20000$. It is important to note that we opted for a relatively large value for $\omega$ due to the complexity of the associated nonlinear problem. A small $\omega$ would not significantly accelerate the training process or help achieve convergence, which is the main advantage over using a large penalty parameter alone. Although we have established that $\omega$ should have an upper bound, as demonstrated in Proposition \ref{prop:omega_upper_bound}, this example illustrates that using a relatively large value for $\omega$ in some instances is practical. The results of this experiment are shown in Figure \ref{fig:2D_Allen_Cahn_adversarial}.

We observe that both the solver and discriminator networks converge to solutions close to the exact solution, with the solver network being notably more accurate than the discriminator network. This is evident in Subfigures (c) and (f), which explicitly show the differences between the solver's and discriminator's predictions compared to the exact solution, respectively. The inaccuracy of the discriminator network stems from the fact that its penalty parameter, set at $\lambda_p^d = \lambda_b^d = 0.2$, is too small to ensure a solution that closely approximates the exact solution of the constrained optimization problem. On the other hand, since $\omega$ is chosen to be large in this example, the expected advantage of the solver network adhering better to the constraint compared to the discriminator network becomes less pronounced, as observed in Subfigures (g) and (h). However, even so, we can still observe that the PDE loss for the solver's network is slightly smaller than that of the discriminator's network, aligning with theoretical expectations.

\begin{figure}
    \centering
    
    \begin{subfigure}[t]{0.32\textwidth}
        \centering
        \includegraphics[width=\textwidth]{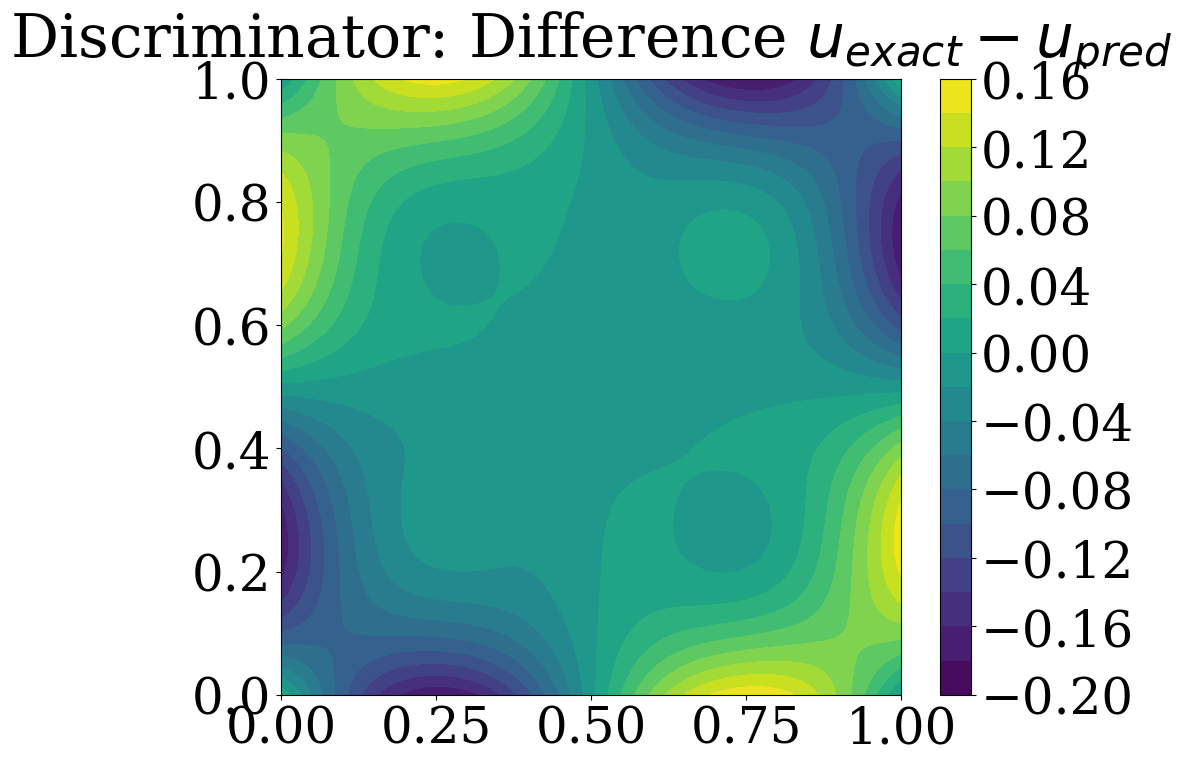}
        \subcaption{Difference at epoch 100000}
    \end{subfigure}
    \begin{subfigure}[t]{0.32\textwidth}
        \centering
\includegraphics[width=\textwidth]{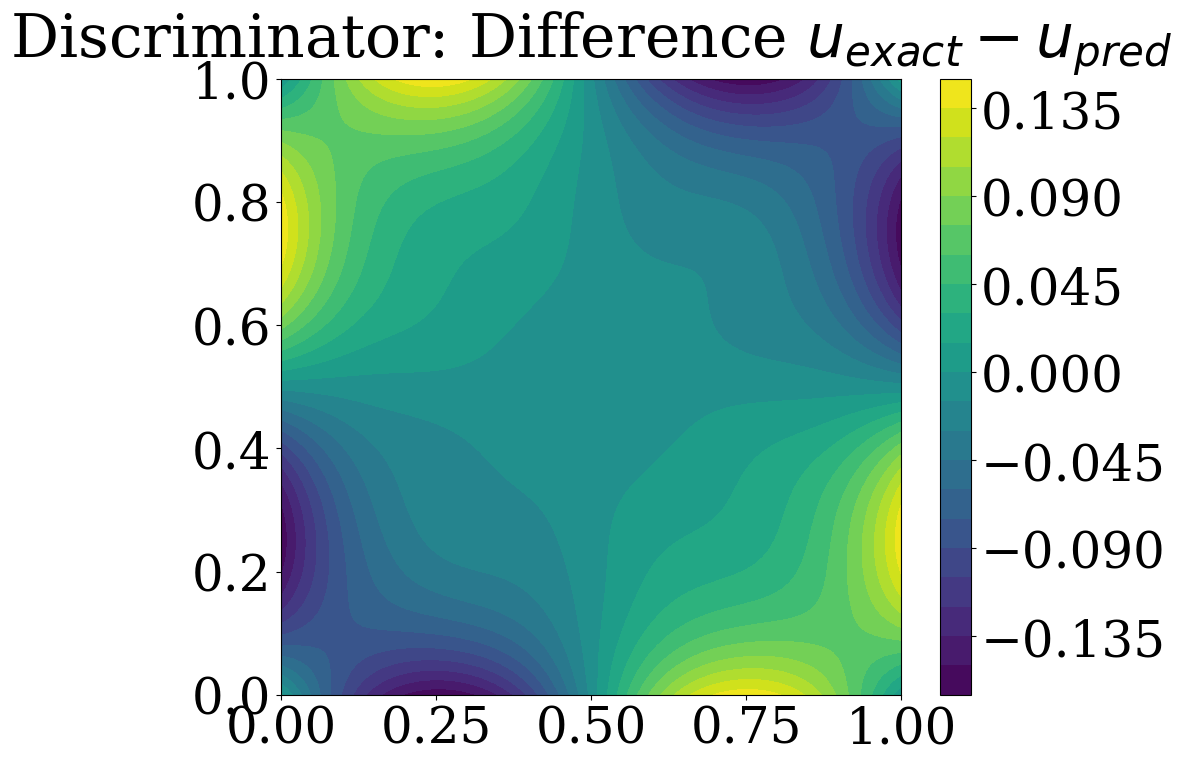}
         \subcaption{Difference at epoch 300000}
    \end{subfigure}
   \begin{subfigure}[t]{0.32\textwidth}
        \centering
\includegraphics[width=\textwidth]{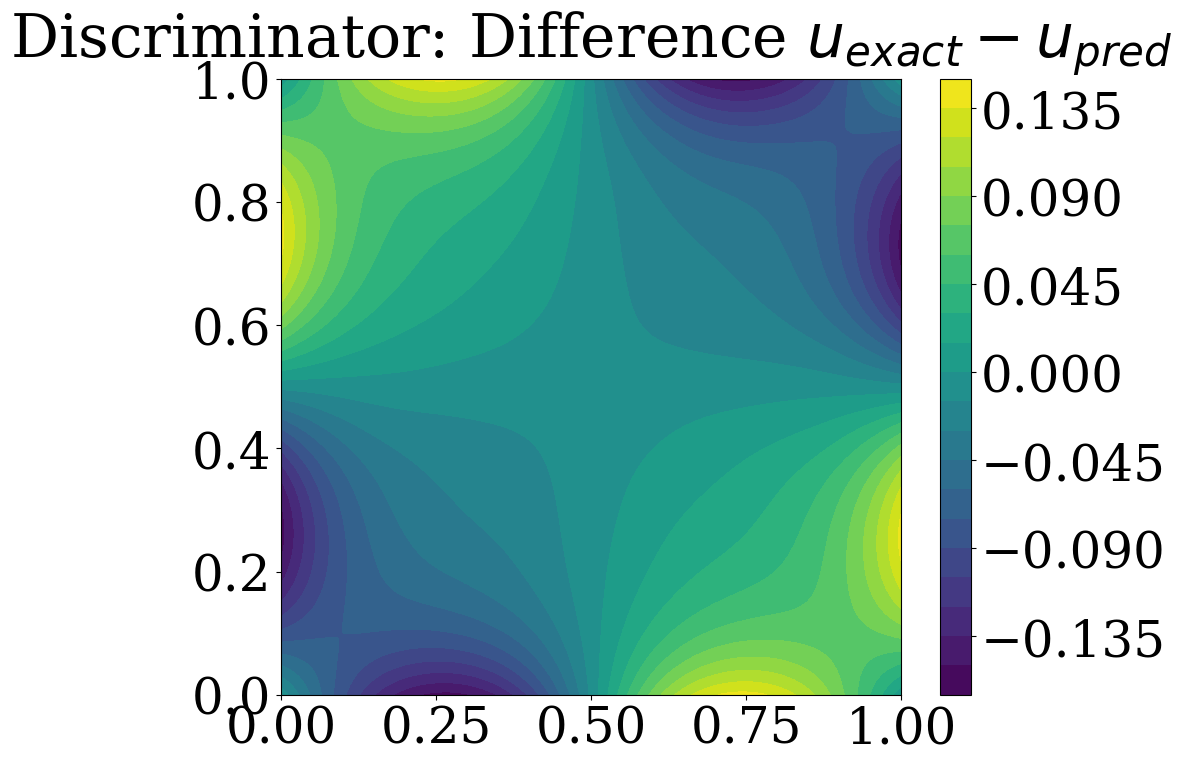}
         \subcaption{Difference at epoch 1500000}
    \end{subfigure}

    \caption{The difference between the exact solution and the prediction by discriminator network}.
    \label{fig:2D_Allen_Cahn_discriminator_over_epochs}
\end{figure}

Another important observation is that the discriminator network converges quickly, as expected, due to its small penalty parameter. As shown in Figure \ref{fig:2D_Allen_Cahn_discriminator_over_epochs}, the discriminator network finds a solution with similar errors to the final result after just 300000 epochs. As we can see, Subfigures (b) and (c) are very similar, indicating that the subsequent 1200000 epochs yield minimal improvement for the discriminator network. Additionally, the discriminator's prediction shows larger inaccuracies on the boundary compared to the interior points. This phenomenon may be related to the limited number of training points assigned to each boundary, with only 8 points on each side.

\begin{figure}
    \centering
    
    \begin{subfigure}[t]{0.30\textwidth}
        \centering
        \includegraphics[width=\textwidth]{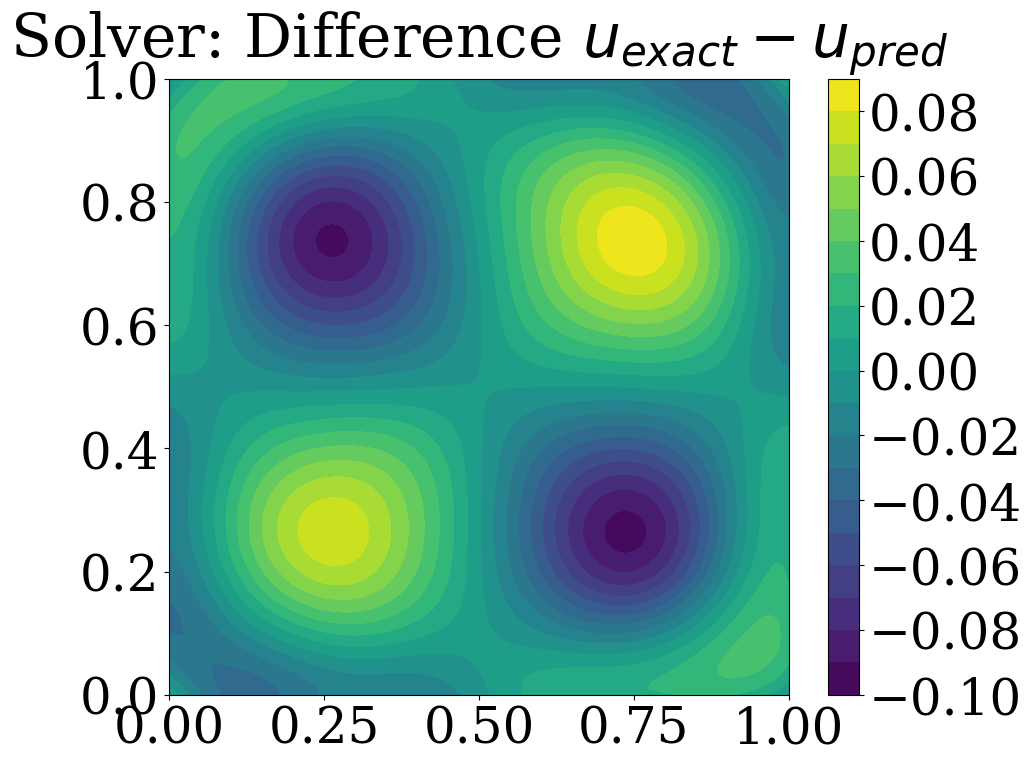}
        \subcaption{Difference at epoch 500000}
    \end{subfigure}
    \begin{subfigure}[t]{0.30\textwidth}
        \centering
\includegraphics[width=\textwidth]{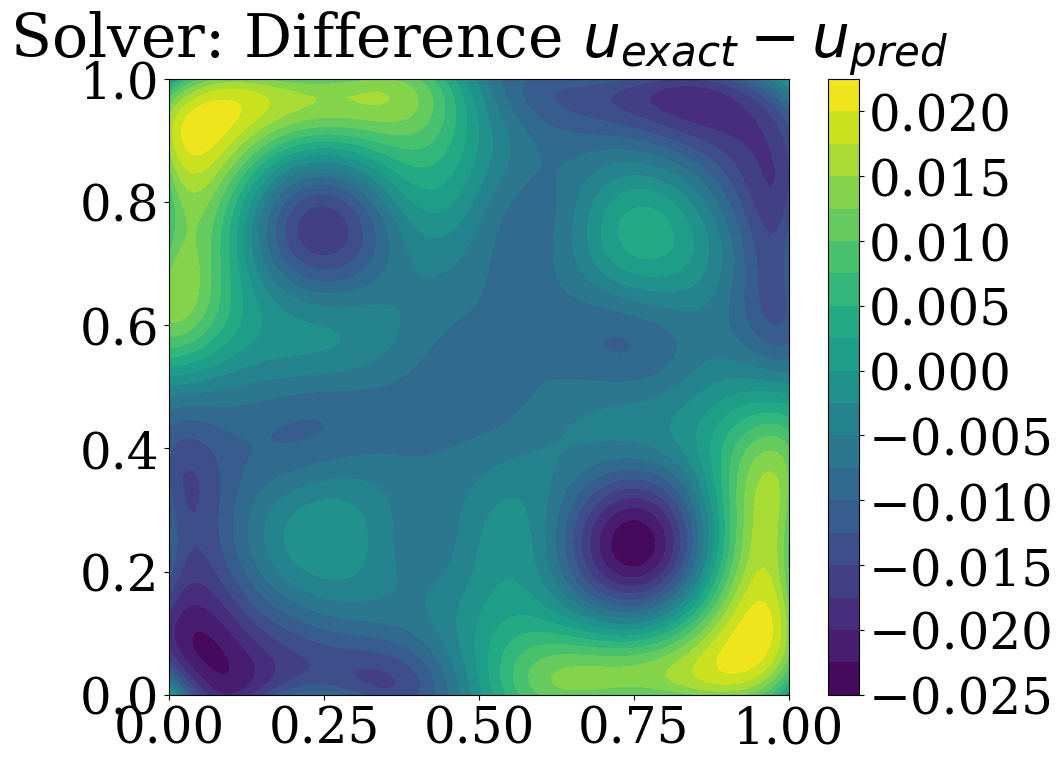}
         \subcaption{Difference at epoch 1000000}
    \end{subfigure}
   \begin{subfigure}[t]{0.30\textwidth}
        \centering
\includegraphics[width=\textwidth]{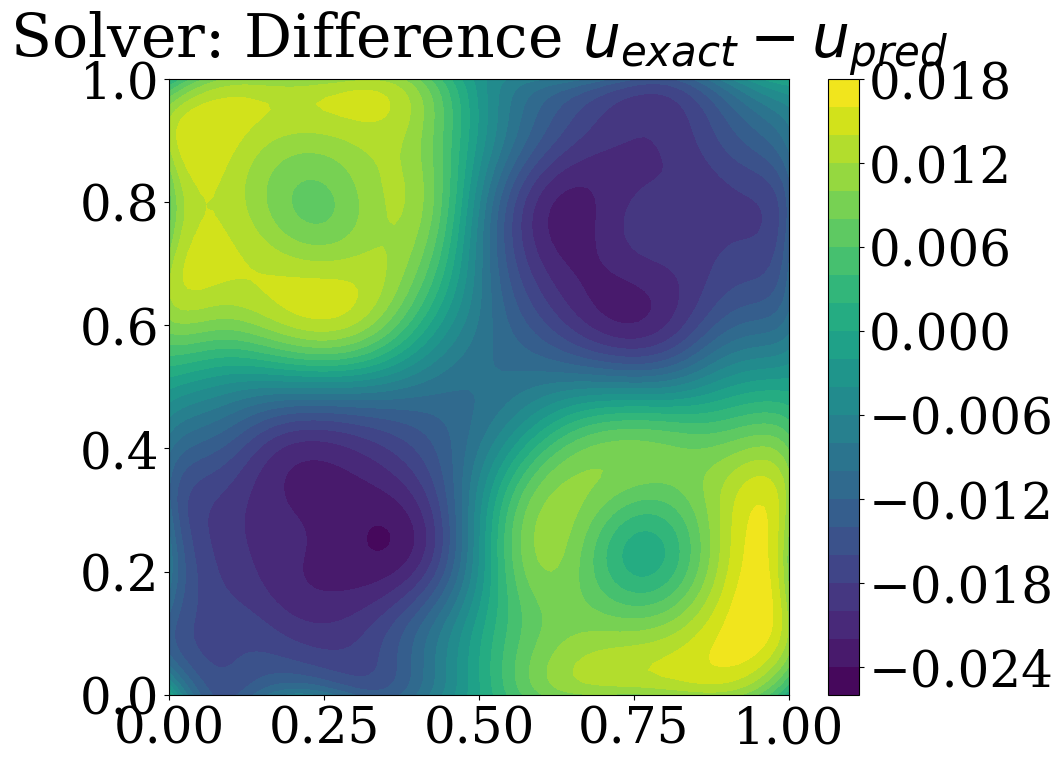}
         \subcaption{Difference at epoch 1500000}
    \end{subfigure}

    \caption{The difference between the exact solution and the prediction by solver network}.
    \label{fig:2D_Allen_Cahn_solver_over_epochs}
\end{figure}

In contrast, the solver network converges much more slowly, as evident in Figure \ref{fig:2D_Allen_Cahn_solver_over_epochs}. The solver network continuously reduces the error between its prediction and the exact solution, ultimately resulting in a more accurate solution than the discriminator network. However, it's important to note that while the solver network requires more time to converge than the discriminator, it is still far more efficient than simply using a traditional penalty formulation with a large penalty parameter, which will still be far from the actual solution after 1.5 million epochs. From Figure \ref{fig:2D_Allen_Cahn_solver_over_epochs}, we observe that after 500000 epochs, the prediction still deviates from the exact solution by a certain margin, but the result after 1000000 epochs is already quite good. The maximum error shown in Subfigure (b) is just slightly larger than in Subfigure (c) while the latter one took an additional 500000 epochs to achieve. This observation suggests a practical tradeoff between efficiency and accuracy that should be considered in real-world applications.

\section*{Acknowledgment}
\subsection*{Declaration of generative AI and AI-assisted technologies in the writing process}
During the preparation of this work, the author(s) used ChatGPT in order to improve the writing of this paper. After using this tool, the author(s) reviewed and edited the content as needed and take(s) full responsibility for the content of the publication.

\bibliographystyle{abbrv}
\bibliography{reference}
\end{document}